\documentclass[11pt, reqno]{amsart}
\pdfoutput=1
\usepackage{amssymb,amsmath,amsthm,amsfonts,latexsym}
\usepackage{epsfig,esint}
\usepackage{mathrsfs}
\usepackage{mathtools,textcomp}
\usepackage{graphicx}
\usepackage{epic,eepic,wrapfig,color, ifthen} 
\usepackage{url,cite,geometry}
\usepackage[colorlinks=true, citecolor=blue, filecolor=cyan, linkcolor=red, urlcolor=magenta]{hyperref}
\usepackage{pmboxdraw}
\usepackage{verbatim}
\usepackage[normalem]{ulem}
\usepackage{enumerate}
\usepackage[bottom]{footmisc}

\newcommand{\cal}[1]{\mathcal{#1}}

\def\1{{\bf 1}}
\def\nn{\nonumber}

\newcommand{\bk}{\color{black}}
\def\sA {{\cal A}} \def\sB {{\cal B}} \def\sC {{\cal C}}
\def\sD {{\cal D}}
\def\sE {{\cal E}} 
\def\sF {{\cal F}}
\def\sG {{\cal G}}  \def\sI {{\cal I}}
 
\def\sL {{\cal L}}

\def\sN {{\cal N}}

\def\sT {{\cal T}} 
\def\sV {{\cal V}}

\def\R {{\mathbb R}} 
 
\def\E {{\mathbb E}}  \def \P{{\mathbb P}}

\def \diam{{\text{\rm diam}}}

\def\p{{\Theta}}

\def \VD {$\mathrm{VD}$}
\def \RVD {${\mathrm {RVD}}$}

\def \KN {${\mathrm{KN}}(\phi)$}
\def \LN {${\mathrm{LN}}_\nu(\phi)$}
\def \HLN {${\mathrm{HLN}}_\nu(\phi)$}

\def \CS {${\mathrm{CS}}^\p(\phi)$}
\def \CSU {${\mathrm{CSU}}^\p(\phi)$}

\def \Gcap {${\mathrm{Gcap}}^\p(\phi)$}
\def \GU {${\mathrm{GU}}^\p(\phi)$}

\def \CSUq {${\mathrm{CSU}}(\phi)$}

\def \Gcapq {${\mathrm{Gcap}}(\phi)$}

\def\Cap{{\mathrm{Cap}}}

\def \FK {${\mathrm{FK}}(\phi)$}
\def \PI {${\mathrm{PI}}(\phi)$}

\def \Tail {${\mathrm{TJ}}^\p_{\le}(\phi)$}

\def \Ju {${\mathrm{J}}_\le^\p(\phi)$}
\def \Jl {${\mathrm{J}}_\ge(\phi)$}

\def \Tailq {${\mathrm{TJ}}_{\le}(\phi)$}

\def \Juq {${\mathrm{J}}_\le(\phi)$}

\def \mEl {${\mathrm{E}}_\ge(\phi)$}
\def \mEu {${\mathrm{E}}_\le(\phi)$}
\def \mE {${\mathrm{E}}(\phi)$}
\def \SP {${\mathrm{SP}}_{\le}(\phi)$}

\def \NDL {${\mathrm{NDL}}(\phi)$}
\def \sNDL {${\mathrm{sNDL}}(\phi)$}

\def \EHR {${\mathrm{EHR}}(\phi)$}
\def \PHR {${\mathrm{PHR}}(\phi)$}

\def \WEHI {${\mathrm{WEHI}}_\p(\phi)$}
\def \WEHIp {${\mathrm{WEHI}}^+_\p(\phi)$}
\def \NT{{\mathrm{Tail \hspace{0.3mm}}}_\p}

\def \Exi {${\mathrm{Exi}}(\alpha)$}

\def\la{{\langle}}
\def\ra{{\rangle}}

\numberwithin{equation}{section}
\def\qed{{\hfill $\Box$ \bigskip}}
\def\eps{\varepsilon}
\def\wh{\widehat}
\def\wt{\widetilde}
\def\pf{\noindent{\bf Proof. }}

\def\vp{{\varphi}}

\DeclareMathOperator*{\esssup}{ess\,sup}
\DeclareMathOperator*{\essinf}{ess\,inf}
\DeclareMathOperator*{\essosc}{ess\,osc}

\theoremstyle{plain}
\newtheorem{thm}{Theorem}[section]
\newtheorem{lem}[thm]{Lemma}
\newtheorem{cor}[thm]{Corollary}
\newtheorem{remark}[thm]{Remark}
\newtheorem{prop}[thm]{Proposition}
\newtheorem{defn}[thm]{Definition}

\newtheorem{example}[thm]{Example}

\theoremstyle{definition}
\newtheorem*{eg*}{Example}
\newtheorem*{egs*}{Examples}
\newtheorem*{def*}{Definition}
\theoremstyle{remark}

\begin{document}
	\title{Stability of H\"older regularity and 
	weighted functional inequalities
}

	\author{Soobin Cho \qquad Panki Kim}\thanks{This research is  supported by the National Research Foundation of Korea(NRF) grant funded by the Korea government(MSIT) (No. RS-2023-00270314).
	}

	\address[Cho]{Department of Mathematics, University of Illinois Urbana-Champaign, Urbana, IL 61801, USA}
	\curraddr{}
	\email{soobinc@illinois.edu}

	\address[Kim]{Department of Mathematical Sciences and Research Institute of Mathematics,
		Seoul National University,
		Seoul 08826, Republic of Korea}
	
	\curraddr{}
	\email{pkim@snu.ac.kr}

\begin{abstract}
We study symmetric  Dirichlet forms on metric measure spaces, which may possess  both strongly local and pure-jump parts.  We introduce a new formulation of a tail condition for jump measures and weighted functional inequalities. Our framework accommodates Dirichlet forms with singular jump measures and those  associated with trace processes of mixed-type stable processes. Using these new weighted functional inequalities, we establish stable, equivalent characterizations of  H\"older regularity for caloric and harmonic functions. As an application of our main result, we prove the H\"older continuity of caloric functions for a large class of symmetric Markov processes  exhibiting boundary blow-up behavior, among other results.

	\medskip
	
	\medskip
	
\end{abstract}

	\maketitle
	
	\allowdisplaybreaks
	
{
	\hypersetup{linkcolor=blue}
	\tableofcontents
}

	\section{Introduction}\label{s:intro}

In recent years, significant foundational work has been conducted on stable and equivalent characterizations across various contexts using Dirichlet forms, including heat kernel estimates, elliptic and parabolic Harnack inequalities, and the regularity of caloric and harmonic functions.
See \cite{AB15, BKKL2, BBK06, BCM, BGK, BM, CKW-elp, CKW-adv, CKW-jems, CKW-memo, Gr1, GH14, GH, GHH18, GHH, GHH2, GHL, GT12,Ki, Sa} and references therein.

  In this paper, we establish stable, equivalent characterizations of H\"older continuity for caloric and harmonic functions within a broad class of symmetric regular Dirichlet forms, which may include both strongly local and pure-jump parts. The Dirichlet forms we consider can be highly singular, with jump measures whose tails may diverge to infinity at the boundary of the state space. 
  This phenomenon arises in trace processes of jump processes, which serve as a special case of our more general framework. To highlight the novelty and scope of our  results, we recall the construction of trace processes and the behavior of the tails of their jump measures  in Euclidean space, as studied in \cite{BGPR, KSV23, KSV}.
		\begin{example}\label{e:tpc}
{\rm	Let $j$ be a positive,  non-increasing function on $(0, \infty)$, and let $Z=(Z_t, \P_x)$ be a unimodal  L\'evy process  in $\R^d$ with  L\'evy measure $j(|x|)dx$. 
For  a smooth open set $D\subset \R^d$  
 such that $U:=\overline D^{\,c}
 $ is non-empty, 
define
$$
A_t:=\int_0^t \1_{\{Z_s\in \overline D\}}\, ds.
$$
Let $\tau_t:=\inf\{s>0: \, A_s>t\}$ 
be the right-continuous inverse of $A_t$. The process 
$X=(X_t)_{t\ge 0}$ 
defined by $X_t:=Z_{\tau_t}$ is a Hunt process with state space $\overline D$, called the trace process of $Z$ on $\overline D$. Let $X^D$ be the  process  $X$ killed upon leaving $D$, or equivalently, upon hitting the boundary $\partial D$. This process  is called the \textit{path-censored process}. For comparisons among the path-censored process $X^D$, the L\'evy process conditioned to stay positive, and the censored stable process, we refer to \cite[Remark 3.3]{KPW14}.

One can describe the behavior of the process $X^D$ as follows: When $Z$ jumps out of $D$ from $x=Z_{\tau_D-}\in D$, it is known that  it almost surely lands at $z\in U$.   The distribution of the returning position $y\in D$ of $Z$ to $D$ is given by the Poisson kernel of $Z$ with respect to $U$: 
$$P_U(z, y)=\int_U G_U(z,w)j(|w-y|)\,  dw, \quad z\in D^c, \, y \in D.$$
Here  $G_U(z,w)$ represents the Green function of the process $Z$ killed upon exiting $U$. 
This implies that when $Z$ jumps out of $D$ from $x$, the process is continued  by resurrecting it in $D$  according to the kernel $q(x,y)dy$, where
$$
q(x,y)=\int_U \int_U j(|x-z|)G_U(z,w)j(|w-y|)\, dz\, 
dw, \quad x,y\in D.
$$

The Dirichlet form $(\sE,\sF)$ associated with $X$ is given by
\begin{align*}
\sE(f,f)=\frac12\int_{D \times D} (f(x)-f(y))(g(x)-g(y)) J(x,y)dxdy,
\end{align*} 
where $J(x,y):=j(x-y)+q(x,y)$, and $\sF$ is the closure of $C_c^2(\overline D)$ under the norm
 $\sE_1(f,f) :=\sE(f,f) + \lVert f \rVert_{L^2(D)}^2$. 
When $Z$ is an isotropic $\alpha$-stable process so that  $j(r)=c  r^{-d-\alpha}$, and $D$ is either a half-space or an exterior smooth open set, it follows from  \cite[Theorems 6.1 and 2.6]{BGPR}  that
$$\int_{\{y\in D:|x-y|>r\}}J(x,y)dy \asymp   r^{-\alpha}\left(1 +\frac{r}{\text{dist}(x, \partial D))}  \right)^{\alpha/2} \quad \text{for all $x\in D$ and $r \in (0, \diam(D^c))$.}
$$
Here for non-negative functions $f$ and $g$, the notation $f\asymp g$  means that $f/g$ is bounded between two positive constants, and $\diam(D^c)$ denotes the diameter of $D^c$.  In particular, we have   $\lim_{y\to \partial D}J(x,y)=\infty$ for all $x\in D$, that is, the jump kernel blows up at the boundary $\partial D$.
}
\end{example}

Another important example of a	Dirichlet form  with a  blow-up jump kernel arises in the study of non-local Neumann problems, introduced in \cite{DRV17}.

While Dirichlet forms with jump parts have been extensively studied, results concerning forms with blow-up jump kernels remain limited. One reason is that standard frameworks for non-local Dirichlet forms often assume uniform upper bounds on the tails of jump measures.    Recently, Audrito, Felipe-Navarro, and Ros-Oton \cite{AFR23} established optimal H\"older regularity up to the boundary for harmonic functions associated with the non-local Neumann problem. Additionally,   Song, Vondra\v{c}ek, and the second-named author \cite{KSV23} studied interior potential-theoretic properties of purely discontinuous Markov processes in $D\subset \R^d$ with potentially degenerate jump kernels at the boundary, including the process in Example \ref{e:tpc}. However, except for certain special cases,   a general regularity theory up to the boundary for Dirichlet forms on metric measure spaces with blow-up jump kernels has remained unknown.

The goal of this paper is to develop a general framework for  establishing regularity results for both caloric and harmonic functions associated with Dirichlet forms whose jump measures do not satisfy uniform upper bounds on their tails.  To achieve this, we introduce the notion of an admissible weight function, and  weighted versions of a tail condition for jump measures, and function inequalities such as
a cut-off Sobolev inequality, and a generalized capacity condition. These conditions are commonly used  in the analysis of Dirichlet forms on general metric measure spaces.

Our approach combines probabilistic and analytic techniques, using functional inequalities and their equivalences to study Dirichlet forms.  We refer to \cite{CKW-elp, CKW-adv, CKW-jems, CKW-memo, GHH18, GHH, GHH2} for unweighted versions of these   inequalities and stability results for Dirichlet forms satisfying a uniform upper tail condition (with respect to a suitable metric). In particular, we establish weighted versions of the Caccioppoli inequality and   
$L^2$-mean value inequalities. Our results apply to forms associated with trace processes, non-local Neumann problems, and other settings, extending recent advances in this area.

The paper is organized as follows. In Section \ref{s:setupP}, we introduce our setup, including a general space-inhomogeneous scale function $\phi$ and the notion of an admissible weight function  $\p$ associated with $\phi$, and state our main results. Section \ref{s:preliminary} presents preliminary results, including relations among the Faber-Krahn inequality, various Nash-type inequalities, and upper bounds for (Dirichlet) heat kernels. Section \ref{s:ex_weight} provides explicit examples of admissible weight functions, and Section \ref{s:CS} presents examples of Dirichlet forms for which the weighted cutoff Sobolev inequality holds. In Section \ref{s:CTCS}, we discuss consequences of the tail condition with weight and the weighted cutoff Sobolev inequality. Section \ref{s:L2} establishes weighted versions of Caccioppoli and $L^2$-mean value inequalities. In Section \ref{s:proof-thm-1}, we prove the equivalence of weighted cutoff Sobolev inequalities, mean exit time estimates, and weighted versions of generalized capacity conditions. Section \ref{s:weHI} proves a weak elliptic Harnack inequality with weighted non-local tails (Definition \ref{d:WEHI}). In Section \ref{s:eHR}, we establish H\"older regularity for harmonic functions. Section \ref{s:NDL} deals with  near-diagonal lower bounds for the heat kernel. In Section \ref{s:PHR}, we prove H\"older continuity for caloric functions and complete the proof of our main result.  Section \ref{s:example} applies our main results to examples of symmetric regular Dirichlet forms with blow-up jump measures, establishing near diagonal lower heat kernel estimates and  H\"older continuity results for caloric and harmonic functions. Finally, Appendix \ref{s:Appendix} contains proofs of preliminary results from Section \ref{s:preliminary}.

\medskip

\textit{Notations:} 
Throughout the paper, 
values of lower case letters with subscripts $c_i$, $i=0,1,2,...$,  are fixed in each statement and proof, and the labeling of these constants starts anew in each proof. We use the symbol ``$:=$'' to denote a definition, 
which is read as ``is defined to be.''  
Recall that the notation $f\asymp g$ means that  $c^{-1}g\le f \le c g$ for some $c\ge 1$. 
We write  $a\wedge b:=\min\{a,b\}$ and  $a\vee b:=\max\{a,b\}$.  
 For $U,V\subset M$, the notation $U \Subset V$ means that  the closure of $U$ is contained in $V$.  
 For   $D \subset M$, $D^c$ denotes the complement of $D$ in $M$, and  
 $\delta_D(x):=\inf\{d(x,y): y \in D^c\}.$    For $p \in [1,\infty]$, we denote by $\lVert \cdot  \rVert_p$ the $L^p$-norm in $L^p(M;\mu)$, and   by  $\la \cdot, \cdot \ra$  the inner product  in  $L^2(M;\mu)$.
  Let $\Sigma_M$ be the $\sigma$-algebra of Borel sets on $M$.    We refer to a function as measurable if it is Borel measurable. 
 For $E \in \Sigma_M$, let $\sB(E)$ be the family  of all measurable  functions on $E$, $\sB_b(E)$ the family of all bounded Borel functions on $E$,  $\sB_c(E)$ the family of all elements of $\sB_b(E)$ with compact supports in $E$, and   $C_c(E)$ the family  of all continuous functions on $E$ with compact supports in $E$.   For $f \in \sB(E)$, we let  supp$[f]$ be  the closure of $\{ x \in E: f(x)\not=0\}$ in $E$.

\section{Setup and main results}\label{s:setupP}

\subsection{Basic setup}
	
	Let $(M,d)$ be a locally compact separable metric space  and $\mu$ be a positive Radon measure on $M$ with full support. The triplet $(M,d,\mu)$ is called a \textit{metric measure space}. Throughout the
	paper, we assume that all balls in $(M,d)$ are relatively compact. 	Denote by  $B(x,r)$  the open ball centered at $x$ with radius $r$, and let $V(x,r):=\mu (B(x,r))$. From now on, we fix a constant $R_0 \in (0, \text{diam}(M)]$,  	where diam$(M)$ is the diameter of $M$.

	We  consider the following (local) volume doubling properties of the  metric measure space.
	\begin{defn}\label{d:VD}
		{\rm (i) We say that $(M,d,\mu)$ satisfies the \textit{volume doubling property}  \VD \ if there exist  constants $d_2,C_\mu>0$ such that
			\begin{align*}
				\frac{V(x,R)}{V(x,r)} \le  C_\mu \left( \frac{R}{r} \right)^{d_2} \quad \text{for all} \;\, x \in M \text{ and } 0<r \le R<2R_0.
			\end{align*}
			(ii) We say that $(M,d,\mu)$ satisfies  \textit{the reverse volume doubling property}   \RVD \ if there exist   constants $d_1,c_\mu>0$ such that 
			\begin{align*}
				\frac{V(x,R)}{V(x,r)} \ge  c_\mu \left( \frac{R}{r} \right)^{d_1} \quad \text{for all} \;\, x \in M \text{ and } 0<r \le R<R_0.
			\end{align*}

	} \end{defn}

	\smallskip
	
	Using  $V(y,R) \le V(x, R+d(x,y))$ for all $x,y \in M$ and $R>0$, we see that \VD \ is equivalent to the following: There exist constants $d_2,C_\mu>0$ such that
	\begin{align}\label{e:VD2}
		\frac{V(y,R)}{V(x,r)} \le C_\mu\bigg(\frac{R + d(x,y)}{r}\bigg)^{d_2}  \quad \text{for all} \;\,  x,y \in M \text{ and } 0<r\le R<2R_0-d(x,y).
	\end{align}

	\medskip

Let  $(\sE, \sF)$ be a  regular symmetric Dirichlet form   on $L^2(M;\mu)$. 
	 According to the Beurling-Deny formula (see \cite[Section 3.2]{FOT}),  $(\sE,\sF)$  can be decomposed into the \textit{strongly local part}, the \textit{pure-jump part}, and the \textit{killing part}. In this paper, we always assume that  $(\sE, \sF)$ has no killing part, namely,
	 \begin{align*}
	 	\sE(u,v)=\sE^{(L)}(f,g) + \sE^{(J)}(f,g), \quad f,g \in \sF,
	 \end{align*}
 where $\sE^{(L)}$ is the strongly local part of $(\sE,\sF)$ (see \cite[(3.2.2)--(3.2.3)]{FOT}) and $\sE^{(J)}$ is the pure-jump part associated with a positive symmetric Radon measure $J$ on $M \times M$:
	 \begin{align*}
	 	\sE^{(J)}(f,g)&= \int_{M \times M} (f(x)-f(y))(g(x)-g(y))J(dx,dy).
	 \end{align*}
We call $J$ the jump measure of $\sE^{(J)}$.  For $\lambda>0$, let $\sE_\lambda(u,u):=\sE(u,u)+\lambda \lVert u \rVert_2^2$. 

 The \textit{$L^2$-generator} $(\sL,\sD(\sL))$ of $(\sE,\sF)$ is defined as follows: The domain $\sD(\sL)$ is the collection of  all  $f \in \sF$ such that there is some (unique) $u \in L^2(M;\mu)$ so that $\sE(f,g)=\la u,g\ra$ for all $g \in \sF$, and  $\sL f := -u$.
  Let  $P_t:=e^{t\sL}$ be the associated semigroup for $(\sL,\sD(\sL))$ (or  $(\sE,\sF)$).

 By a general theory, there exists a $\mu$-symmetric Hunt process $X=\{X_t,t \ge 0;\P^x, x \in M \setminus \sN\}$, where $\sN \subset M$ is a properly exceptional set for $X$, which is associated  with the form $(\sE, \sF)$ in the sense that for 
 any $f \in \sB_b(M)$ and $t>0$, 
 $$\E^x[f(X_t)] = P_t f(x) \quad \text{for $\mu$-a.e. $x\in M$}.$$
We refer to \cite{FOT} for details. Let us fix   $\sN$ and denote $M_0=M \setminus \sN$.

A Borel measurable function $\phi:M \times (0,\infty)\to (0,\infty)$ is called a \textit{scale function} on $M$, if it satisfies the following three properties:
\begin{flalign}
\indent	\text{(1)	For each $x \in M$, the map $r \mapsto \phi(x,r)$ is strictly increasing and  continuous.}&&\label{e:phi-ic}
\end{flalign} 

\vspace{-0.07in}

(2)  There exist constants $\beta_2 \ge \beta_1>0$ and  $c_L, c_U>0$  such that
	\begin{align}\label{e:phi-scale}
c_L\left( \frac{R}{r} \right)^{\beta_1} \le 		\frac{\phi(x,R)}{\phi(x,r)} \le c_U\left( \frac{R}{r} \right)^{\beta_2} \quad \text{for all} \;\, x \in M \text{ and } 0<r\le R< \infty.
	\end{align}

(3)    There exists a constant $C\ge 1$ such that
\begin{align}\label{e:phi-comp}
\phi(y,r)\le C\phi(x,r)  \quad \text{for all} \;\, x,y \in M  \text{ with }  d(x,y) \le r.
\end{align}

The simplest example of a scale function is $\phi_\beta=r^\beta$,  where $\beta>0$.

\smallskip

For an open set $D\subset M$, let $\sF^D$  be the $\sE_1$-closure of $\sF \cap C_c(D)$ in $\sF$. Then $(\sE, \sF^D)$ is a regular Dirichlet form on $L^2(D;\mu)$, which is called the restriction of $(\sE,\sF)$ to $D$.  Let $(P_t^D)_{t \ge 0}$ denote the semigroup   associated with $(\sE,\sF^D)$. Define the first exit time 
$\tau_D:= \inf \{ t>0 : X_t \notin D \}.$ 
For any $f\in \sB_b(M)$, a quasi-continuous version of $	P^D_tf$    is given by
\begin{align*}
	P^D_tf(x)= \E^x[f(X_t): t<\tau_D], \quad x \in D\setminus \sN.
\end{align*}

 The \textit{heat kernel} of  the semigroup $(P^D_t)_{t \ge 0}$ (if it exists) is a non-negative measurable function $p^D(t,x,y)$ on $(0,\infty) \times (D\setminus \sN') \times (D\setminus \sN')$ for a properly exceptional set $\sN'\supset \sN$ satisfying the following properties: 

\medskip

\setlength{\leftskip}{4.2mm}

\noindent (1) $p^D(t,x,y)=p^D(t,y,x) $ for all $t>0$ and $ x,y \in D \setminus \sN'$.

\noindent (2) For all $t>0$, $x \in D \setminus \sN'$ and any $f \in \sB_b(M)$,
$$
P^D_tf(x)
= \E^x[f(X_t): t<\tau_D]= \int_{D} p^D(t,x,y) f(y) \mu(dy).
$$

\noindent (3) For all $s,t>0$ and $x,y  \in D \setminus \sN'$,
$$
p^D(s+t,x,y) =\int_{D} p^D(s,x,z) p^D(t,z,y) \mu(dz).
$$

\noindent (4) There exists an $\sE$-nest $(K_n)_{n \ge 1}$  of compact sets so that $\sN'=D \setminus \cup_{n =1}^\infty K_n$ and that for all $t>0$ and $y \in E \setminus \sN'$, the map $x \mapsto p^D(t,x,y)$ is continuous on each $K_n$.

\setlength{\leftskip}{0mm}

\smallskip

\begin{defn}
	\rm   We say that  
	 \textit{a near-diagonal lower estimate for (Dirichlet) heat kernel} \NDL \ holds if 	there exist constants  $q_0 \ge 1$, $\eta_1,\eta_2 \in (0,1)$ and   $C>0$ such that for all $x_0 \in M$ and $r\in (0,R_0/q_0)$, the  heat kernel $p^{B(x_0,r)}$ of  $(P^{B(x_0,r)}_t)_{t \ge 0}$ exists and  satisfies
	\begin{equation*}
		p^{B(x_0,r)}(t,x,y) \ge \frac{C}{V(x_0, \phi^{-1}(x_0,t))}
	\end{equation*}
	for  all $0<t \le \eta_1\phi(x_0, r)$ and $\mu$-a.e. $x,y \in B(x_0, \eta_2 \phi^{-1}(x_0,t))$.
\end{defn}

\subsection{Admissible weight function and new form of the tail condition for  jump measure}

We now introduce an admissible weight function and, using this, present  a new condition on the tail of the jump measure, which is weaker than the standard one.

\begin{defn}\label{d:weight}
	\rm
	We say that $\p:M \times (0,\infty)\to [1,\infty]$ is an \textit{admissible weight function} for $(\sE,\phi)$, if 
	there exist constants $\gamma_0,\gamma_1>0$ and $C\ge 1$ such that 	 the following  hold:
	\begin{flalign}
		\indent	\text{(1)  For each $x \in M$, the map $r \mapsto \p(x,r)$ is  non-decreasing and  continuous.}&&\label{e:theta-ic}
	\end{flalign} 
	
	\vspace{-0.07in}

	\setlength{\leftskip}{4.2mm}

\noindent 	(2)  For all $x\in M$ and $0<r\le R<\infty$, 
	\begin{align}\label{e:blow-up-scale}
		\frac{\p(x,R)}{\p(x,r)}     \le C \frac{\phi(x,R)}{\phi(x,r)}  \bigg( \frac{R}{r} \bigg)^{-\gamma_0}.
	\end{align}
	
\noindent 	(3) 	(Volume compatibility) For all $x\in M$ and $r\in (0,R_0)$,
	\begin{align}\label{e:blow-up-integral-general}
		\int_{B(x,r)} \p(y,r)^{1+\gamma_1} \mu(dy)\le C V(x,r).
	\end{align}
	
	\noindent (4) (Hardy-type inequality) For all $x\in M$, $r\in (0,R_0)$ and $u\in \sB_c(M)$,
	\begin{align}\label{e:blow-up-Hardy}
		\int_{B(x,r)} u(y)^2 \p(y,r)^{1+\gamma_1}\mu(dy)&\le C\big( \phi(x,r)\sE(u,u)+ \lVert u \rVert_{2}^2 \big).
	\end{align}		
	\setlength{\leftskip}{0mm}
\end{defn}

\smallskip

The inequality \eqref{e:blow-up-Hardy} is an $L^2$-Hardy-type inequality with the weight  $\p(y,r)^{1+\gamma_1}$, which will plays a crucial role in establishing the Caccioppoli and $L^2$-mean value inequalities in Section \ref{s:L2}.

Throughout this paper,  $\phi$ always stands for a scale function on $M$ and $\p$ for an admissible  weight function for  $(\sE,\phi)$.

\begin{remark}
	\rm
	(i) For any $c \ge 1$, the constant map $\p \equiv c$ is an admissible weight function for any scale function.
	
	\noindent (ii) By \eqref{e:phi-scale}, if there exist constants $\beta_3 \in [0,\beta_1)$ and $C\ge 1$ such that\begin{align*}	\frac{\p(x,R)}{\p(x,r)} \le C\left( \frac{R}{r} \right)^{\beta_3} \quad \text{for all} \;\, x \in M \text{ and } 0<r\le R,\end{align*}then \eqref{e:blow-up-scale} holds with $\gamma_0=\beta_1-\beta_3$, where $\beta_1>0$ is the constant in \eqref{e:phi-scale}.
	
	\noindent (iii) Since $\p \ge 1$, the inequality in  \eqref{e:blow-up-integral-general}  is, in fact, a comparability. 
	
	\noindent (iv) It follows from \eqref{e:blow-up-integral-general} and H\"older inequality that 
	for all $A \in \Sigma_M$, $x\in M$  and $r\in (0,R_0)$,
	\begin{equation}\label{e:blow-up-integral0}
		\int_{B(x,r)\cap A} \p(y,r) \mu(dy)\le
		C^{1/(1+\gamma_1)}
		\mu(B(x,r)\cap A)^{\gamma_1/(1+\gamma_1)}
		V(x,r)^{1/(1+\gamma_1)}.
	\end{equation}
	In particular, it holds that 
	\begin{equation}\label{e:blow-up-integral}
		\int_{B(x,r)} \p(y,r) \mu(dy)\le C^{1/(1+\gamma_1)} V(x,r) \quad \text{for all} \;\, x \in M \text{ and } r\in (0,R_0).
	\end{equation}
\end{remark}
\begin{defn}\label{d:Tail}
	\rm  We say that 		the jump measure $J$ satisfies \Tail\
	if	 there exist a non-negative function $J(x,dy)$ on $M \times \Sigma_M$  and a constant $C_J \in [0,\infty)$ such that the following hold:
	
	\medskip
	
	\setlength{\leftskip}{4.2mm}
	
	\noindent 	(1)			 For each  $x \in M$, $J(x,dy)$ represents a  measure on $\Sigma_M$, and  for every $E \in \Sigma_M$,  $x \mapsto J(x,E)$ is a measurable function on $M$. Moreover, 
	\begin{equation}\label{e:Tail-transition}
		J(dx,dy)=J(x,dy) \mu(dx)  \quad \text{in  $M\times M$.}
	\end{equation}

	\noindent	(2) For all $x \in M$ and $r\in (0,R_0)$,
	\begin{equation}\label{e:Tail}
		J(x, B(x,r)^c) \le C_J \frac{ \p(x,r)}{\phi(x,r)}.
	\end{equation}
\end{defn}

\subsection{Functional inequalities}

For an open set $D\subset M$, let $\sL^D$ be the $L^2$-generator of $(\sE, \sF^D)$ and let $\lambda_1(D)$ denote the bottom of the spectrum of $-\sL^D$. Then $\lambda_1(D)$ has the following variational formula:
\begin{align}\label{e:eigen}
	\lambda_1(D)= \inf \left\{ \sE(f,f) : f \in \sF^D \mbox{ with } \Vert f \Vert_{2} = 1  \right\}.
\end{align}

\begin{defn}\label{d:FK}	{\rm We say that the \textit{Faber-Krahn inequality} \FK \ holds if there exist constants  $C, \nu>0$ and $q_1 \ge 1$  such that for any $x_0 \in M$, $r \in (0, R_0/q_1)$ and  non-empty open set $D \subset B(x_0,r)$,
		\begin{equation*}
			\lambda_1(D) \ge \frac{C}{\phi(x_0,r)} \left( \frac{V(x_0,r)}{\mu(D)}\right)^\nu.
		\end{equation*}
	}
\end{defn}

\smallskip

Let $\sF_b:=\sF \cap L^\infty(M;\mu)$. For 
$B \in \Sigma_M$ with  $\mu(B) \in (0,\infty)$ and  $f\in L^1(B;\mu)$, define $$ \overline f_B:=\frac{1}{\mu(B)} \int_B f d\mu.$$

\begin{defn}\label{d:PI}
	\rm
	We say that \textit{Poincar\'e inequality} \PI \   holds if there exist constants $C>0$  and $q_2 \ge1$ such that for any $x_0 \in M$, $r \in (0, R_0/q_2)$ and  any $f \in \sF_b$,
	\begin{align*}
		&\int_{B(x_0,r)} (f- \overline{f}_{B(x_0,r)})^2 \, d\mu \\
		&\le C\phi(x_0,r) \bigg(\int_{B(x_0,q_2 r)} d\Gamma^{(L)}(f,f) + \int_{B(x_0,q_2r) \times B(x_0,q_2 r)} (f(x) - f(y))^2 J(dx,dy) \bigg).
	\end{align*}
\end{defn}

\smallskip

We note that conditions \FK \ and \PI \ are independent of the weight function $\p$.

 For any $f\in\sF_b$, there exists a unique positive Radon measure  $\Gamma(f,f)$  on $M$ such that  
\begin{align*}
	\int_M g \,d \Gamma(f,f)=\sE(f,fg)-\frac12 \sE(f^2,g) \quad \text{ for all $g \in \sF\cap C_c(M)$.}
\end{align*}
The measure $\Gamma(f,f)$  can be uniquely extended to any $f \in \sF$ as the increasing limit of $\Gamma(f_n,f_n)$, where $f_n:=((-n) \vee f) \wedge n$. See \cite[Section 3.2]{FOT}.  $\Gamma(f,f)$  is called the \textit{energy measure} or the \textit{carr\'e du champ} of $f$ for $\sE$. For $f \in \sF$, we denote by  $\Gamma^{(L)}(f,f)$  and  $\Gamma^{(J)}(f,f)$ the energy measures for $\sE^{(L)}$ and $\sE^{(J)}$, respectively.

Let  $\sF':=\{ f + a : f \in \sF, \, a \in \R\}$ and $\sF_b':=\sF' \cap L^\infty(M;\mu)$. Since $(\sE,\sF)$ has no killing part, the  form $\sE$ and the energy measure $\Gamma(f,f)$ can be extended to functions from $\sF'$ by
\begin{align*}
	&\sE(f,g):=	\sE(f+a,g+b), \qquad d\Gamma(f,f):=d\Gamma(f+a,f+a)
\end{align*}
for  $f,g \in \sF'$ with $f-a$, $g-b\in \sF$.
Similarly,  $\Gamma^{(L)}(f,f)$  and  $\Gamma^{(J)}(f,f)$  can be extended to functions from $\sF'$.

For a constant $\kappa \ge 1$ and open sets  $U,V\subset M$ with $U \Subset V$, a function $\vp$ on $M$ is called a \textit{$\kappa$-cutoff function} for $U \Subset V$, if $0\le \varphi \le \kappa$ in $M$, $\varphi \ge 1$ in $U$ and $\varphi=0$ in $V^c$. Any $1$-cutoff function is  referred to as a \textit{cutoff function}.

We introduce weighted versions of  cutoff Sobolev inequalities.

\begin{defn}\label{d:CS}
	\rm (i) We say that condition \CS \ holds if there exist constants $C_0\in (0,1]$ and $C_1,C_2>0$ such that for any $x_0 \in M$, $0<r\le R$ with $R+2r<R_0$  and any 	$f\in \sF'_b$, there exists a cutoff function $\vp \in \sF_b$ for $B(x_0,R) \Subset B(x_0,R+r)$ so that the following  holds:
	\begin{align}\label{e:CS}
		&\int_{B(x_0,R+(1+C_0)r)} f^2 d\Gamma(\varphi,\varphi) \nn\\
		&\le   C_1\bigg(\int_{B(x_0,R+r)} \varphi^2 d\Gamma^{(L)}(f,f) + \int_{B(x_0,R+r) \times B(x_0,R+(1+C_0)r)} \vp(x)^2 (f(x) - f(y))^2 J(dx,dy) \bigg) \nn\\
		&\quad   + \sup_{z \in B(x_0, R +(1+C_0)r)} \frac{C_2}{\phi(z,r)} \int_{B(x_0,R+(1+C_0)r)} f(x)^2 \p(x,r) \mu(dx). 
	\end{align}	
	(ii) We say that condition \CSU \ holds if there exist constants $C_0\in (0,1]$ and $C_1,C_2>0$ such that for any $x_0 \in M$ and  $0<r\le R$ with $R+2r<R_0$,  there exists a cutoff function $\vp\in \sF_b$ for $B(x_0,R) \Subset B(x_0,R+r)$ so that \eqref{e:CS} holds for any $f\in \sF'$.
\end{defn}

\begin{remark}\label{r:remark-CS}
	\rm  
	When $\p=1$ and $R_0=\infty$, the  cutoff Sobolev inequality  ${\mathrm{CS}}^1(\phi)$ with both strongly local and pure jump Dirichlet forms,  which is  a combination of CSA$(\phi)$ for strongly local Dirichlet forms  (see \cite{AB15})  and CSJ$(\phi)$  for pure jump Dirichlet forms (see \cite{CKW-memo}), was introduced and studied in \cite{CKW-adv}.  The condition \CSU, which is an abbreviation for \textit{cutoff Sobolev inequality with universal cutoff function}, is an adaptation of SCSJ$(\phi)$  for pure jump Dirichlet forms introduced in \cite{CKW-memo}.
\end{remark}

\subsection{Main results}

For an open set $D \subset M$, let $(G^D_\lambda)_{\lambda>0}$  be the  resolvent corresponding to  $(\sE,\sF^D)$.
The \textit{Green operator} $G^D$ is define  by
\begin{align}\label{e:def-Green}
	G^D f:=\lim_{\lambda \to 0}	G^D_\lambda f = \lim_{\lambda \to 0} \int_0^\infty e^{-\lambda t} P^D_tf dt
\end{align}
for  any non-negative measurable function $f$.  Recall that  
$\tau_D= \inf \{ t>0 : X_t \notin D \}$ denotes the first exit time from $D$. The mean exit time $\E^x\tau_D$ is known to satisfy  (see \cite[Lemma 3.1]{GT12}):
\begin{align}\label{e:mean-exit-time}
	G^D \1_D(x)=\E^x\tau_D \quad \text{for a.e.} \;\, x \in M.
\end{align} 
\begin{defn}\label{d:E}
	{\rm  We say that condition \mEu \ holds  if there  exist constants $C>1$ and $q_3 \ge 1$ such that
		\begin{equation*}
			\esssup_{y \in B(x_0,r)} \E^{y}[\tau_{B(x_0,r)}]  \le C \phi(x_0,r) \quad \text{for all} \;\, x_0 \in M, \; r\in (0, R_0/q_3),
		\end{equation*}
		and  condition \mEl \ holds  if there  exists a constant $C>1$ such that
		\begin{equation*}
			\essinf_{y \in B(x_0,r/4)}	 \E^{y}[\tau_{B(x_0,r)}]  \ge C^{-1} \phi(x_0,r)  \quad \text{for all} \;\, x_0 \in M, \; r\in (0, R_0).
		\end{equation*}
		We say that condition \mE \ holds if both \mEu \ and \mEl  \ hold.
	}
\end{defn}

We give  the probabilistic definitions of harmonic and caloric functions.

Let $Z:=(T_s,X_s)_{s \ge 0}$ be the time-space process corresponding to $X$ where $T_s:=T_0-s$. The law of the time-space process $s \mapsto Z_s$ starting from $(t,x)$ will be denoted by $\P^{(t,x)}$.  For an open subset $U$ of $[0,\infty) \times M$, define $\wh \tau_U:=\inf\{t>0:Z_t \notin U\}$.

A set $E \subset [0,\infty) \times M$ is said to be \textit{nearly Borel measurable} if for any probability measure $m$ on $[0,\infty) \times M$, there are Borel  subset $E_1,E_2$ of $[0,\infty) \times M$ such that $E_1\subset E \subset E_2$ and $\int_M \P^y(\cup_{t \ge 0} \{Z_t \in E_2 \setminus E_1\})m(dy)=0$.

\begin{defn}
	{\rm (i)  A nearly Borel 
		function $u$ on $M$ is said to be \textit{harmonic} (resp. \textit{subharmonic, superharmonic}) in an open set $V\subset M$ (with respect to  $X$), if there is a properly exceptional set $\sN_0\supset \sN$ such that  for any relatively compact  subset $U \subset V$,  $t\mapsto u(X_{t \wedge \tau_U})$ is a uniformly integrable martingale (resp. submartingale, supermartingale) under $\P^x$ for all $x \in U\setminus \sN_0$. 
		
		\noindent (ii) A nearly Borel 
		function $q$ on $[0,\infty) \times M$ is said to be \textit{caloric} in $(a,b] \times B(x_0,r)$ (with respect to $X$), if there is a properly exceptional set $\sN_0 \supset \sN$ such that for any relatively compact open subset $U \subset (a,b] \times B(x_0,r)$, it holds that $q(t,x)=\E^{(t,x)} q(Z_{\wh \tau_U})$ for all $(t,x) \in U \cap ([0,\infty) \times (M \setminus \sN_0))$.
		
	}
\end{defn}

\begin{defn}
	{\rm  (i) We say that \textit{elliptic H\"older regularity} \EHR \ holds for $(\sE,\sF)$, if there exist constants $C,\theta>0$ and $\delta \in (0,1)$ such that for all $x_0 \in M$,  $R \in (0,R_0)$ and  any  $u \in \sB_b(M)$		that is harmonic in $B(x_0,R)$, 
		\begin{align}\label{e:EHR}
			|u(x)-u(y)| \le C \bigg(\frac{d(x,y)}{R} \bigg)^\theta \lVert u \rVert_{\infty} \quad \text{for $\mu$-a.e.} \;\, x,y \in B(x_0, \delta R). 
		\end{align}

		\noindent(ii) We say that \textit{parabolic H\"older regularity} \PHR \ holds for $(\sE,\sF)$, if there exist constants $C,\theta>0$ and $\delta \in (0,1)$ such that for all $x_0 \in M$, $R \in (0,R_0)$,  $t_0 \ge 0$ and any bounded measurable function $q$ on $[t_0,t_0+\phi(x_0,R)] \times M$ that is caloric in $[t_0,t_0+\phi(x_0,R)] \times B(x_0,R)$, 
		\begin{align}\label{e:PHR}
			|q(s,x)-q(t,y)| \le C \bigg(\frac{\phi^{-1}(x_0,|s-t|)+d(x,y)}{R} \bigg)^\theta \esssup_{[t_0, t_0+\phi(x_0,R)] \times M} |q|
		\end{align}
		for all $s,t \in [t_0+  \phi(x_0, R) - \phi(x_0, \delta R), t_0 +  \phi(x_0,  R)]$ and $\mu$-a.e. $x,y \in B(x_0, \delta R)$.}
\end{defn}

\begin{remark}\label{r:EHR}
	{\rm In the definition of \EHR \ (resp. \PHR), if  \eqref{e:EHR} (resp.  \eqref{e:PHR}) holds for some $\delta \in (0,1)$, then it holds for any $\delta \in (0,1)$, with a possibly different constant $C$  by a standard covering argument. See \cite[Remark 1.14(iv)]{CKW-jems} for example.}
\end{remark}

We now present our  main result that establishes  the stability of \NDL \ and \PHR \ $+$ \mE \  for regular Dirichlet forms that may admit blow-up jump kernels.

\begin{thm}\label{t:main}
	Suppose that \VD, \RVD \ and \Tail \  hold. Then the following equivalences hold:
	\begin{align}\label{e:main}
		\text{\rm \NDL  \ } 	\;	  \Leftrightarrow \ \;
		\text{\rm \PHR \ + \mE \ }   \;
		\Leftrightarrow \ \;	\text{\rm  \EHR  \ + \mE \   }  \; \Leftrightarrow  \ \;	\text{\rm  \CS \ + \PI}.
	\end{align}
\end{thm}

\smallskip

Theorem \ref{t:main} follows from Theorem \ref{t:main-2}, which also establishes additional equivalence hypotheses for \eqref{e:main}.

There is an alternative functional inequality, often used in the study of Dirichlet forms on metric measure spaces, that plays essentially the same role as the cutoff Sobolev inequalities. This is referred to as the  \textit{generalized capacity condition}, introduced in \cite{GHL} for strongly Dirichlet forms, and in \cite{GHH18, GHH} for non-local forms and more general settings.
We  introduce weighted versions of the generalized capacity condition and the one with universal cutoff functions. 	Since $\p \ge 1$,  \Gcap \ and \GU \ below are weaker than the original  conditions (Gcap) and  (GU) in \cite{GHH18, GHH}, respectively.

\begin{defn}\label{d:Gcap}
	\rm (i) We say that condition \Gcap \ holds if there exist constants $\kappa \ge 1$ and $C>0$ such that for any $x_0 \in M$, $0<r\le R$ with $R+r<R_0$ and any $f \in \sF'_b$, there exists a $\kappa$-cutoff function $\vp \in \sF_b$ for $B(x_0,R) \Subset B(x_0,R+r)$ so that  the following holds:
	\begin{align}\label{e:Gcap}
		\sE(f^2 \vp, \vp) \le  \sup_{z \in B(x_0, R +r)} \frac{C}{\phi(z,r)} \int_{B(x_0,R+r)} f(x)^2 \p(x,r) \mu(dx). 
	\end{align}

	\noindent (ii)	We say that condition \GU \ holds if there exist constants $\kappa \ge 1$ and $C>0$ such that for any $x_0 \in M$, $0<r\le R$ with $R+r<R_0$, there exists a $\kappa$-cutoff function $\vp \in \sF_b$ for $B(x_0,R) \Subset B(x_0,R+r)$ so that  \eqref{e:Gcap}  holds for  all $f \in \sF'_b$.
\end{defn}

\smallskip

In establishing Theorem \ref{t:main}, we also prove the equivalence of mean exit time estimates,  cutoff Sobolev inequalities and  generalized capacity conditions. For regular Dirichlet forms with a non-blow-up jump kernel $(\p=1)$, these relations first appeared in \cite{CKW-memo, GHH18} and were later extended to a more general setting in \cite{GHH}.
\begin{thm}\label{t:main-1}
	Suppose that \VD, \Tail \ and \FK \  hold. Then the following equivalences hold:
	\begin{align}\label{e:main-1}
		\text{\rm  \CS \  } \; \Leftrightarrow  \	\;\text{\rm  \CSU \  }  \; \Leftrightarrow  \	\; \text{\rm  \Gcap \  } 
		\;	\Leftrightarrow  \ \;
		\text{\rm  \GU \  }
		\;	\Leftrightarrow  \ \;
		\text{\rm \mE}.
	\end{align}
	In particular, if  \VD, \RVD, \Tail \ and \PI \  hold, then \eqref{e:main-1} holds.
\end{thm}
The proof of Theorem \ref{t:main-1} will be presented in Section \ref{s:proof-thm-1}.

  \section{Preliminary results}\label{s:preliminary}
  
  In this section, we introduce a condition on the pointwise upper bound of the jump  density and   state some preliminary (known) results. 
  
  Recall that $\phi$ is a scale function satisfying \eqref{e:phi-ic}--\eqref{e:phi-comp}, and $\p$ is an admissible weight function for $(\sE,\phi)$ satisfying \eqref{e:theta-ic}--\eqref{e:blow-up-Hardy}. Note that  $\phi$ satisfies
  \begin{align}\label{e:phi2}
  	\frac{\phi(y,R)}{\phi(x,r)} \le C\bigg(\frac{R + d(x,y)}{r}\bigg)^{\beta_2}  \quad \text{for all} \;\, x,y\in M \text{ and } 0<r \le R<\infty.
  \end{align}
  Indeed, for all $x,y \in M$ and $0<r \le R$, by \eqref{e:phi-scale} and \eqref{e:phi-comp}, we have
  \begin{align*}
  	&\frac{\phi(y,R)}{\phi(x,r)} \le  c_U\frac{\phi(y,R+d(x,y))}{\phi(x,R+d(x,y))} \bigg(\frac{R + d(x,y)}{r}\bigg)^{\beta_2}   \le  c_1c_U\bigg(\frac{R + d(x,y)}{r}\bigg)^{\beta_2}.
  \end{align*}
  \begin{defn}
  	\rm (i) We say that condition  \Ju \ holds if	 there exist a non-negative symmetric measurable function $J(x,y)$ on $M \times M$ and a constant $C>0$ such  that 
  	\begin{equation}\label{e:jump-density}	J(dx,dy)=J(x,y) \,\mu(dx)\mu(dy),
  	\end{equation}
  	and for all $x,y \in M$,
  	\begin{equation*}
  		J(x,y) \le C\frac{\p(x,d(x,y)) \p(y,d(x,y))}{V(x,d(x,y))\phi(x,d(x,y))}.
  	\end{equation*}
  	(ii)	We say that condition \Jl \ holds,  if \eqref{e:jump-density} holds and there exists a constant $C>0$ such  that  for all $x,y \in M$,
  	\begin{equation*}
  		J(x,y) \ge \frac{C}{V(x,d(x,y))\phi(x,d(x,y))}.
  	\end{equation*}
  \end{defn}

  \begin{defn}
  	\rm  When \Tail \ (resp.  \CSU, \Gcap \ or \Ju) holds with the trivial weight function $\p \equiv 1$, we simply say that \Tailq \ (resp. \CSUq,  \Gcapq \ or \Juq) holds.
  \end{defn}

  An unweighted version of the next lemma is standard. See \cite[Lemma 2.1]{CK08} or  \cite[Lemma 2.1]{CKW-memo}. We extend this result to our setting, which involves the admissible weight function $\p$.

  \begin{lem}\label{l:Ju-Tail} 
  Assume  \VD \  with $R_0=\diam(M)$. Then \Ju\ implies \Tail \  with $R_0=\diam(M)$.
  \end{lem}
  \pf   \eqref{e:Tail-transition} follows from \eqref{e:jump-density}.  Let $x \in M$ and $r\in (0,\diam(M))$. 
  Using \Ju \ in the first line below, \eqref{e:blow-up-scale} in the second, and the monotonicity of $\p$ in the third, we get
  \begin{align}\label{e:Tail-check-1}
  	J(x, B(x,r)^c)& \le c_1\int_{B(x,r)^c}   \frac{\p(x,d(x,y))  \p(y,d(x,y)) }{V(x,d(x,y))\phi(x,d(x,y))} \mu(dy) \nn\\
  	&\le \frac{c_2 r^{\gamma_0} \p(x,r)}{\phi(x,r)}\int_{B(x,r)^c}    \frac{\p(y,d(x,y))}{V(x,d(x,y)) d(x,y)^{\gamma_0}}\mu(dy) \nn\\
  	&	\le\frac{c_2 \p(x,r)}{\phi(x,r)}  \sum_{n \ge 0: \,2^n r <\diam(M)} \frac{2^{-n\gamma_0}}{V(x,2^n r)}\int_{B(x,2^{n+1}r) \setminus B(x,2^nr)}   \p(y,2^{n+1}r)\mu(dy). 
  \end{align}
  By \eqref{e:blow-up-integral} and \VD, we have
  \begin{align}\label{e:Tail-check-2}
  	&\sum_{n \ge 0: \,2^{n+1} r <\diam(M)} \frac{2^{-n\gamma_0}}{V(x,2^n r)}\int_{B(x,2^{n+1}r) \setminus B(x,2^nr)}   \p(y,2^{n+1}r)\mu(dy)\nn\\
  	&\le c_3\sum_{n \ge 0: \,2^{n+1} r <\diam(M)} \frac{2^{-n\gamma_0} V(x,2^{n+1}r)}{V(x,2^n r)}\le c_4\sum_{n \ge 0} 2^{-n\gamma_0} = c_5.
  \end{align}
  Moreover, note that by \eqref{e:blow-up-scale} and \eqref{e:phi-scale}, there exists $c_6>0$ such that $\p(y,2r) \le c_6\p(y,r)$ for all $y \in M$ and $r>0$. Hence,  when $\diam(M)<\infty$, for $n$ satisfying $2^{n} r <\diam(M) \le   2^{n+1} r $, using \eqref{e:blow-up-integral} and \VD,
  we obtain 
  \begin{align}\label{e:Tail-check-3}
  	&\frac{2^{-n\gamma_0}}{V(x,2^n r)}\int_{B(x,2^{n+1}r) \setminus B(x,2^nr)}   \p(y,2^{n+1}r)\mu(dy)\le \frac{1}{V(x,2^n r)}\int_{M}   \p(y,2^{n+1}r)\mu(dy)\nn\\
  	&\le \lim_{s \to \diam(M)} \frac{c_6}{V(x,2^nr)}\int_{B(x,s)}   \p(y,s)\mu(dy) \le \lim_{s \to \diam(M)} \frac{c_7V(x,s)}{V(x,\diam(M)/2)} \le c_8.
  \end{align}
  By \eqref{e:Tail-check-1}, \eqref{e:Tail-check-2} and \eqref{e:Tail-check-3}, we deduce that \eqref{e:Tail} holds with $C_J=c_2(c_5+c_8)$. \qed 
  
  The following is a simple  consequence of  \Tail.
  \begin{lem}\label{l:Tail-integral}
  	Suppose that   \Tail \ holds. For any $\gamma>\beta_2$, where $\beta_2$ is the constant in \eqref{e:phi-scale}, there exists $C>0$ depending on $\gamma$ such that for all $x \in M$ and $r \in (0,R_0)$,
  	\begin{align*}
  		&	\int_{M}  \bigg(1 \wedge \frac{d(x,y)}{r} \bigg)^\gamma J(x,dy)  \le   C\frac{\p(x,r)}{\phi(x,r)}.
  	\end{align*}
  \end{lem}
  \pf Let $x \in M$ and $r \in (0,R_0)$. By \Tail,  $ J(x,B(x,r)^c)  \le   c_1\p(x,r)/\phi(x,r).$ Further, using   \Tail \ in the second inequality below and the monotonicity of $\p$ and \eqref{e:phi-scale} in the third, we obtain
  \begin{align*}
  	&r^{-\gamma}\int_{B(x,r)} d(x,y)^\gamma  J(x,dy) 
  	\le \sum_{n=0}^\infty  2^{-\gamma n} \int_{B(x,2^{-n}r) \setminus B(x,2^{-n+1}r)}  J(x,dy)\\
  	&\le c_2 \sum_{n=0}^\infty 2^{-\gamma n} \frac{\p(x,2^{-n+1}r)}{\phi(x,2^{-n+1}r)} \le \frac{c_3\p(x,r)}{\phi(x,r)} \sum_{n=0}^\infty 2^{-(\gamma-\beta_2)n} = \frac{c_5\p(x,r)}{\phi(x,r)}.
  \end{align*} \qed
  
  We now discuss some results related to \FK.  Recall that for an open set $D \subset M$,  $\lambda_1(D)$ is  given by \eqref{e:eigen}.   The following relation between the bottom of the spectrum and mean exit time was established in \cite[Lemma 6.2]{GH}:
  \begin{align}\label{e:lambda1}
  	\lambda_1(D) \ge \Big(\esssup_{y \in D} \,\E^y \tau_D \Big)^{-1}.
  \end{align}
  
  \begin{defn}
  	\rm We say that  \textit{homogeneous local Nash inequality} \HLN \ holds if there exist  $C,\nu>0$ and $q_1 \ge 1$ such that for all $x_0\in M$,  $r\in (0,R_0/q_1)$ and  $f \in \sF^{B(x_0,r)}\cap L^\infty(M;\mu)$,
  	\begin{align*}
  		\lVert f \rVert_{2}^{2+2\nu}\, \lVert f \rVert_{1}^{-2\nu} \le C \frac{\phi(x_0,r) }{V(x_0,r)^\nu}  \sE(f,f). 
  	\end{align*}
  \end{defn}
 
  \begin{defn}\label{d:SP}
  	{\rm  We say that \textit{survival probability estimate} \SP \ holds  if there  exist constants $c_1,c_2>0$ and $q_3 \ge 1$ such that for all $x_0 \in M$,  $r \in (0, R_0/q_3)$ and $t>0$,
  		\begin{equation*}
  			\esssup_{y \in B(x_0,r)} \P^y (\tau_{B(x_0,r)} >t)  \le c_1 e^{-c_2 t/\phi(x_0,r)}.
  		\end{equation*}
  	}
  \end{defn}
  
  \smallskip

  Several results have established the equivalence between Faber-Krahn-type inequalities, various Nash-type inequalities, and upper bounds for heat kernels (see \cite{BCS, GH, CKW-memo}). In our context, we obtain the following result.
  
  \begin{prop}\label{p:FK}
  	The following statements are equivalent:

  	\smallskip
  	
  	\setlength{\leftskip}{4mm}
  	
  	\noindent	{\rm (1)} \FK \ holds.
  	
  	\noindent	{\rm (2)} \HLN \ holds.
  	
  	\noindent 	{\rm (3)} For all $x_0\in M$ and  $r \in (0,R_0/q_1)$,   the heat kernel $p^{B(x_0,r)}$ associated with the semigroup $(P^{B(x_0,r)}_t)_{t\ge0}$ exists and has the following estimate: There exist constants $C,\nu>0$ and $q_1 \ge 1$ such that
  	\begin{align}\label{e:FK-Dirichlet}
  		\esssup_{y,z \in B(x_0,r)} p^{B(x_0,r)}(t,y,z) \le C\frac{\phi(x_0,r)^{1/\nu}}{V(x_0,r)} \,t^{-1/\nu} \quad \text{for all} \;\, t>0.
  	\end{align}
  	
  	\setlength{\leftskip}{0mm}
  	
  	Moreover, any of the above {\rm (1)--(3)} implies \SP, and thus \mEu.
  \end{prop}
  \begin{cor}\label{c:PI-FK}
  	Suppose that \VD\ and \RVD \ hold. Then  we have
  	\begin{align*}
  		\text{\PI \, \ $\Rightarrow$ \;  \HLN \ \, $\Leftrightarrow$   \, \FK}.
  	\end{align*}
  \end{cor}
  We prove Proposition \ref{p:FK} and \ref{c:PI-FK} using two variants of Nash-type inequalities, with the proof provided in Appendix \ref{s:A}.

  Recall that $Z=(T_s,X_s)_{s\ge 0}$ is the time-space process corresponding to $X$.
  For a Borel subset $E\subset [0,\infty) \times M$, define  $\sigma_E:=\inf\{t>0:Z_t \in E\}$.   Denote by $\wt\mu$ the product of the Lebesgue measure on $\R_+$ and $\mu$ on $M$. The next two lemmas are standard and their proofs are given in Appendix \ref{s:A2}.
  \begin{lem}\label{l:PHR-1}
  	Suppose that \VD \ and \NDL \ hold. 
  	There exist constants $0<\delta_1< \delta_2\le 1/4$ and $C>0$ such that for all $x_0 \in M$ and  $r \in (0,R_0)$, we have $3\phi(x_0,\delta_2 r) < \phi(x_0,r)$ and the following  holds: For any  $t \ge  \phi(x_0,r)$, $s \in [t-\phi(x_0,\delta_2 r),t]$  and any compact set $E \subset [t-3\phi(x_0,\delta_2 r), t-  2\phi(x_0,\delta_2 r)] \times B(x_0, \delta_1 r)$, it holds that 
  	\begin{align*}
  		\P^{(s,y)} \big(  \sigma_E < \wh \tau_{[t-\phi(x_0,r), t] \times B(x_0,r)} \big) \ge \frac{ C\,\wt \mu (E)}{V(x_0,r)\phi(x_0,r)} \quad \text{for all } y \in B(x_0, \delta_1r).
  	\end{align*} 
  \end{lem}
  \begin{lem}\label{l:NDL-1}
  	Suppose that \VD \ and \NDL \ hold. Then \PI \ and \mEl \ hold. Moreover, if \RVD \ also holds, then \FK \ and \mE \ hold.
  \end{lem}

  We end this section with an analytic formulation of harmonicity and the equivalence between probabilistic and analytic notions of harmonicity.
  
  Let $D\subset M$ be an open set. We say that a function $u$ is locally in $\sF^D$, denoted as $u \in \sF_{\rm loc}^D$, if for any relatively compact open set $U\subset D$, there exists $f \in \sF^D$ such that $u=f$ a.e. on $U$. 
   A nearly Borel  function $u$ on $M$ is said to be $\sE$-\textit{harmonic} (resp. \textit{$\sE$-subharmonic, $\sE$-superharmonic}) in $D$,  if  $u \in \sF_{\rm loc}^D$ is locally bounded on $D$ and satisfies the following two properties:
  
  \medskip

  \setlength{\leftskip}{4.2mm}
  
  \noindent (1) For any relatively compact open subset $U$ and $V$ of $D$ with $\overline U \subset V \subset \overline V \subset D$,
  $$
  \int_{U \times V^c} |u(y)| \, J(dx,dy)<\infty.
  $$
  
 \noindent (2) For any non-negative $\vp \in \sF\cap C_c(D)$, it holds that
  $$
  \sE(u,\vp)=0 \quad\;\;(\text{resp. } \sE(u,\vp)\le 0, \; \sE(u,\vp)\ge 0).
  $$
  
  \setlength{\leftskip}{0mm}

  \smallskip
  
  The following equivalence between probabilistic and analytic notions of harmonicity for bounded functions was proved in \cite{Ch} for harmonic functions and extended in \cite{ChKu} to subharmonic functions.
  
  \begin{thm}\label{t:harmonicity}
  	Let $D\subset M$ be an open set and $u$ be a bounded function on $M$. Then $u$ is harmonic (resp. subharmonic) in $D$ if and only if $u$ is $\sE$-harmonic (resp. $\sE$-subharmonic) in $D$.
  \end{thm}

\section{Examples of  admissible weight functions}\label{s:ex_weight}

In this section, we discuss explicit examples of admissible weight functions. To do this, we first introduce condition \Exi \ and some geometric concepts.

\subsection{Condition \Exi}
\begin{defn}
 	\rm 
 	For a given  $\alpha>1$, we say that 	 $(M,d,\mu)$ satisfies 	
	 \Exi,  if there exists a strongly local regular Dirichlet form $(\sE^{(c)}, \sF^{(c)})$ on $L^2(M;\mu)$ with the corresponding heat kernel $q(t,x,y)$ satisfying the following (sub)-Gaussian estimates: there exist $c_i>0$, $i=1, \dots, 4$, such that for all $t>0$ and $x,y \in M$,
 		\begin{align}\label{e:Gaussian}
 			\frac{c_1}{V(x, t^{1/\alpha})}  \1_{\{d(x,y) \le c_2 t^{1/\alpha}\}} \le 	q(t,x,y) \le \frac{c_3}{V(x,t^{1/\alpha})} \exp \left[ - c_4 \bigg( \frac{d(x,y)^\alpha}{t}\bigg)^{1/(\alpha-1)} \right].
 		\end{align}
 \end{defn}
 
 \smallskip

The Euclidean space $\R^d$ satisfies Exi$(2)$.  For  the $d$-dimensional unbounded  Sierpi\'nski gasket, \Exi \ holds with $\alpha=\log (d+3)/\log 2$. See \cite{BP88}. For the unbounded Vicsek set, \Exi \ holds with $\alpha=\log 15/\log 3$. See \cite{FHK94}.  We refer to \cite{Bar98, Te06} for further details and related results on \Exi.

 For $\beta>0$, let
 \begin{align}\label{e:def-Sobolev}
	W_{\beta/2}(u):=  \int_{M\times M} \frac{(u(x)-u(y))^2}{V(x,d(x,y))d(x,y)^\beta} \mu(dx)\mu(dy)
\end{align} 
and define  a Sobolev-Slobodeckij space $W^{\beta/2,2}(M)$ by
  \begin{align*}
 	W^{\beta/2,2}(M)&:=\left\{ u \in L^2(M;\mu): W_{\beta/2}(u)<\infty \right\}.
 \end{align*}

 A non-decreasing L\'evy process $\xi=(\xi_t)_{t\ge 0}$ on $[0,\infty)$ is called a \textit{subordinator}. 
 It is known that the Laplace exponent $\psi^\xi$ of a subordinator $\xi$, defined by
 \begin{align*}
 	\E\left[ \exp(-\lambda \xi_t)\right] = \exp\big(-t \psi^\xi(\lambda)\big), \quad \lambda,t>0,
 \end{align*}
is a Bernstein function.
That is, $\psi^\xi$ is a $C^\infty$ function on $(0,\infty)$ such that  $(-1)^n (\psi^\xi)^{(n)} \le 0$ for all $n \ge 0$. 
Every Bernstein function has the following representation: 
There exist constants $a,b\ge 0$ and a measure $\nu$ on $(0,\infty)$ satisfying $\int_0^\infty (1\wedge t)\nu(dt)<\infty$ such that
 \begin{align*}
 	\psi^\xi(\lambda)=a+b\lambda + \int_0^\infty (1-e^{-\lambda t})\nu(dt).
 \end{align*} 
 We refer to \cite{SSV12} for fundamental results on subordinators and Bernstein functions.

   Suppose that \Exi \ holds  with a strongly local regular Dirichlet form  $(\sE^{(c)},\sF^{(c)})$ and let $Z=\{Z_t,t \ge 0;\P^x, x \in M\}$ be the Hunt process associated with $(\sE^{(c)},\sF^{(c)})$.
  Let $\beta \in (0,\alpha)$  and  $\xi$ be a subordinator with the Laplace exponent 
$\psi^\xi(\lambda) = \int_0^\infty (1-e^{-\lambda t}) t^{-\beta/\alpha-1}dt$, 
 independent of $Z$. Define a process  $Y=\{Y_t,t \ge 0;\P^x, x \in M\}$ by $Y_t=Z_{\xi_t}$ for $t\ge 0$.  According to \cite[Theorems 2.1 and 3.1(i)]{Ok02}, the time-changed process $Y$ is a conservative  Hunt process on $M$ associated with  a  regular Dirichlet form $(\sE^\beta,\sF^\beta)$ defined by
 \begin{align}\label{e:def-sE-beta}
 	\begin{split} 
 		\sE^\beta(u,v)&=\frac12\int_{M\times M} (u(x)-u(y))(v(x)-v(y)) J_{\beta}(x,y)\mu(dx)\mu(dy), \\
 		 \sF^\beta&=\text{the closure of $\sF^{(c)}$ with respect to $(\sE^\beta_1(\cdot,\cdot))^{1/2}$-norm},
 	\end{split}
 \end{align}
 where 
 \begin{align*}
 	J_{\beta}(x,y):=\int_0^\infty q(t,x,y) t^{-\beta/\alpha-1}dt.
 \end{align*} 
 Moreover, since
 \begin{align*}
 	\int_0^\infty (t \wedge r) t^{-\beta/\alpha-1}dt = \frac{\alpha^2}{\beta(\alpha-\beta)}r^{1-\beta/\alpha} \quad \text{for all} \;\, r>0,
 \end{align*}
 by \cite[(2.3)]{Ok02} and Young's inequality,  we get that for all $r>0$ and $u \in \sF^{(c)}$,
 \begin{align}\label{e:sE-beta-sE-c}
 	r^\beta \sE^\beta(u,u) \le  c_1 r^\beta\sE^{(c)} (u,u)^{\beta/\alpha} \lVert u \rVert_2^{2(1-\beta/\alpha)}  \le \frac{c_1\beta}{\alpha} r^\alpha \sE^{(c)}(u,u)    + \frac{c_1(\alpha-\beta)}{\alpha}\lVert u \rVert_2^{2}.
 \end{align}

 Recall  that  $\phi_\beta(x,r)=r^\beta$.  By applying the stability results in  \cite{BKKL2, CKW-memo, GHH18}, we obtain the following proposition.

 \begin{prop}\label{p:subordinate-Exi}
 	Suppose that  \VD \ and \RVD \  hold with $R_0=\infty$, and \Exi \ holds. Then for every $\beta \in (0,\alpha)$, the Dirichlet form  $(\sE^\beta,\sF^\beta)$ defined in \eqref{e:def-sE-beta} satisfies     {\rm J}$_\le(\phi_\beta)$,  {\rm J}$_\ge(\phi_\beta)$,   {\rm Gcap}$(\phi_\beta)$ and {\rm PI}$(\phi_\beta)$, and has a  jointly continuous heat kernel $q^\beta(t,x,y)$ on $(0,\infty) \times M \times M$ satisfying
 \begin{align}\label{e:HKE-beta}
 	q^\beta(t,x,y)   \asymp \frac{1}{V(x,t^{1/\beta})} \wedge \frac{t}{V(x,d(x,y))d(x,y)^\beta}, \qquad t>0, \; x,y \in M.
 \end{align}
 Moreover, we have
 	\begin{align}\label{e:subordinate-norm}
 	\sE^{\beta}(u,u) \asymp W_{\beta/2}(u) \quad \text{for all} \;\, u \in  \sF^\beta \quad \text{ and } \quad \sF^\beta= W^{\beta/2,2}(M).
 \end{align}
 \end{prop}
 \pf  Let $\beta \in (0,\alpha)$.  According to  \cite[Section 6.1]{CKW-memo}, $(\sE^\beta, \sF^\beta)$ satisfies both {\rm J}$_\le(\phi_\beta)$ and  {\rm J}$_\ge(\phi_\beta)$, and has a heat kernel $q^\beta(t,x,y)$ which satisfies \eqref{e:HKE-beta}.
 Note that these results are proved under an additional assumption that $q(t,x,y)$ satisfies certain off-diagonal estimates in  addition to \eqref{e:Gaussian}, but \eqref{e:Gaussian} alone is sufficient to guarantee these results (cf.  \cite[Lemm 4.2]{BKKL2}). By \cite[Theorem 1.13, Remark 3.7 and Lemma 5.6]{CKW-memo} and \cite[Proposition 3.5]{CKW-jems}, it follows that   $(\sE^\beta, \sF^\beta)$ satisfies   {\rm Gcap}$(\phi_\beta)$ and {\rm PI}$(\phi_\beta)$, and that $q^\beta(t,x,y)$ is jointly continuous. 

Now, let us prove \eqref{e:subordinate-norm}. Since  {\rm J}$_\le(\phi_\beta)$ and  {\rm J}$_\ge(\phi_\beta)$ hold for $(\sE^\beta,\sF^\beta)$,  	$\sE^{\beta}(u,u) \asymp W_{\beta/2}(u)$ for all $u \in \sF^\beta$  and $\sF^\beta \subset W^{\beta/2,2}(M)$. 
To prove the opposite inclusion,   let   $u \in W^{\beta/2,2}(M)$. Denote by $(Q^\beta_t)_{t\ge 0}$ the semigroup of $(\sE^\beta,\sF^\beta)$.  For all $t>0$, by using the  conservativeness  and the symmetry of $Y$, and \eqref{e:HKE-beta},  we have
   \begin{align*}
   	&\frac1t  \la u- Q^\beta_t u, u\ra =	\frac1t \int_M u(x) \left( \int_M (u(x)-u(y)) q^\beta(t,x,y)\mu(dy) \right) \mu(dx)\nn\\
   &= \frac1{2t} \int_{M\times M}   (u(x)-u(y))^2 q^\beta(t,x,y)\mu(dy)  \mu(dx)\le c_1 \int_{M\times M}   \frac{(u(x)-u(y))^2}{V(x,d(x,y))d(x,y)^\beta} \mu(dy)  \mu(dx).
   \end{align*}
 Hence, $\sup_{t>0} t^{-1}\la u- Q^\beta_t u, u\ra\le c_1W_{\beta/2}(u) <\infty$. By  \cite[Lemma 1.3.4]{FOT}, this implies $u \in \sF^\beta$. The proof is complete.  \qed

\subsection{Some geometric concepts}
In this subsection, we recall some known geometric concepts. We begin with the Assouad dimension \cite{Ass83} and the Aikawa codimension   \cite{Ai91}.
\begin{defn}
\rm Let  $E \subset M$ be a non-empty subset.

\smallskip

\noindent (i) Let $\mathfrak I_1(E)$ be the set of all $\lambda>0$ satisfying the following property: There exists $C>0$ such that for all $x\in E$ and  $0<r<R<2\,\diam(E)$, the set $B(x,R) \cap E$ can be covered by at most  $N\le C(R/r)^\lambda$ open balls of radius $r$ centered in $E$. The \textit{Assouad dimension} of $E$ is defined by
\begin{align*}
{\rm dim}_{AS}(E):= \inf\{\lambda: \lambda \in \mathfrak I_1(E)\}.
\end{align*} 

\noindent (ii) Let $\mathfrak I_2(E)$ be the set of all $\lambda>0$ satisfying the following property: There exists  $C>0$ such that for all $x\in E$ and $0<r<\diam(E)$, 		\begin{align*}	\frac{1}{V(x,r)}	\int_{B(x,r)}  {\rm dist}(y,E)^{-\lambda} \mu(dy) \le Cr^{-\lambda}.		\end{align*}		 The \textit{Aikawa codimension} of $E$ is defined by		\begin{align*}		{\rm codim}_{AI}(E):= \sup\{\lambda: \lambda \in \mathfrak I_2(E)\}.		\end{align*} 
\end{defn}

\smallskip

 Note that if $M=\R^d$ and $D\subset \R^d$ is a  bounded Lipschitz open set, then the Assouad dimension of  the boundary $\partial D$ of $D$  is $d-1$.
We refer to \cite{KLV13} for relations among the Assouad dimension and other notions of dimension. In particular, by \cite[Remark 3.3 and Lemma 3.4]{KLV13}, we have:
\begin{lem}\label{l:dimensions}
	Suppose that \VD \ and \RVD \ hold with $R_0=\diam(M)$. 
For any non-empty subset $E\subset M$, we have
	 \begin{align*}
	{\rm dim}_{AS}(E)+	{\rm codim}_{AI}(E)\ge d_1,
	\end{align*}
	where $d_1>0$ is the constant in \RVD.
\end{lem}

We recall the notion of $\kappa$-fatness. A $\kappa$-fat open set is also known as  a $\kappa$-plump open set.

\begin{defn}\label{d:fat}
	{\rm Let $\kappa \in (0,1)$. An open set $D\subset M$  is \textit{$\kappa$-fat} if for any $x \in \overline D$ and $0<r<\diam(D)$, there exists $z \in \overline{B(x,r)}$ such that $B(z, \kappa r) \subset D$.}	
\end{defn}
By  \cite[Theorem 2.15]{Va88}, any uniform domain (thus any Lipschitz domain) in $\R^d$ is  $\kappa$-fat. Any $\kappa$-fat open set $D\subset \R^d$ is an Ahlfors $d$-regular set, namely, there exists $C>1$ such that
\begin{align}\label{e:d-set}	C^{-1}r^d\le 	|D \cap B(x,r)| \le Cr^d \quad \text{ for all} \;\, x \in M \text{ and } 0<r<\diam(D). \end{align}

\subsection{Examples}
In this subsection, we discuss examples of admissible weight functions. 
 We first consider admissible weight functions that blow  up to infinity at a fixed point.

Let $o \in M$. For $\theta>0$, we define 
\begin{align}\label{e:def-Hardy-theta}
\p_\theta(x,r)	=\p_\theta(x,r;o):=\bigg(1+ \frac{r}{d(x,o)}\bigg)^\theta, \qquad x \in M, \; r>0.
\end{align}

Recall that $\phi_{\beta}(x,r)=r^\beta$ for $\beta>0$.
\begin{prop}\label{ex:Hardy-origin} 
Suppose that $(M,d,\mu)$ satisfies 	
	 \Exi\ with a strongly local regular Dirichlet form   $(\sE^{(c)}, \sF^{(c)})$, and that
	  \VD \ and \RVD \ hold with $R_0=\infty$. Suppose also that 	  $\sF=\sF^{(c)}$ and 
	  there exists $C>1$ such that 
		\begin{equation}\label{e:diffusion-part}
	 C^{-1} \sE^{(c)}(u,u) \le 	\sE^{(L)}(u,u)\le C \sE^{(c)}(u,u) \quad \text{for all}\;\, u \in \sF.
		\end{equation}
	Then for any $o \in M$ and $\theta \in [0, \alpha \wedge d_1)$,    where $d_1$ is the constant in \RVD, the function $\p_\theta(x,r;o)$ defined in \eqref{e:def-Hardy-theta} 	is an admissible  weight function for  $(\sE^{(L)}, 	\phi_\alpha)$ with $R_0=\infty$.
\end{prop}
\pf  For each $x \in M$,  the map $r \mapsto \p_\theta(x,r)$ is non-decreasing and continuous, and  $\p_\theta$ satisfies \eqref{e:blow-up-scale} with $\gamma_0=\alpha-\theta$. 
  Let $\beta \in (0,\alpha \wedge d_1)$ be such that $\beta>\theta$. Set  $\gamma_1:=(\beta-\theta)/\theta$. 
By Lemma \ref{l:dimensions}, since ${\rm dim}_{AS}(\{o\})=0$, we have ${\rm codim}_{AI}(\{o\}) \ge d_1>\beta$. It follows that  for all $x\in M$ and $r>0$,
\begin{align*}
	&\frac{1}{V(x,r)}	\int_{B(x,r)} \p_\theta(y,r)^{1+\gamma_1} \mu(dy)\le   \frac{2^{\beta}}{V(x,r)} \bigg[ V(x,r)+ \int_{B(x,r)} \bigg( \frac{r}{d(y,o)}\bigg)^{\beta} \mu(dy) \bigg] \le c_0,
\end{align*}
proving that \eqref{e:blow-up-integral-general} holds.

 We now prove that \eqref{e:blow-up-Hardy} holds with $\gamma_1$.  Let $x \in M$, $r>0$ and  $u\in \sB_c(M)$ be such that $\sE(u,u)<\infty$. 
 Set $B':=B(x,r)$ and $B:=B(x,2r)$.  
Let $(\sE^\beta,\sF^\beta)$ be the Dirichlet form defined in \eqref{e:def-sE-beta}. By  applying \eqref{e:sE-beta-sE-c} with $r=1$ and \eqref{e:diffusion-part}, we obtain 
 $\sE^\beta(u,u) \le c_1 \sE^{(L)}_1(u,u) \le c_1\sE_1(u,u)<\infty$. Hence, $u \in \sF^\beta$.
 By Proposition \ref{p:subordinate-Exi},  Gcap$(\phi_\beta)$ holds for $(\sE^\beta,\sF^\beta)$. Thus, there exist  constants $\kappa \ge 1$, $c_1>0$  and  a $\kappa$-cutoff function $\vp \in \sF^\beta$ for $B' \Subset B$ independent of $u$ such that 
\begin{align}\label{ex:e:Gcap-1}
	\sE^{\beta}(u^2\vp,\vp) \le \frac{c_1}{r^\beta}\int_{B} u^2 d\mu.
\end{align} 
Since $\beta<d_1$, by Proposition \ref{p:subordinate-Exi} and \cite[Lemma 2.4]{CGL21},  condition $\mathbf{(G)_{\beta}}$ in \cite{CGL21} holds for $(\sE^\beta,\sF^\beta)$.  Applying   \cite[Theorem 5.6]{CGL21}, we deduce that
\begin{align}\label{ex:e:Gcap-2}
	\int_{B'}  \frac{ u(y)^2 \vp(y)^2}{d(y,o)^{\beta}} \mu(dy)\le c_2  \sE^{\beta}(u\vp,u\vp),
\end{align} 
where  $c_2>0$ is a constant independent of $x,r,u$ and $\vp$.   Besides, by \cite[(3.18)]{GHH18},  it holds that 
\begin{align}\label{ex:e:Gcap-3}
	\sE^{\beta}(u\vp,u\vp) =	\sE^{\beta}(u^2\vp,\vp) +  \int_{B \times B}  \vp(z)\vp(y)(u(z)-u(y))^2 J_\beta(x,y) \mu(dz)\mu(dy).
\end{align}
Using \eqref{ex:e:Gcap-2}, \eqref{ex:e:Gcap-3} and J$_\le (\phi_\beta)$ for $(\sE^\beta,\sF^\beta)$ in the first inequality below, \eqref{ex:e:Gcap-1} and $|\vp|\le \kappa$ in the second,   J$_\ge (\phi_\beta)$ for $(\sE^\beta,\sF^\beta)$ in the third,  and \eqref{e:sE-beta-sE-c} and \eqref{e:diffusion-part} in the fourth,   we obtain
\begin{align}\label{ex:e:Gcap-4}
	r^\beta\int_{B'}  \frac{ u(y)^2 \vp(y)^2}{d(y,o)^{\beta}} \mu(dy)
	&\le c_3 \left( 	r^\beta 	\sE^{\beta}(u^2\vp,\vp) + r^\beta \int_{B \times B}  \frac{\vp(z)\vp(y)(u(z)-u(y))^2}{V(z,d(z,y))  d(z,y)^{\beta}} \mu(dz)\mu(dy)\right)\nn\\
	&\le c_3 \left( \kappa^2 r^\beta \int_{B \times B}  \frac{(u(z)-u(y))^2}{ V(z,d(z,y)) d(z,y)^{\beta}} \mu(dz)\mu(dy) +  c_1\lVert u \rVert_2^2 \right)\\
	& \le  c_4\kappa^2 r^\beta \sE^\beta(u,u) +	c_5\lVert u \rVert_2^2\le c_6 \big( r^\alpha \sE^{(L)} (u,u) +  \lVert u \rVert_2^{2}\big).\nn
\end{align}
Therefore, since  $\vp\ge 1$ in $B'$, we arrive at
\begin{align*}
	&	\int_{B'} u(y)^2 \p_\theta(y,r)^{1+\gamma_1} \mu(dy)  \le 2^\beta \bigg( r^{\beta }	\int_{B'}  \frac{ u(y)^2 \vp(y)^2}{d(y,o)^{\beta}} \mu(dy)+ \lVert u \rVert_2^2\bigg) \le c_7 \big( r^\alpha \sE^{(L)}(u,u)+ \lVert u \rVert_2^2 \big).
\end{align*} 
The proof is complete.\qed

\begin{prop}\label{ex:Hardy-origin-2}
Suppose that $(M,d,\mu)$ satisfies 	 
	  \Exi,  and that
	  \VD \ and \RVD \ hold with $R_0=\infty$. Assume that 
	 $(\sE, \sF)$ satisfies  
  {\rm J}$_\ge (\phi_\beta)$ with $R_0=\infty$ for some $\beta\in (0,\alpha)$. 	
  Then for any $o\in M$ and $\theta \in [0, \beta \wedge d_1)$, where $d_1$ is the constant in \RVD, the function $\p_\theta(x,r;o)$ defined in \eqref{e:def-Hardy-theta} 
		is an admissible weight function for $(\sE^{(J)},	\phi_\beta)$ with $R_0=\infty$.
\end{prop}
\pf For each $x \in M$, $r \mapsto \p_\theta(x,r)$ is non-decreasing and continuous. Further, $\p_\theta$ satisfies \eqref{e:blow-up-scale} with $\gamma_0:=\beta-\theta$. Set 
$$\lambda:= \frac12( \theta + (\beta \wedge d_1)) \quad \text{and} \quad  \gamma_1:= \frac{\lambda-\theta}{\theta}.$$ 
Since $	{\rm codim}_{AI}(\{o\}) \ge d_1>\lambda$ by  Lemma \ref{l:dimensions},   \eqref{e:blow-up-integral-general} holds with $\p=\p_\theta$ and $\gamma_1$.

 Let $x \in M$, $r>0$ and  $u \in \sB_c(M)$ satisfying $\sE(u,u)<\infty$.    Write $B':=B(x,r)$ and $B:=B(x,2r)$. Let $(\sE^\lambda,\sF^\lambda)$ be the Dirichlet form defined in \eqref{e:def-sE-beta} with $\beta=\lambda$. 
    By Proposition \ref{p:subordinate-Exi},  J$_\le(\phi_\lambda)$ and J$_\ge(\phi_\lambda)$ hold for $(\sE^\lambda,\sF^\lambda)$.
     Let $R>0$ be  such that $d(z,w)>1$ for all $z \in \textrm{supp$[u]$}$ and $w \in B(x,R)^c$. By \eqref{e:VD2}, we have $V(z,d(z,y))\asymp V(y,d(y,z))$ for all $y,z \in M$.
    Using  this in the first inequality below,    and  \cite[Lemma 2.1]{CKW-memo}  and 
      J$_\ge(\phi_\beta)$  for $(\sE,\sF)$  in the second, we get
   \begin{align*}
   W_{\lambda/2}(u) &=  \int_{B(x,R) \times B(x,R)} \frac{(u(z)-u(y))^2}{V(z,d(z,y)) d(z,y)^\lambda} \mu(dz)\mu(dy) \\
   &\quad +  \int_{B(x,R)} u(z)^2 \int_{B(x,R)^c} \bigg( \frac{\mu(dy)}{V(z,d(z,y))  d(z,y)^\lambda}   + \int_{B(x,R)^c}\frac{\mu(dy)}{V(y,d(y,z))  d(y,z)^\lambda} \bigg) \, \mu(dz) \\
   &\le 	(2R)^{\beta-\lambda}\int_{B(x,R) \times B(x,R)} \frac{(u(z)-u(y))^2}{V(z,d(z,y))  d(z,y)^\beta} \mu(dz)\mu(dy)\\
   &\quad  +  c_1 \int_{\textrm{supp$[u]$}} u(z)^2 \int_{B(z,1)^c}\frac{\mu(dy)}{V(z,d(z,y))  d(z,y)^\lambda}  \mu(dz)\\
   &\le c_2 (R^{\beta-\lambda}\sE(u,u) + \lVert u \rVert_2^2)<\infty.
   \end{align*}   
   Hence,  $u \in W^{\lambda/2,2}(M) = \sF^\lambda$. 
   Since $\lambda<d_1$ and   Gcap$(\phi_\lambda)$ holds for $(\sE^\lambda,\sF^\lambda)$ by Proposition \ref{p:subordinate-Exi}, repeating the arguments for \eqref{ex:e:Gcap-4} in the proof of Proposition \ref{ex:Hardy-origin}, we obtain
    \begin{align*}
    	r^\lambda\int_{B'}  \frac{ u(y)^2 \vp(y)^2}{d(y,o)^{\lambda}} \mu(dy)&\le c_3 \left(  r^\lambda \int_{B \times B}  \frac{(u(z)-u(y))^2}{ V(z,d(z,y)) d(z,y)^{\lambda}} \mu(dz)\mu(dy) +  \lVert u \rVert_2^2 \right).
    \end{align*}
Therefore, since  $\vp\ge 1$ in $B'$, $\lambda<\beta$  and  J$_\ge (\phi_\beta)$ holds for $(\sE,\sF)$, we arrive at 
\begin{align*}
	\int_{B'} u(y)^2 \p_\theta(y,r)^{1+\gamma_1} \mu(dy) & \le 2^\lambda \bigg( r^{\lambda }	\int_{B'}  \frac{ u(y)^2 \vp(y)^2}{d(y,o)^{\lambda}} \mu(dy)+ \lVert u \rVert_2^2\bigg) \\
	&\le c_4 \left(  r^\lambda (4r)^{\beta-\lambda}\int_{B \times B}  \frac{(u(z)-u(y))^2}{ V(z,d(z,y)) d(z,y)^{\beta}} \mu(dz)\mu(dy) +  \lVert u \rVert_2^2 \right)\\& \le c_5 \big( r^\beta \sE(u,u)+ \lVert u \rVert_2^2 \big).
\end{align*} \qed

We now assume that 
 the state space is a closed subset of the Euclidean space and that the scale function is independent of the spatial variable. We  then discuss admissible weight functions that blow up to infinity at the boundary.
\begin{prop}\label{ex:Hardy-boundary}
	 Let $D$ be a non-empty $\kappa$-fat open subset of $\R^d$ satisfying
		\begin{align}\label{ex:e:AS}
			{\rm dim}_{AS}(\partial D)\le d-\eta
		\end{align}
		for some $\eta \in (0,d]$, 	and $\phi(x,r)=\phi(r)$ be a scale function independent of the spatial variable.
		Set  $M:=\overline{D}$. 
		Suppose that \Jl \ holds for a regular Dirichlet form $(\sE,\sF)$ on $L^2(M;dx)$. Then for every $\theta \in [0, \beta_2 \wedge \eta \wedge 2)$,   where $\beta_2$ is the constant in \eqref{e:phi-scale}, the function
	\begin{align*}
		\p(x,r):=\bigg(1+ \frac{\phi(r)}{\phi(\delta_D(x))}\bigg)^{\theta/\beta_2}
	\end{align*}
	is an admissible weight function for $(\sE, \phi)$ with $R_0=\diam(M)$.
\end{prop}
\pf For each $x \in M$,   $r \mapsto \p(x,r)$ is non-decreasing and continuous. Further, $\p$ satisfies \eqref{e:blow-up-scale} with $\gamma_0=(\beta_2-\theta)\beta_1/\beta_2$. Set 
$\lambda:=[\theta +  (\beta_2 \wedge \eta \wedge 2)]/2$ and $ \gamma_1:=(\lambda-\theta)/{\theta}.$
By the monotonicity of $\phi$ and \eqref{e:phi-scale},  it holds that
\begin{equation}\label{e:example-Theta-bound}
	\p(y,r)^{1+\gamma_1} \le \bigg( 2 + \1_{\delta_D(y)<r} \frac{\phi(r)}{\phi(\delta_D(y))}\bigg)^{\lambda/\beta_2} \le c_0+ c_1\bigg(\frac{r}{\delta_D(y)}\bigg)^{\lambda}  \quad \text{for all} \;\, y \in M, \, r>0.
\end{equation}
By \eqref{e:d-set}, \eqref{ex:e:AS} and Lemma \ref{l:dimensions}, we have $	{\rm codim}_{AI}(\partial D) \ge \eta>\lambda$. Thus, by \eqref{e:example-Theta-bound}, we get that for all $x\in M$ and $r \in (0,\diam(M))$,
\begin{align*}
	\frac{1}{V(x,r)}	\int_{B(x,r)} \p(y,r)^{1+\gamma_1} dy
	&\le   c_0 + \frac{c_1}{V(x,r)} \int_{B(x,r)} \bigg( \frac{r}{\delta_D(y)}\bigg)^{\lambda} dy \le c_2,
\end{align*}
proving that \eqref{e:blow-up-integral-general} holds.

 To establish \eqref{e:blow-up-Hardy},  it suffices to prove that  \eqref{e:blow-up-Hardy} holds for all $u \in C_c(M)$ satisfying $\sE(u,u)<\infty$. 
Let $x \in M$, $r>0$ and $u\in C_c(M)$ such that $\sE(u,u)<\infty$. Set $B':=B(x,r)$ and $B:=B(x,2r)$. If $\diam(M)=\infty$, then define 
  $M_\infty:=M$. If $\diam(M)<\infty$, then define  $M_\infty:=M \cup (\R^d\setminus B(0, R_1))$ for  $R_1>0$ satisfying $M \subset B(0,R_1/4)$. 
We extend $u$ to  $M_\infty$ by letting $u(z)=0$ for  $z \in M_\infty \setminus M=\R^d\setminus B(0, R_1)$.   Define a symmetric bilinear form $(\sC,\sG)$ on $L^2(M_\infty;dx)$ by
\begin{align*}
\sC(f,g)&:= \int_{M_\infty \times M_\infty} \frac{(f(z)-f(y))(g(z)-g(y))}{|z-y|^d \, \phi(|z-y|)^{\lambda/\beta_2}} dzdy, \quad 
\sC_1(f,f):=\sC_1(f,f)+ \lVert f \rVert_2^2
\\
\sG&:=\overline{\big\{f \in C_c(M_\infty):\sC(f,f)<\infty
\big\}}^{\sC_1}.
\end{align*}
By \cite[Proposition 2.2]{CK08},  $(\sC,\sG)$ is a regular Dirichlet form on  $L^2(M_\infty;dx)$.  Let $R>0$ be  such that $d(z,w)>1$ for all $z \in \textrm{supp$[u]$}$ and $w \in B(0,R)^c$.
   Using \Jl \ for $(\sE,\sF)$ and  \cite[Lemma 2.1]{CKW-memo} in the last inequality below,  we see that
\begin{align*}
\sC(u,u) &\le  \int_{B(0,R) \times B(0,R)} \frac{(u(z)-u(y))^2}{|z-y|^d \, \phi(|z-y|)^{\lambda/\beta_2}} dzdy\\
	&\quad  + 2\int_{B(0,R)} u(z)^2 \int_{B(0,R)^c} \frac{dy}{|z-y|^d \, \phi(|z-y|)^{\lambda/\beta_2}}  dz\\
	&\le 	\phi(2R)^{(\beta_2-\lambda)/\beta_2}\int_{B(0,R) \times B(0,R)} \frac{(u(z)-u(y))^2}{|z-y|^d \, \phi(|z-y|)} dzdy\\
	&\quad  +  2\int_{\textrm{supp$[u]$}} u(z)^2 \int_{B(z,1)^c} \frac{dy}{|z-y|^d \, \phi(|z-y|)^{\lambda/\beta_2}}  dz\le c_3 \sE_1(u,u)<\infty.
\end{align*}   
Hence,  $u \in 	\sG$. By \eqref{e:phi-scale}, \eqref{e:d-set} and \cite[Remarks 1.7  and  3.7]{CKW-memo}, since $\lambda<2$,   Gcap$(\phi^{\lambda/\beta_2})$  holds  for $(	\sC,\sG)$. Thus, there exist   $\kappa \ge 1$, $c_4>0$ and  a $\kappa$-cutoff function $\vp \in  \sG$ for $B' \Subset B$ such that
\begin{align}\label{ex:e:Gcap-distance}
\sC(u^2\vp,\vp) \le \frac{c_4}{\phi(r)^{\lambda/\beta_2}}\int_{B} u^2 dz.
\end{align}
Since $	{\rm dim}_{AS}(\partial D)< d-\lambda$, using \cite[Theorem 1]{DV14}  and  \cite[(3.18)]{GHH18}, we get
\begin{align*}
\int_{M_\infty}  \frac{ u(y)^2 \vp(y)^2}{\phi(\delta_D(y))^{\lambda/\beta_2}} dy	&\le c_5 \sC(u\vp,u\vp)= c_5\bigg( \sC(u^2\vp,\vp) + \int_{B \times B}  \frac{\vp(z)\vp(y)(u(z)-u(y))^2}{|z-y|^d \phi(|z-y|)^{\lambda/\beta_2}}dzdy \bigg).
\end{align*}  Using this in the second inequality below,  \eqref{ex:e:Gcap-distance} in the third and  \eqref{e:phi-scale} in the fifth, since $\vp\ge 1$ in $B'$ and   $0\le \vp\le \kappa$, we arrive at
\begin{align*}
	&\int_{B'} u(y)^2 \p(y,r)^{1+\gamma_1} dy 
	\le c_6 \bigg( \phi(r)^{\lambda/\beta_2}	\int_{M_\infty}  \frac{ u(y)^2 \vp(y)^2}{\phi(\delta_D(y))^{\lambda/\beta_2}} dy+ \lVert u \rVert_2^2\bigg) \\
	&\le c_7 \bigg( \phi(r)^{\lambda/\beta_2} \sC(u^2\vp,\vp)+ \kappa^2\phi(r)^{\lambda/\beta_2} \int_{B \times B}\frac{(u(z)-u(y))^2}{|z-y|^d \phi(|z-y|)^{\lambda/\beta_2}}dzdy + \lVert u \rVert_2^2\bigg) \\
	&\le  c_8 \left( \phi(r)^{\lambda/\beta_2} \int_{B \times B}\frac{(u(z)-u(y))^2}{|z-y|^d \phi(|z-y|)^{\lambda/\beta_2}}dzdy+ \lVert u \rVert_2^2\right)\\
		&\le  c_8 \left( \phi(r)^{\lambda/\beta_2} \phi(4r)^{(\beta_2-\lambda)/\beta_2} \int_{B \times B}\frac{(u(z)-u(y))^2}{|z-y|^d \phi(|z-y|)}dzdy+ \lVert u \rVert_2^2\right)\\
			&\le  c_8 \left( c_9\phi(r) \int_{B \times B}\frac{(u(z)-u(y))^2}{|z-y|^d \phi(|z-y|)}dzdy+ \lVert u \rVert_2^2\right).
\end{align*} 
Since \Jl \ holds for $(\sE,\sF)$, this yields that \eqref{e:blow-up-Hardy} holds. The proof is complete. \qed

\section{Examples of weighted  cutoff Sobolev inequalities}\label{s:CS}

In this section, we provide several examples where the weighted version \CS \  of  a cutoff Sobolev inequality holds for an admissible weight function $\p$.

For  $D \subset M$, let ${\rm Lip_c}(D)$ denote the family of all Lipschitz continuous functions with compact support in $D$.

\begin{prop}\label{ex:beta<2}
	Suppose that  $\sE^{(L)}=0$,    $\sF$ contains all functions $f \in C_c(M)$ with $\sE(f,f)<\infty$, that $\phi$  satisfies \eqref{e:phi-scale} with $\beta_2<2$ and that \Tail \ holds.   Then    ${\rm Lip_c}(M) \subset \sF$ and  \CSU \ holds for $(\sE,\sF)$.
\end{prop}
\pf Let $u \in {\rm Lip_c}(M)$. 
Fix $x_0\in M$ and $R\ge1$ such that supp$[u] \subset B(x_0,R)$.  There exists $c_1>0$ such that 
\begin{align*}
	|u(x)-u(y)|^2\le c_1\left(1 \wedge \frac{d(x,y)}{R}\right)^2 \quad \text{for all} \;\, x,y \in M.
\end{align*}
Hence, since $\beta_2<2$, by  Lemma \ref{l:Tail-integral}, we get  that for all $x \in M$,
\begin{align}\label{e:example-CS-01}
	\int_M (u(x)-u(y))^2 J(x,dy) \le  \frac{c_2\p(x,R)}{\phi(x,R)}.
\end{align}
Using   \eqref{e:example-CS-01} in the second  inequality below,  \eqref{e:phi-comp} in the third and \eqref{e:blow-up-integral} in the fourth, we obtain
\begin{align}\label{e:example-CS-02}
	\sE(u,u)& \le  2\int_{B(x_0,R)} \int_M (u(x)-u(y))^2 J(x,dy) \mu(dx) \le 2c_2 \int_{B(x_0,R)} \frac{\p(x,R)}{\phi(x,R)} \mu(dx)\nn\\
	&\le  \frac{c_3}{\phi(x_0,R)} \int_{B(x_0,R)} \p(x,R) \mu(dx) \le\frac{c_4 V(x_0,R)}{\phi(x_0,R)} <\infty.
\end{align}
Therefore, $u \in \sF$, proving that  ${\rm Lip_c}(M) \subset \sF$.

 We now prove that \CSU \ holds. Let $x_0 \in M$ and $0<r \le R$ with $R+2r<R_0$.  Set $V_0:=B(x_0,R)$,  $V_1:=B(x_0,R+r)$ and $V_2:=B(x_0,R+2r)$.  Define
\begin{align}\label{e:Lipschitz-cutoff}
\vp(x):=1 \wedge \frac{(R+r-d(x,x_0))_+}{r}.
\end{align}
Then $\vp \in C_c(M)$ and $\vp$ is a cutoff function for $V_0 \Subset V_1$. Further, $|\vp(x)-\vp(y)|^2 \le (1 \wedge (d(x,y)/r))^2$ for all $x,y \in M$. Hence, since $\beta_2<2$, using Lemma \ref{l:Tail-integral}, we get
\begin{align}\label{e:example-CS-1}
	\int_M (\vp(x)-\vp(y))^2 J(x,dy) \le \frac{c_1\p(x,r)}{\phi(x,r)} \quad \text{for all} \;\, x\in M.
\end{align}
Using \eqref{e:example-CS-1} and repeating the argument for \eqref{e:example-CS-02}, we can deduce that  $\sE(\vp,\vp)<\infty$. Thus, $\vp \in \sF_b$.  Moreover, using \eqref{e:example-CS-1}, we get that for all $f \in \sF'_b$,
\begin{align*}
	\int_{V_2} f^2\, d\Gamma(\varphi,\varphi) &= \int_{V_2} \int_M f(x)^2 (\vp(x)-\vp(y))^2 J(x,dy)\mu(dx)\le \sup_{z \in V_2}\frac{c_1}{\phi(z,r)} \int_{V_2} f(x)^2 \p(x,r)\mu(dx).
\end{align*}
The proof is complete. \qed

Recall that under \Exi,  the regular Dirichlet form 
$(\sE^\beta,\sF^\beta)$ is defined as in \eqref{e:def-sE-beta} for  $\beta \in (0,\alpha)$.
For the next result, we will use Lemma \ref{l:CS-continuity} which will be proved in Section \ref{s:proof-thm-1}.

\begin{prop}\label{ex:CS-Hardy}
Let $o \in M$. Suppose that \VD \ and \RVD \ hold with $R_0=\infty$, and \Exi \ holds. 	 Suppose also that  $\sE^{(L)}=0$,  and   {\rm J}$_\ge(\phi_\beta)$ and  {\rm J}$^{\p}_\le(\phi_\beta)$  hold with  $R_0=\infty$,  $\beta  \in (0, \alpha \wedge d_1)$,  
 where $d_1$ is the constant in \RVD,     and   $\p$ satisfying
		\begin{align}\label{e:CS-Hardy}
				\p(x,r)\le \bigg(1+ \frac{r}{d(x,o)}\bigg)^{\theta_0}, \quad \;\; x\in M, \;r>0,
		\end{align} 
		for some 		$\theta_0 \in [0, \beta)$.    Then
		\begin{align}\label{e:CS-Hardy-domain}
			 \sE(u,u) \asymp \sE^\beta(u,u) \quad \text{for all $u \in \sF^\beta$},
		\end{align}
		and	 there exists 
		$\theta \in  [\theta_0, \beta)$ such that  {\rm J}$^{\p_\theta}_\le(\phi_\beta)$  and  {\rm CS}$^{\p_\theta}(\phi_\beta)$ hold  with $R_0=\infty$ and the  function  $\p_\theta$ defined in \eqref{e:def-Hardy-theta}.
\end{prop}
\pf   We first prove \eqref{e:CS-Hardy-domain}. Since  {\rm J}$_\ge(\phi_\beta)$ holds for $\sE$ and  {\rm J}$_\le(\phi_\beta)$ holds for $\sE^\beta$, we have $\sE(u,u) \ge c_1\sE^\beta(u,u)$ for all  $u \in \sF^\beta$. 
Define
\begin{align*}
	E&=\left\{ (x,y)\in M\times M: d(x,y)\ge 2(d(x,o)\wedge d(y,o)) \right\}.
\end{align*} 
By  {\rm J}$^{\p}_\le(\phi_\beta)$ and \eqref{e:CS-Hardy},   we have $J(x,y) \le  c_2V(x,d(x,y))^{-1}d(x,y)^{-\beta}$ for $(x,y) \in E^c$. Using this and the symmetry of $J$, we see that for all $u \in \sF^\beta$,
\begin{align}\label{e:CS-Hardy-domain-1}
	\sE(u,u)&\le  c_3 \sE^\beta(u,u) +  4\int_{E} u(x)^2 J(x,y)\mu(dx)\mu(dy).
\end{align} 
Observe that
\begin{align}\label{e:CS-Hardy-range}
	d(x,y) \ge d(x,o) \quad \text{for all $(x,y)\in E$}.
\end{align}
Indeed, for any $(x,y)\in E$, if $d(x,o)\ge 2d(y,o)$, then $d(x,y) \ge d(x,o)-d(y,o) \ge d(x,o)$ and if $d(x,o)<2d(y,o)$, then $d(x,y) \ge 2(d(x,o) \wedge (d(x,o)/2)) = d(x,o)$.  Using \eqref{e:CS-Hardy-range},   Lemma \ref{l:Ju-Tail} and \eqref{e:CS-Hardy}, we obtain
\begin{align}\label{e:CS-Hardy-domain-2}
&\int_{E} u(x)^2 J(x,y)\mu(dx)\mu(dy)\le \int_M u(x)^2 \int_{B(x,d(x,o))^c} J(x,y)\mu(dy) \mu(dx) \nn\\
&\le c_4\int_M \frac{u(x)^2 \p(x,d(x,o))}{d(x,o)^{\beta}} \mu(dx) \le 2^{\theta_0}c_4 \int_M \frac{u(x)^2}{d(x,o)^{\beta}} \mu(dx).
\end{align}
Recall from the proof of Proposition \ref{ex:Hardy-origin} that condition $\mathbf{(G)_{\beta}}$ in \cite{CGL21} holds for $(\sE^\beta,\sF^\beta)$.  By   \cite[Theorem 5.6]{CGL21}, we have $\int_M u(x)^2d(x,o)^{-\beta} \mu(dx)\le c_5\sE^\beta(u,u)$. Combining this with \eqref{e:CS-Hardy-domain-1} and  \eqref{e:CS-Hardy-domain-2}, we conclude  \eqref{e:CS-Hardy-domain}. 

\bk

 By Proposition \ref{p:subordinate-Exi},   $(\sE^\beta, \sF^\beta)$ satisfies condition HK$(\phi_\beta)$ in \cite{CKW-jems}. 
 Hence, by \cite[Theorem 1.18 and Corollary 1.3]{CKW-jems},  E$(\phi_\beta)$ and   PHR$(\phi_\beta)$  hold for $(\sE^\beta, \sF^\beta)$. Thus, by Lemma \ref{l:CS-continuity},  CSU$(\phi_\beta)$ holds for $(\sE^\beta, \sF^\beta)$ with  cutoff functions satisfying \eqref{e:CS-continuity} for some $\chi \in (0,1]$. Let 
 $
 \theta:=\theta_0 \vee (\beta -\chi).$
By  \eqref{e:CS-Hardy},  $\p(x,r)\le \p_\theta(x,r)$ for all $x \in M$ and $r>0$. Hence, J$_\le^{\p_\theta}(\phi_\beta)$  holds for $(\sE,\sF)$.

 We now establish  CS$^{\p_\theta}(\phi_\beta)$ for $(\sE,\sF)$. Let $x_0 \in M$, $0<r\le R$ and $f \in \sF'_b$. Write $B_s:=B(x_0,s)$ for $s>0$.  Let $a>0$ be such that $f-a\in \sF$.
  By J$_\ge(\phi_\beta)$, we have 
 $W_{\beta/2}(f-a)\le c_6 \sE(f,f)<\infty
 $. Hence,  $ f-a\in W^{\beta/2,2}(M)=\sF^\beta$ so that $f \in (\sF^\beta)'_b$.   We consider the following two cases separately.\bk

\smallskip

Case 1: $d(x_0,o)\ge R+r/8$. Let $\vp\in \sF^\beta_b$ be a cutoff function for $B_R \Subset B_{R+r/32}$ satisfying   \eqref{e:CS-continuity} and CSU$(\phi_\beta)$ for $(\sE^\beta, \sF^\beta)$ with  $x_0,R$ and $r$ replaced by $r/32$. By Lemma \ref{l:Ju-Tail}, TJ$_\le^{\p_\theta}(\phi_\beta)$ holds for $(\sE, \sF)$. Hence, we have
\begin{align}\label{e:CS-Hardy-case1}
	&\int_{B_{R+2r} \times B(x,r/32)^c} f(x)^2 (\vp(x)-\vp(y))^2 J(x,y) \mu(dx)\mu(dy)\nn\\
	& \le \int_{B_{R+2r}} f(x)^2 J(x, B(x,r/32)^c)  \mu(dx)\le \frac{c_7}{r^\beta} \int_{B_{R+2r}} f(x)^2 \p_\theta(x,r) \mu(dx).
\end{align}
Further, since $\vp = 0$ in $B_{R+r/32}^c$, it holds that
\begin{align}\label{e:CS-Hardy-case12}
	&\int_{(B_{R+2r} \setminus B_{R+r/16}) \times B(x,r/32)} f(x)^2 (\vp(x)-\vp(y))^2 J(x,y) \mu(dx)\mu(dy)=0.
\end{align}
Note that $d(x,o)\wedge d(y,o) > r/32$ for all $(x,y) \in B_{R+r/16} \times B(x,r/32)$ in this case. Hence, by J$^{\p}_\le(\phi_\beta)$ and \eqref{e:CS-Hardy}, we get that for all  $(x,y) \in B_{R+r/16} \times B(x,r/32)$,
\begin{align*}
	J(x,y) \le  \frac{c_8}{ V(x,d(x,y)) d(x,y)^{  \beta}} \bigg(1+ \frac{r/32}{d(x,o)}\bigg)^{\theta_0} \bigg(1+ \frac{r/32}{d(y,o)}\bigg)^{\theta_0} \le  \frac{2^{2\theta_0}c_8}{ V(x,d(x,y)) d(x,y)^{  \beta}}.
\end{align*}
Using this and  CSU$(\phi_\beta)$ for $(\sE^\beta, \sF^\beta)$, since J$_\ge (\phi_\beta)$ and J$_\le (\phi_\beta)$ hold for $(\sE^\beta,\sF^\beta)$,  we obtain
\begin{align}\label{e:CS-Hardy-case13}
&	\int_{B_{R+r/16}  \times B(x,r/32)} f(x)^2 (\vp(x)-\vp(y))^2 J(x,y) \mu(dx)\mu(dy)\nn\\
			&\le  2^{\theta_0} c_8	\int_{B_{R+r/16}  \times M}  \frac{ f(x)^2 (\vp(x)-\vp(y))^2}{ V(x,d(x,y)) d(x,y)^{  \beta}}\mu(dx)\mu(dy)\nn\\
			&\le   c_9 \int_{B_{R+r} \times B_{R+2r}}  \frac{ \vp(x)^2 (f(x) - f(y))^2 }{  V(x,d(x,y)) d(x,y)^{  \beta}} \mu(dx)\mu(dy)  +  \frac{c_{10}}{r^\beta} \int_{B_{R+2r}} f(x)^2 \mu(dx). 
\end{align}
Combining  \eqref{e:CS-Hardy-case1}, \eqref{e:CS-Hardy-case12} and  \eqref{e:CS-Hardy-case13}, since 
J$_\ge(\phi_\beta)$ holds for $(\sE,\sF)$, we are done in this case.

Case 2: $d(x_0,o)<R+r/8$. Let $\vp\in \sF^\beta_b$ be a cutoff function for $B_R \Subset B_{R+r}$ satisfying   \eqref{e:CS-continuity} and CSU$(\phi_\beta)$ for $(\sE^\beta, \sF^\beta)$.  Note that  \eqref{e:CS-Hardy-case1} remains valid. Since $d(x_0,o)<R+r/8$, by 
J$_\le ^{\p}(\phi_\beta)$,  \eqref{e:CS-Hardy} and CSU$(\phi_\beta)$ for $(\sE^\beta, \sF^\beta)$, we have
\begin{align}\label{e:CS-Hardy-case2-1}
	&	\int_{(B_{R+2r}\setminus B_{R+r/4}) \times B(x,r/32)} f(x)^2 (\vp(x)-\vp(y))^2 J(x,y) \mu(dx)\mu(dy)\nn\\
	&\le c_{11} 	\int_{(B_{R+2r}\setminus B_{R+r/4}) \times B(x,r/32)} \frac{ f(x)^2 (\vp(x)-\vp(y))^2}{  V(x,d(x,y)) d(x,y)^{  \beta}}  \bigg(1+ \frac{r/32}{d(x,o) }\bigg)^{\theta_0}  \bigg(1+ \frac{r/32}{ d(y,o)}\bigg)^{\theta_0} \mu(dx)\mu(dy)\nn\\
		&\le (4/3)^{2\theta_0} c_{11}	\int_{(B_{R+2r}\setminus B_{R+r/4}) \times B(x,r/32)} \frac{ f(x)^2 (\vp(x)-\vp(y))^2}{  V(x,d(x,y)) d(x,y)^{  \beta}}   \mu(dx)\mu(dy)\nn\\
			&\le   c_{12} \int_{B_{R+r} \times B_{R+2r}}  \frac{ \vp(x)^2 (f(x) - f(y))^2 }{  V(x,d(x,y)) d(x,y)^{  \beta}} \mu(dx)\mu(dy)  +  \frac{c_{13}}{r^\beta} \int_{B_{R+2r}} f(x)^2 \mu(dx). 
\end{align}
Besides,  since $\theta_0\le \theta$, using J$^{\p}_\le (\phi_\beta)$,   \eqref{e:CS-Hardy}, \eqref{e:CS-continuity}  and CSU$(\phi_\beta)$ for $(\sE^\beta, \sF^\beta)$, we see that
\begin{align}\label{e:CS-Hardy-case2-2}
	&	\int_{B_{R+r/4} \times B(x,r/32)} f(x)^2 (\vp(x)-\vp(y))^2 J(x,y) \mu(dx)\mu(dy)\nn\\
	&\le c_{14}	\int_{B_{R+r/4}  \times B(x,r/32)} \frac{ f(x)^2 (\vp(x)-\vp(y))^2}{  V(x,d(x,y)) d(x,y)^{  \beta}} \mu(dx)\mu(dy)\nn\\
	&\quad  + c_{14}	\int_{B_{R+r/4}  \times B(x,r/32)} \frac{ f(x)^2 (\vp(x)-\vp(y))^2}{  V(x,d(x,y)) d(x,y)^{\beta}}\bigg( \frac{d(x,y)}{d(x,o)}\bigg)^{\theta}\bigg( \frac{d(x,y)}{d(y,o)}\bigg)^{\theta} \mu(dx)\mu(dy)\nn\\
		&\le   c_{15} \int_{B_{R+r} \times B_{R+2r}}  \frac{ \vp(x)^2 (f(x) - f(y))^2 }{  V(x,d(x,y)) d(x,y)^{  \beta}} \mu(dx)\mu(dy)  +  \frac{c_{16}}{r^\beta} \int_{B_{R+2r}} f(x)^2 \mu(dx)\nn\\
	&\quad  + \frac{c_{17}}{r^{2\chi}}	\int_{B_{R+r/4}  \times B(x,r/32)} \frac{ f(x)^2 d(x,y)^{2\chi}}{  V(x,d(x,y)) d(x,y)^{\beta}}\bigg( \frac{d(x,y)}{d(x,o)}\bigg)^{\theta}\bigg( \frac{d(x,y)}{d(y,o)}\bigg)^{\theta} \mu(dx)\mu(dy).
\end{align}
 Recall that ${\rm codim}_{AI}(\{o\}) \ge d_1$ by Lemma \ref{l:dimensions}. Hence, since $\theta<d_1$ and $2\chi + \theta>\beta$, using \VD, we deduce that  for any $x \in B_{R+r/4}$, 
\begin{align*}
	& \frac{1}{r^{2\chi}}	\int_{ B(x,r/32)} \frac{  d(x,y)^{2\chi+2\theta}}{  V(x,d(x,y)) d(x,y)^{\beta} d(y,o)^\theta} \mu(dy)\nn\\
	&\le r^{\theta -\beta}\sum_{n\ge 5} \frac{ 2^{-n(2\chi+2\theta)+(n+1)\beta} }{  V(x,2^{-n-1}r) } 	\int_{ B(x,2^{-n}r)\setminus B(x,2^{-n-1}r)} \bigg(\frac{r}{  d(y,o)}\bigg)^\theta \mu(dy)\nn\\
		& \le c_{18}r^{\theta -\beta} \sum_{n\ge 5} \frac{ 2^{-n(2\chi+\theta)+(n+1)\beta} V(x,2^{-n}r) }{  V(x,2^{-n-1}r) } \le c_{19}r^{\theta -\beta} \sum_{n\ge 5} 2^{-n(2\chi+\theta-\beta)}  = c_{c_{20}}r^{\theta-\beta}.
\end{align*}
Using this, we obtain
\begin{align}	\label{e:CS-Hardy-case2-3}
		&\frac{1}{r^{2\chi}}	\int_{B_{R+r/4}  \times B(x,r/32)} \frac{ f(x)^2 d(x,y)^{2\chi}}{  V(x,d(x,y)) d(x,y)^{\beta}}\bigg( \frac{d(x,y)}{d(x,o)}\bigg)^{\theta}\bigg( \frac{d(x,y)}{d(y,o)}\bigg)^{\theta} \mu(dx)\mu(dy)\nn\\
			&\le \frac{c_{21}}{r^\beta} \int_{B_{R+r/4}} f(x)^2  \bigg(\frac{r}{d(x,o)}\bigg)^\theta \mu(dx)\le \frac{c_{21}}{r^\beta} \int_{B_{R+r/4}} f(x)^2  \p_\theta(x,r) \mu(dx).\end{align}
Combining   \eqref{e:CS-Hardy-case1},  \eqref{e:CS-Hardy-case2-1},   \eqref{e:CS-Hardy-case2-2} and \eqref{e:CS-Hardy-case2-3}, since  J$_\ge(\phi_\beta)$ holds for $(\sE,\sF)$ and $\p_\theta \ge 1$, we get the desired result.

\smallskip

The proof is complete. \qed

\section{Consequences of \Tail \ and \CS} \label{s:CTCS}

In this section, we discuss several consequences of \Tail \ and \CS  \ including self-improvement properties of \CS, which will be used later sections.

\begin{lem}\label{l:Tail}
		Suppose that \Tail \ holds. There exists $C>0$ such that
	\begin{align*}
		J\left(B(x_0,R)\times B(x_0,R+r)^c\right) 	 \le C\frac{V(x_0,R)}{\phi(x_0,R)} \bigg(\frac{R}{r}\bigg)^{\beta_2} \quad \text{for all} \;\, x_0 \in M, \; 0<r\le R<R_0.
	\end{align*}
\end{lem}
\pf  
Using \Tail, \eqref{e:phi2} and \eqref{e:blow-up-integral}, we get that for all $x_0\in M$ and $0<r\le R<R_0$, 
\begin{align*}
		&	J\left(B(x_0,R)\times B(x_0,R+r)^c\right)   \le 	\int_{B(x_0,R)} 	J(x,B(x,r)^c) \mu(dx)\\
		& \le c_1\int_{B(x_0,R)} \frac{\p(x,r)}{\phi(x,r)}\mu(dx) \le \frac{c_2}{\phi(x_0,R)} \bigg(\frac{R}{r}\bigg)^{\beta_2}\int_{B(x_0,R)} \p(x,r)\mu(dx)\le \frac{c_3V(x_0,R)}{\phi(x_0,R)} \bigg(\frac{R}{r}\bigg)^{\beta_2}.
\end{align*}
\qed

In the next lemma, we prove the self-improvement property of the constant $C_0$ in \CS. In the case of  $\p=1$, this property was obtained in \cite[Proposition 2.3(4)]{CKW-memo}.
\begin{lem}\label{l:CS-basic}
	Suppose that \VD \ and \Tail \ hold. 	If \CS \ holds for some $C_0 \in (0,1]$, then for any $C_0' \in [C_0,1]$, \CS \ remains valid for $C_0'$ with a redefined $C_2 > 0$.
\end{lem}
\pf  Let $x_0 \in M$ and write $B_s:=B(x_0,s)$ for $s>0$. Assume that \eqref{e:CS} holds for $C_0 \in (0,1]$ with a cutoff function $\vp \in \sF_b$ for $B_R \Subset B_{R+r}$. Let $C_0' \in (C_0,1]$ and denote $U:=B_{R+(1+C_0')r} \setminus B_{R+(1+C_0)r}$. Since $\vp=0$ in $U$,  we have $\int_{U} f^2 d \Gamma^{(L)} (\vp,\vp)=0$ by the strong local property. Thus, by \Tail, \eqref{e:phi-scale} and the monotonicity of $\p$, we get
\begin{align*}
	&\int_{U} f^2 d \Gamma (\vp,\vp)  =\int_U f(x)^2 \int_{B_{R+r}} \vp(y)^2 J(x,dy) \mu(dx)\le \int_U f(x)^2 J(x,B(x,C_0r)^c) \mu(dx)\\
		& \le  c_1 \int_U f(x)^2 \frac{\p(x,C_0r)}{\phi(x,C_0r)} \mu(dx)\le  \sup_{z \in B(x_0, R +(1+C_0')r)} \frac{c_2}{\phi(z,r)} \int_U f(x)^2 \p(x,r) \mu(dx).
\end{align*}
Therefore,  \CS  \ holds for $C_0'$. \qed

Henceforth, in view of Lemma \ref{l:CS-basic}, we assume that \CS \ holds with $C_0=1$
whenever \VD, \Tail \ and \CS \ hold.

We next show the  self-improvement property of the leading constant $C_1$ in  \CS. 
 The  proof of the next lemma adopts the method from \cite[Lemma 2.9]{GHH18} and \cite[Proposition 2.4]{CKW-memo}, initially established in \cite[Lemma 5.1]{AB15} within the context of local Dirichlet forms.   

\begin{lem}\label{l:CS-improv}
	Suppose that \VD,  \Tail \ and \CS \ hold. Then for any $\eps>0$, there exists a constant $C(\eps)>0$ such that  for any $x_0 \in M$,  $0<r\le R$ with $R+2r<R_0$ and 	$f\in \sF'_b$, there exists a cutoff function $\vp \in \sF_b$ for $B(x_0,R) \Subset B(x_0,R+r)$ satisfying the following inequality:
		\begin{align}\label{e:CS-improv}
		&\int_{B(x_0,R+2r)} f^2\, d\Gamma(\varphi,\varphi) \nn\\
		&\le   \eps \Big(\int_{B(x_0,R+r)} \varphi^2 d\Gamma^{(L)}(f,f) + \int_{B(x_0,R+r) \times B(x_0,R+2r)} \varphi(x)^2 (f(x) - f(y))^2 J(dx,dy) \Big) \nn\\
		&\quad + \sup_{z \in B(x_0,R+2r)}  \frac{C(\eps)}{\phi(z,r)} \int_{B(x_0,R+2r)} f(x)^2 \p(x,r) \mu(dx). 
	\end{align}	
\end{lem}
\pf Let $x_0 \in M$,  $0<r\le R$ with $R+2r<R_0$ and $f \in \sF'_b$. Write $B_s:=B(x_0,s)$ for $s>0$.  Without loss of generality, we assume $\int_{B_{R+2r}} f(x)^2 \p(x,r)\mu(dx) \in (0,\infty)$. Set
$$
\theta:=\Big(\frac{1}{V(x_0,R+2r)}\int_{B_{R+2r}}  f(x)^2 \p(x,r) \mu(dx)  \Big) ^{1/2}  \quad \text{and} \quad g:=|f|+\theta.
$$
  For  a constant $\lambda>0$  whose exact value is to be determined later, set $r_0:=0$,
\begin{align*}
s_n:=c_1 r e^{-\lambda n/\beta_2} \quad \text{ and } \quad    r_n:=\sum_{k=1}^n s_k \quad  \text{for}\;\;n \ge 1,
\end{align*}
where $\beta_2$ is the constant in \eqref{e:phi-scale} and  $c_1:=(\sum_{n=0}^\infty e^{-\lambda n/\beta_2})^{-1}$. Note that $|g(x)-g(y)| \le |f(x)-f(y)|$ for all $x,y \in B_{R+2r}$, and by \cite[Proposition 4.3.1, Theorems 4.3.7 and 4.3.8]{CF},  $d\Gamma^{(L)}(g,g)=d\Gamma^{(L)}(f,f)$ on $B_{R+2r}$.  Thus, for each $n \ge 1$, by \CS \ (with $C_0=1$), there exists a cutoff function $\vp_n\in \sF_b$ for $B_{R+r_n} \Subset B_{R+r_{n+1}}$ so that 
	\begin{align}\label{e:CS-1}
	&\int_{B_{R+r_{n+1}+s_{n+1}}} g^2\, d\Gamma(\varphi_n,\varphi_n) \nn\\
	&\le   c_1\int_{B_{R+r_{n+1} }} \varphi_n^2\, d\Gamma^{(L)}(f,f) +c_1 \int_{B_{R+r_{n+1} } \times B_{R+r_{n+1}+s_{n+1}}} \varphi_n(x)^2 (f(x) - f(y))^2 J(dx,dy)  \nn\\
	&\quad + \sup_{z \in B_{R+2r}} \frac{c_2}{\phi(z,s_{n+1})} \int_{B_{R+r_{n+1}+s_{n+1}}} g(x)^2 \p(x,s_{n+1})  \mu(dx). 
\end{align}	
Define
\begin{align*}
	\vp = \sum_{n=1}^\infty \, (e^{-\lambda(n-1)}-e^{-\lambda n})\vp_n.
\end{align*}
We will show that $\vp \in \sF_b$ in the last paragraph of this proof. Since each $\vp_n$ is a cutoff function for $B_{R+r_n} \Subset B_{R+r_{n+1}}$,  $\vp\in \sF_b$ is a cutoff function for $B_R \Subset B_{R+r}$ and for all $k \ge 1$,
\begin{align}\label{e:vp>ek}
	\vp \ge \sum_{n=k}^\infty(e^{-\lambda(n-1)}-e^{-\lambda n})=  e^{-\lambda (k-1)} \quad  \text{on}\;\, B_{R+r_k}.
\end{align}
Define $\Phi_{n,m}(x,y)=(\vp_n(x)-\vp_n(y))(\vp_m(x)-\vp_m(y))$ for $n,m \ge 1$. We have
\begin{align}\label{e:CS-improv-1}
	\int_{B_{R+2r}} g^2 \,d \Gamma^{(J)}(\vp,\vp)	&=\int_{B_{R+2r} \times M} g(x)^2  \bigg( \sum_{n=1}^\infty (e^{-\lambda(n-1)}-e^{-\lambda n})(\vp_n(x)-\vp_n(y))\bigg)^2 J(dx,dy)\nn\\
	&\le2 (e^\lambda -1)^2\sum_{n=3}^\infty \sum_{m=1}^{n-2}\int_{B_{R+2r}\times M} e^{-\lambda (n+m)} g(x)^2\Phi_{n,m}(x,y)\, J(dx,dy)\nn\\
		&\quad  + 2(e^\lambda -1)^2 \sum_{n=2}^\infty e^{-\lambda(2n-1)} \int_{B_{R+2r}\times M}  g(x)^2\Phi_{n,n-1}(x,y)\, J(dx,dy)\nn\\
		&\quad+ (e^\lambda-1)^2\sum_{n=1}^\infty  e^{-2\lambda n}  \int_{B_{R+2r}\times M}   g(x)^2\Phi_{n,n}(x,y) \, J(dx,dy)\nn\\
	&=:I_1+I_2+I_3.
\end{align}

Let $m \ge 1$ and $n \ge m+2$. Note that $\Phi_{n,m} \le 1$. Further,  $\Phi_{n,m}(x,y) \neq 0$ only if $x \in B_{R+r_{m+1}}$ and $y \in B_{R+r_n}^c$, or $x \in B_{R+r_n}^c$ and $y \in B_{R+r_{m+1}}$. Using these facts  in the first line below, the fact that $r_n - r_{m+1} \ge s_{m+2}$ in the second, and \Tail \ in the last, we get
\begin{align*}
	&\int_{B_{R+2r}\times M}  g(x)^2\Phi_{n,m}(x,y)J(dx,dy)\nn\\
	 &\le  \int_{B_{r_{m+1}}} g(x)^2 \int_{B_{R+r_n}^c} J(x,dy) \mu(dx) + \int_{B_{R+2r} \setminus B_{R+r_n}} g(x)^2 \int_{B_{R+r_{m+1}}}  J(x,dy)  \mu(dx)\\
	 &\le   \int_{B_{r_{m+1}}} g(x)^2  J(x,B(x, s_{m+2})^c) \mu(dx) + \int_{B_{R+2r} \setminus B_{R+r_n}} g(x)^2 J(x,B(x, s_{m+2})^c) \mu(dx) \\
	 &\le c_3\int_{B_{R+2r}} g(x)^2 \frac{\p(x,s_{m+2}) }{\phi(x,s_{m+2})} \mu(dx).
\end{align*}
Using this in the first inequality below,  \eqref{e:phi-scale} and the monotonicity of $\p$ in the second, and the definition $s_n=c_1 r e^{-n\lambda/(2\beta_2)}$ in the third, we obtain
\begin{align}\label{e:CS-improv-2}
	I_1	&\le 2c_3(e^{\lambda}-1)^2 \sum_{n=3}^\infty \sum_{m=1}^{n-2} e^{-\lambda(n+m)} \int_{B_{R+2r}}  g(x)^2 \frac{\p(x,s_{m+2}) }{\phi(x,s_{m+2})} \mu(dx)\nn\\
	&\le \sup_{z\in B_{R+2r}} \frac{ c_4(e^\lambda -1)^2}{\phi(z,r)}\sum_{n=3}^\infty \sum_{m=1}^{n-2} e^{-\lambda (n+m)} \bigg( \frac{r}{s_{m+2}}\bigg)^{\beta_2}  \int_{B_{R+2r}}  g(x)^2 \p(x,r) \mu(dx) \nn\\
	&\le  \sup_{z \in B_{R+2r}} \frac{c_5( e^{\lambda}-1)^2 }{\phi(z,r)} \bigg(  \int_{B_{R+2r}}  g(x)^2 \p(x,r)  \mu(dx)\bigg)  \sum_{n=3}^\infty \sum_{m=1}^{n-2} e^{-\lambda(n+m/2)} \nn\\
	& \le\sup_{z \in B_{R+2r}} \frac{c_6(\lambda)}{\phi(z,r)}  \int_{B_{R+2r}}  g(x)^2 \p(x,r)  \mu(dx),
\end{align}
where $c_6(\lambda)$ is a constant depending on $\lambda$.  

By the Cauchy-Schwarz inequality, we have $I_2 \le 2I_3$. Observe that
\begin{align*}
	&\int_{B_{R+r_{n+1}+s_{n+1}} \times M}   g(x)^2 \Phi_{n,n}(x,y) J(dx,dy)\le \int_{B_{R+r_{n+1}+s_{n+1}}} g^2\, d\Gamma^{(J)}(\vp_n,\vp_n).
\end{align*}
Further, since $\Phi_{n,n} \le 1$ and $\Phi_{n,n}(x,y)=0$ for $x,y \in B_{R+r_{n+1}}^c$, by \Tail,  we have
\begin{align*}
	&\int_{(B_{R+2r} \setminus B_{R+r_{n+1}+s_{n+1}})\times M}  g(x)^2 \Phi_{n,n}(x,y)J(dx,dy)\\
	&\le
	\int_{(B_{R+2r} \setminus B_{R+r_{n+1}+s_{n+1}}) }  g(x)^2 \int_{B(x,s_{n+1})^c}J(x,dy) \mu(dx)\\
	& \le \sup_{z \in B_{R+2r}} \frac{c_7}{\phi(z,s_{n+1})} \int_{B_{R+2r}}  g(x)^2 \p(x,s_{n+1})  \mu(dx).
\end{align*}
Thus, by \eqref{e:CS-1}, we have
\begin{align}\label{e:CS-improv-4}
	I_3 &\le c_8(e^\lambda-1)^2 \sum_{n=1}^\infty e^{-2\lambda n}  \int_{B_{R+r_{n+1} }} \vp_n^2 d\Gamma^{(L)}(f,f)  \nn\\
	&\quad + c_8(e^\lambda-1)^2 \sum_{n=1}^\infty e^{-2\lambda n} \int_{B_{R+r_{n+1} } \times B_{R+r_{n+1}+s_{n+1}}} \vp_n(x)^2 (f(x) - f(y))^2 J(dx,dy) \nn\\
	&\quad + c_8(e^\lambda-1)^2 \sum_{n=1}^\infty  \sup_{z \in B_{R+2r}} \frac{e^{-2\lambda n}}{\phi(z,s_{n+1})}  \int_{B_{R+2r}}  g(x)^2 \p(x,s_{n+1})  \mu(dx)   \nn\\
	& =:c_8(e^\lambda-1)^2 (I_{3,1} + I_{3,2} + I_{3,3}).
\end{align}
Using the fact that $\vp_n^2\le 1$ in the first line below,  and \eqref{e:vp>ek} and the Fubini's theorem in the second, we get that
\begin{align}\label{e:I-3-1}
	I_{3,1}&\le    \sum_{n=1}^\infty  e^{-2\lambda n} \int_{B_R}  d\Gamma^{(L)}(f,f) +  \sum_{n=1}^\infty \sum_{k=1}^{n+1}  e^{-2\lambda n} \int_{B_{R+r_{k} } \setminus B_{R+r_{k-1}}}  d\Gamma^{(L)}(f,f)  \nn\\
	&\le  \sum_{n=1}^\infty  e^{-2\lambda n}\int_{B_R} \varphi^2 d\Gamma^{(L)}(f,f)   +    \sum_{k=1}^\infty \int_{B_{R+r_{k} } \setminus B_{R+r_{k-1}}} \varphi^2 d\Gamma^{(L)}(f,f)  \sum_{n=k-1}^{\infty}  e^{-2\lambda (n-k)} \nn\\
		&\le  \frac{1}{e^{2\lambda}-1} \int_{B_R} \varphi^2 d\Gamma^{(L)}(f,f) +  \frac{e^{4\lambda}}{e^{2\lambda}-1} \sum_{k=1}^\infty \int_{B_{R+r_{k} } \setminus B_{R+r_{k-1}}}  \varphi^2 d\Gamma^{(L)}(f,f)  \nn\\
		&\le \frac{e^{4\lambda}}{e^{2\lambda}-1}  \int_{B_{R+r }}  \varphi^2 d\Gamma^{(L)}(f,f).
\end{align}
Using \eqref{e:vp>ek}, we also see that
\begin{align}\label{e:I-3-2}
	I_{3,2} &\le    \int_{B_{R} \times B_{R+2r}}  \vp(x)^2(f(x) - f(y))^2 J(dx,dy)  \sum_{n=1}^\infty e^{-2\lambda n}\nn\\
	&\quad +   \sum_{k=1}^\infty \int_{(B_{R+r_{k}} \setminus B_{R+r_{k-1}}  ) \times B_{R+2r} } \vp(x)^2(f(x) - f(y))^2 J(dx,dy)  \sum_{n=k-1}^{\infty} e^{-2\lambda (n-k)}    \nn\\
	 &\le  \frac{e^{4\lambda}}{e^{2\lambda}-1} \int_{ B_{R+r} \times B_{R+2r} } \vp(x)^2(f(x) - f(y))^2 J(dx,dy).
\end{align}
For $I_{3,3}$, using \eqref{e:phi-scale} and the monotonicity of $\p$, we get
\begin{align}\label{e:I-3-3}
I_{3,3} &\le \sup_{z \in B_{R+2r}} \frac{c_9}{\phi(z,r)}  \bigg(\int_{B_{R+2r}}  g(x)^2 \p(x,r) \mu(dx) \bigg) \sum_{n=1}^\infty e^{-2\lambda n} \bigg(\frac{r}{s_{n+1}}\bigg)^{\beta_2}\nn\\
&= \sup_{z \in B_{R+2r}} \frac{c_{10}}{\phi(z,r)} \bigg(\int_{B_{R+2r}}  g(x)^2 \p(x,r) \mu(dx) \bigg) \sum_{n=1}^\infty e^{-\lambda (n-1)} \nn\\
&=\sup_{z \in B_{R+2r}} \frac{c_{10}e^\lambda }{(e^\lambda-1)\phi(z,r)}\int_{B_{R+2r}}  g(x)^2 \p(x,r) \mu(dx).
\end{align}
Besides, by the strong local property of $d\Gamma^{(L)}$ and   \eqref{e:CS-1}, we see that
\begin{align}\label{e:CS-improv-5}
	\int_{B_{R+2r}} g^2 \,d \Gamma^{(L)}(\vp,\vp)&= \sum_{n=1}^\infty  (e^{-\lambda(n-1)}-e^{-\lambda n})^2	\int_{B_{R+r_{n+1}}} g^2 \,d \Gamma^{(L)}(\vp_n,\vp_n) \nn\\
	&\le c_{11}(e^{\lambda}-1)^2(I_{3,1} + I_{3,2} +I_{3,3}).
\end{align}

 Now since $f^2 \le g^2$, combining  \eqref{e:CS-improv-1}--\eqref{e:CS-improv-5},  we obtain
	\begin{align}\label{e:CS-improv-6}
	& \int_{B_{R+2r}} f^2\, d\Gamma(\varphi,\varphi) \le 	\int_{B_{R+2r}} g^2 \,d \Gamma^{(L)}(\vp,\vp) + 	\int_{B_{R+2r}} g^2 \,d \Gamma^{(J)}(\vp,\vp)\nn\\
	&\le  \frac{c_{12}e^{4\lambda}(e^{\lambda}-1)^2}{e^{2\lambda}-1}\bigg(\int_{B_{R+r}} \varphi^2 d\Gamma^{(L)}(f,f) + \int_{B_{R+r} \times B_{R+2r}} \varphi(x)^2 (f(x) - f(y))^2 J(dx,dy) \bigg) \nn\\
	&\quad +  \sup_{z \in B_{R+2r}} \frac{c_{13}(\lambda)}{\phi(z,r)}\int_{B_{R+2r}}  g(x)^2 \p(x,r) \mu(dx),
\end{align}	
where $c_{13}(\lambda)$ is a constant depending on $\lambda$. 
 By \eqref{e:blow-up-integral} and the choice of $\theta$, we have
\begin{align*}
	\int_{B_{R+2r}}  g(x)^2 \p(x,r) \mu(dx)& \le 2\int_{B_{R+2r}}  f(x)^2 \p(x,r) \mu(dx) +  2\theta^2 \int_{B_{R+2r}} \p(x,r)\mu(dx)\\
	&\le 2(1+c_{14})\int_{B_{R+2r}}  f(x)^2 \p(x,r) \mu(dx).
\end{align*}
Note that $\lim_{\lambda \to 0}e^{4\lambda}(e^\lambda-1)^2/(e^{2\lambda}-1) =0$. By choosing $\lambda>0$ small enough so that  $4c_{12}e^{4\lambda}(e^\lambda-1)^2/(e^{2\lambda}-1) \le \eps$, we conclude from \eqref{e:CS-improv-6} that \eqref{e:CS-improv} holds.

To complete the proof,  we prove that $\vp \in \sF_b$. Define $h_i=\sum_{n=1}^i (e^{-\lambda(n-1)}-e^{-\lambda n}) \vp_n$. Then $h_i \in \sF_b$ for all $i \ge 1$ and $\lim_{i \to \infty} h_i (x)=\vp(x)$. To establish $\vp \in \sF_b$, it suffices to prove that $(h_i)_{i \ge 1}$ is a $\sE$-Cauchy sequence.  For any $i>j$, following the arguments above for \eqref{e:CS-improv-6}, we get
\begin{align}\label{e:CS-check-1}
&\int_{B_{R+2r}} d\Gamma(h_i-h_j,h_i-h_j)\nn\\
&\le  \theta^{-2}\int_{B_{R+2r}} g^2 \,  d\Gamma\bigg(\sum_{n=j+1}^i (e^{-\lambda(n-1)}-e^{-\lambda n})\vp_n,\sum_{n=j+1}^i (e^{-\lambda(n-1)}-e^{-\lambda n})\vp_n \bigg) \nn\\
	&\le c_{15}\theta^{-2}e^{- \lambda j} \bigg(\int_{B_{R+r}} \varphi^2 d\Gamma^{(L)}(f,f) + \int_{B_{R+r} \times B_{R+2r}} \varphi(x)^2 (f(x) - f(y))^2 J(dx,dy) \bigg) \nn\\
&\quad + \sup_{z \in B_{R+2r}}  \frac{c_{16}(\lambda)\theta^{-2}e^{-\lambda j}}{\phi(z,r)} \int_{B_{R+2r}} g(x)^2\p(x,r) \mu(dx) .
\end{align}
On the other hand,  since supp$[h_i-h_j] \subset B_{R+r}$ and $h_i(y)-h_j(y)\le \sum_{n=j+1}^i (e^{-\lambda(n-1)}-e^{-\lambda n}) \le e^{-\lambda j}$ for all $y \in B_{R+r}$, using Lemma \ref{l:Tail}, \VD \ and \eqref{e:phi-scale}, we get
\begin{align*}
	&\int_{B_{R+2r}^c} d\Gamma(h_i-h_j,h_i-h_j)
	\le e^{-2\lambda j} \int_{B_{R+2r}^c \times B_{R+r}} J(dx,dy)\le \frac{c_{17}e^{-2\lambda j}V(x_0,R)}{\phi(x_0,R)} \bigg(\frac{R}{r}\bigg)^{\beta_2}.
\end{align*}
Combining this with \eqref{e:CS-check-1}, we conclude that $\sE(h_i-h_j,h_i-h_j)=\int_{M} d\Gamma(h_i-h_j,h_i-h_j)$ converges to $0$ uniformly as $i,j \to \infty$. The proof is complete. \qed

The following inequalities were established in  \cite[Propositions 5.2 and 5.4]{GHH}.

\begin{prop}
	Let $U \subset M$ be an open subset. For any $u \in \sF'_b$ and $\vp \in \sF_b$ such that $\vp=0$ in $U^c$, we have
	\begin{align}\label{e:EnergyofProduct-1}
		&\int_{U} \vp^2 d\Gamma^{(L)}(f,f) + 	\int_{U\times U} \vp(x)^2 (f(x)-f(y))^2J(dx,dy) \nn\\
		&\le 2 \sE(f,f\vp^2) + 4 \int_{U} f^2 d\Gamma(\vp,\vp)+  4 \int_{U\times U^c} \vp(x)^2f(x)f(y)J(dx,dy)
	\end{align}
	and
\begin{equation}\label{e:EnergyofProduct-2}
	\sE(f\vp,f\vp) \le \sE(f,f\vp^2) + \int_{U} f^2 d\Gamma(\vp,\vp)   + 2 \int_{U\times U^c} \vp(x)^2f(x)f(y)J(dx,dy).
\end{equation}
\end{prop}

\begin{lem}\label{l:EnergyofProduct}
	Suppose that \VD, \Tail \ and  \CS \ hold. Then there exists $C>0$ such that  for any $x_0 \in M$, $0<r\le R$ with $R+2r<R_0$ and any $f \in \sF'_b$, there exists a cutoff function $\vp \in \sF$ for $B(x_0,R) \Subset B(x_0,R+r)$ so that  the following holds:
\begin{align}\label{e:EnergyofProduct}
	\sE(f\vp, f\vp) &\le  \frac32\sE(f,f\vp^2) + 3 \int_{B(x_0,R+r)\times B(x_0,R+2r)^c} \vp(x)^2 f(x)f(y)  J(dx,dy) \nn\\
	&\quad + \sup_{z \in B(x_0, R +2r)} \frac{C}{\phi(z,r)} \int_{B(x_0,R+2r)} f(x)^2 \p(x,r) \mu(dx).
\end{align}	
\end{lem}
\pf Write  $B_s:=B(x_0,s)$ for $s>0$. 	By Lemma \ref{l:CS-improv}, there exists a cutoff function $\vp \in \sF_b$  for $B_R \Subset B_{R+r}$ satisfying \eqref{e:CS-improv}  with $\eps=1/8$.  Since $\vp=0$ in $B_{R+r}^c$, using \eqref{e:CS-improv} and \eqref{e:EnergyofProduct-1} (with $U=B_{R+2r})$, we obtain
	\begin{align*}
	&\int_{B_{R+2r}} f^2 d\Gamma(\vp,\vp) = 2\int_{B_{R+2r}} f^2 d\Gamma(\vp,\vp)  - \int_{B_{R+2r}} f^2 d\Gamma(\vp,\vp)\nn\\
	&\le  \frac14 \bigg[  2\sE(f,f\vp^2)  + 4\int_{B_{R+2r}} f^2 d\Gamma(\vp,\vp) + 4 \int_{B_{R+r}\times B_{R+2r}^c} \vp(x)^2f(x)f(y)J(dx,dy)\bigg]\\
	& \quad +  \sup_{z \in B_{R+2r}}  \frac{c_1}{\phi(z,r)} \int_{B_{R+2r}} f(x)^2 \p(x,r) \mu(dx) - \int_{B_{R+2r}} f^2 d\Gamma(\vp,\vp)\\
	&=\frac12 \sE(f,f\vp^2)  +  \int_{B_{R+r}\times B_{R+2r}^c} \vp(x)^2f(x)f(y)J(dx,dy)+  \sup_{z \in B_{R+2r}}  \frac{c_1}{\phi(z,r)} \int_{B_{R+2r}} f(x)^2 \p(x,r) \mu(dx).
\end{align*}	
Combining this and  \eqref{e:EnergyofProduct-2} (with $U=B_{R+2r}$), we arrive at \eqref{e:EnergyofProduct}. \qed

For open subsets $A$ and $B$ of $M$ with $A \Subset B$, the relative capacity $\Cap(A,B)$ is defined by
\begin{align*}
	\Cap(A,B)=\inf\left\{\sE(\vp,\vp):\vp \in \sF_b \text{ is a cutoff function for $A \Subset B$}\right\}.
\end{align*}

\begin{lem}\label{l:Cap}
		Suppose that \VD, \Tail  \ and 
	\CS \ hold.  Then there exists  $C>0$ such that for any $x_0 \in M$ and  $0<r \le R$ with $R+2r<R_0$, 
	\begin{align*}
		\Cap	\left(B(x_0,R), B(x_0,R+r)\right)\le C \frac{V(x_0,R)}{\phi(x_0,R)}  \bigg( \frac{R}{r}\bigg)^{\beta_2}.
	\end{align*}
\end{lem}
\pf Write $B_s:=B(x_0,s)$ for $s > 0$. By Lemma \ref{l:EnergyofProduct}, there exists a cutoff function  $\vp\in\sF_b$ for $B_R \Subset B_{R+r}$ satisfying \eqref{e:EnergyofProduct} with $f=1$.   Using \eqref{e:EnergyofProduct} in the first inequality below, $\vp^2\le 1$, \eqref{e:phi2} and the monotonicity of $\p$ in the second, 
Lemma \ref{l:Tail} and \eqref{e:blow-up-integral} in the third, 
and \VD \ and \eqref{e:phi-scale} in the last, we get
\begin{align*}
	\sE(\vp, \vp) &\le  c_1 \int_{B_{R+r}\times B_{R+2r}^c} \vp(x)^2  J(dx,dy)  + \sup_{z \in B_{R+2r}} \frac{c_2}{\phi(z,r)} \int_{B_{R+2r}} \p(x,r) \mu(dx)\\
	&\le  c_1 \int_{B_{R+r}\times B_{R+2r}^c}   J(dx,dy)  +  \frac{c_3}{\phi(x_0,R)}\bigg( \frac{R}{r}\bigg)^{\beta_2} \int_{B_{R+2r}} \p(x,R+2r) \mu(dx)\\
	&\le  \frac{c_4V(x_0,R+r)}{\phi(x_0,R+r)}\bigg( \frac{R+r}{r}\bigg)^{\beta_2} +  \frac{c_5V(x_0,R+2r)}{\phi(x_0,R)}\bigg( \frac{R}{r}\bigg)^{\beta_2}\le \frac{c_6V(x_0,R)}{\phi(x_0,R)}\bigg( \frac{R}{r}\bigg)^{\beta_2}.
\end{align*}	
Since $\Cap	\left(B_R, B_{R+r}\right) \le 	\sE(\vp, \vp) $, this proves the lemma.
\qed

\section{Weighted  non-local tail and corresponding Caccioppoli and $L^2$-mean value inequalities}\label{s:L2}

Let $x_0 \in M$, $0<r\le R < R_0$ and $u$ be a measurable function on $M$. When \eqref{e:Tail-transition} holds,  we define a 
\textit{weighted non-local tail}  $\NT(u;x_0,r,R)$ of $u$ with respect to $\p$ by
\begin{align}\label{e:def-tail}
	\NT(u;x_0,r,R):=\sup_{x \in B(x_0,r)} \frac{1}{\p(x,r)}\int_{ B(x_0,R)^c} |u(y)|\,J(x,dy).
\end{align}

In this section, we first establish a weighted version of Caccioppoli-type inequality in terms of the weighted non-local tail. Cf.,  \cite[Lemma 4.6]{CKW-memo}.
\begin{lem}\label{l:Caccio}
	Suppose that \VD, \Tail \  and \CS \ hold. Let  $x_0 \in M$, $R \in (0,R_0)$, $r \in (0,R/2)$ with $R+r<R_0$, and   	$u$ be a bounded subharmonic function  on $B(x_0,R+r)$. Define  $u_\theta=(u-\theta)_+$ for some $\theta>0$. 	Then there exists a cutoff function $\vp$ for $B(x_0,R-r) \Subset B(x_0,R)$ such that
	\begin{align}\label{e:Caccio}
		&\sE(u_\theta\vp,u_\theta\vp) \le   C\bigg(  \frac{1}{\phi(x_0,R)} \bigg(\frac{R}{r}\bigg)^{\beta_2}  +  \frac{\NT(u_+;x_0,R,R+r)}{\theta}  \bigg)
		 \int_{B(x_0,R+r)} u(x)^2 \p(x,R+r) \mu(dx),
	\end{align}
	where   $C>0$ is a constant independent of $x_0,r,R,u$ and $\theta$. 
\end{lem}
\pf Write  $B_s:=B(x_0,s)$ for $s>0$. By Lemma \ref{l:EnergyofProduct}, there exists a cutoff function $\vp$ for $B_{R-r} \Subset B_R$ satisfying \eqref{e:EnergyofProduct} with $R$ replaced by $R-r$ and $f$ by $u_\theta$. 
 Note that $u_\theta \vp \in \sF'_b$ and $0\le u_\theta\vp^2 \in \sF^{B_R}$. Further, since $u-\theta$ and the constant function $0$ are subharmonic in $D$, by \cite[Corollary 2.10(iv)]{ChKu}, $u_\theta$ is harmonic in $D$. By Theorem \ref{t:harmonicity},  it follows that  $
\sE(u_\theta,u_\theta\vp^2)\le 0$. Combining this with  \eqref{e:EnergyofProduct}, we obtain 
	\begin{equation}\label{e:Caccio1}
	\sE(u_\theta\vp, u_\theta\vp) \le 3\int_{B_{R}\times B_{R+r}^c} \vp(x)^2 u_\theta(x)u_\theta(y)  J(dx,dy) + \sup_{z \in B_{R+r}} \frac{c_1}{\phi(z,r)} \int_{B_{R+r}} u_\theta(x)^2 \p(x,r) \mu(dx).
\end{equation}	
Since $u_\theta^2\le u^2$, using \eqref{e:phi2} and the monotonicity of $\p$, we get
\begin{equation}\label{e:Caccio2}
	 \sup_{z \in B_{R+r}} \frac{1}{\phi(z,r)} \int_{B_{R+r}} u_\theta(x)^2 \p(x,r) \mu(dx) \le  \frac{c_2}{\phi(x_0,R)} \bigg(\frac{R}{r}\bigg)^{\beta_2}  \int_{B_{R+r}} u(x)^2 \p(x,R+r) \mu(dx).
\end{equation}
Moreover, since $\vp^2\le 1$, $|u_\theta| \le u_+$ and   $|u_\theta| \le  u^2/\theta$, we have
\begin{align}\label{e:Caccio3}
	\int_{B_{R}\times B_{R+r}^c} \vp(x)^2 u_\theta(x)u_\theta(y)  J(dx,dy) &\le \frac{1}{\theta}\int_{B_R}u(x)^2     \frac{\p(x,R+r) }{\p(x,R)}\int_{ B_{R+r}^c} u_+(y) J(x,dy)\,\mu(dx)\nn\\
	&\le \frac{\NT(u_+;x_0,R,R+r) }{\theta}  \int_{B_R} u(x)^2\p(x,R+r)\mu(dx). \qquad 
\end{align}
Combining \eqref{e:Caccio1} with  \eqref{e:Caccio2} and \eqref{e:Caccio3}, we arrive at the desired result.\qed

Recall that $\lambda_1(D)$ denotes the bottom of the spectrum of $-\sL^D$. 

\begin{lem}\label{l:L2-onestep}
	Let $x_0 \in M$, $0<r\le R$,  $\vp \in \sF_b$ be a cutoff function for $B(x_0,R)\Subset B(x_0,R+r/2)$,  $u \in \sF^{B(x_0,R+r)}_{\rm loc}$ be bounded in $B(x_0,R+r)$, and $u_\theta:=(u-\theta)_+$ for $\theta>0$.
	There exists  $C>0$ independent of $x_0,r,R,\vp,u$ and $\theta$ such that  for  any open set $D\subset M$ satisfying
	\begin{align}\label{e:L2-onestep-ass}
		\left\{x \in B(x_0,R+r/2):u(x)>\theta\right\} \subset D \subset  B(x_0,R+r),
	\end{align}
	it holds that 
	\begin{align*}
	\frac{\lambda_1(D)}{ (1+\phi(x_0,R)\lambda_1(D))^{1/(1+\gamma_1)}} \int_{B(x_0,R)} u_\theta(x)^2 \p(x,R) \mu(dx)\le  C   \sE(u_\theta\vp,u_\theta\vp).
	\end{align*}
\end{lem}
\pf  Write  $B_s:=B(x_0,s)$ for $s>0$.  By \eqref{e:blow-up-Hardy}, since $\vp=1$ in $B_R$, we have 
\begin{align}
	\int_{B_{R}} u_\theta(x)^2\p(x,R)^{1+\gamma_1}\mu(dx) \le c_1\left(  \phi(x_0,R) \sE(u_\theta\vp, u_\theta\vp) + \lVert u_\theta \vp \rVert_2^2 \right). \label{e:L2-onestep-1}
\end{align}
 By \eqref{e:L2-onestep-ass}, $D$ contains 	$\left\{x \in B_{R+r/2}:u_\theta(x)\neq 0\right\}$ and $u_\theta\vp=0$ in $B_{R+r} \subset D^c$. Thus, 
\begin{align}\label{e:L2-onestep-0}
	\int_{B_R}u_\theta^2\, d\mu \le  	\int_{D}u_\theta^2 \vp^2 d\mu  =  \lVert u_\theta \vp \rVert_2^2.
\end{align}
Further, by  \eqref{e:eigen},  we get
\begin{align}\label{e:L2-onestep-2} 
\lambda_1(D)\lVert u_\theta \vp \rVert_2^2\le    \sE(u_\theta\vp,u_\theta\vp).
\end{align}
Using H\"older inequality,   \eqref{e:L2-onestep-1}, \eqref{e:L2-onestep-0} and \eqref{e:L2-onestep-2}, we conclude that 
\begin{align*}
	&	\lambda_1(D)	 \int_{B_R} u_\theta(x)^2 \p(x,R) \mu(dx) \\
	&\le \lambda_1(D)\bigg(\int_{B_{R}} u_\theta^2 \,d\mu\bigg)^{\gamma_1/(1+\gamma_1)} \bigg( \int_{B_{R}} u_\theta(x)^2\p(x,R)^{1+\gamma_1}\mu(dx) \bigg)^{1/(1+\gamma_1)}   \nn\\
	&\le  c_2\lambda_1(D) \lVert u_\theta \vp \rVert_2^{2\gamma_1/(1+\gamma_1)}\left(  \phi(x_0,R) \, \sE(u_\theta\vp,u_\theta\vp)  \right)^{1/(1+\gamma_1)}    +  c_2\lambda_1(D) \lVert u_\theta \vp \rVert_2^{2}  \nn\\
	&\le c_3\left( 1+(\phi(x_0,R)\lambda_1(D))^{1/(1+\gamma_1)}  \right) \sE(u_\theta\vp,u_\theta\vp).
\end{align*} \qed 

The following lemma is a weighted version of \cite[Lemma 4.8]{CKW-memo}.

\begin{lem}\label{l:Gr0-Lemma3.2}
	Suppose that \VD, \Tail, 
	\CS \  and \FK \ hold.	Let $x_0 \in M$   and $0<r\le R$ with $R+r <R_0/q_1$, where $q_1\ge 1$ is the constant in \FK.
		Let $u$ be  a bounded subharmonic function on $B(x_0,R+r)$ and $\theta>0$ be a constant. Set $u_\theta:=(u-\theta)_+$, 
	\begin{align*}
		\sI_0:= \int_{B(x_0,R+r)} u(x)^2 \p(x,R+r) \mu(dx) \quad \text{and} \quad \sI_1:= \int_{B(x_0,R)} u_\theta(x)^2 \p(x,R) \mu(dx).
	\end{align*}
	Then there exists  $C>0$ independent of $x_0,R,r,u$ and $\theta$ such that 
	\begin{align}\label{e:Gr0-Lemma3.2}
		\sI_1\le  \frac{C\phi(x_0,R)}{\theta^{2\nu_0}V(x_0,R)^{\nu_0}} \bigg( \frac{1}{\phi(x_0,R)} \bigg(\frac{R}{r}\bigg)^{\beta_2}   + \frac{\NT(u_+;x_0,R+ r/2,R+r)}{\theta}\bigg) \,\sI_0^{1+\nu_0},
	\end{align}
	where 
	\begin{align}\label{e:def-nu0}
		\nu_0:=\frac{\gamma_1\nu}{1+\gamma_1},
	\end{align}
and $\gamma_1$ and $\nu$ are the constants in \eqref{e:blow-up-integral-general} and  \FK \ respectively.
\end{lem}
\pf  Write $B_s:=B(x_0,s)$ for $s>0$. By Lemma \ref{l:Caccio},  there exists a cutoff function  $\vp \in \sF_b$   for $B_{R} \Subset B_{R+r/2}$   satisfying \eqref{e:Caccio}  with $R$ replaced by $R+r/2$ and $r$ by $r/2$. Set
$$
D_0:=	\left\{x \in B_{R+r/2}:u(x)>\theta\right\} = 	\left\{x \in B_{R+r/2}:u_\theta(x)\neq 0\right\}.
$$
If $\mu(D_0)=0$, then $\sI_1=0$ so that \eqref{e:Gr0-Lemma3.2} clearly holds. Suppose that $\mu(D_0)>0$. Since $\mu$ is regular, there exists 
an open set $D\subset M$  satisfying $D_0 \subset D\subset B_{R+r}$ and $\mu(D)\le 2\mu(D_0)$.  Since $\p \ge 1$, using Markov's inequality, we have $\mu(D) \le 2\mu(D_0)  \le2 \sI_0/\theta^2.$ 
Hence, by \FK, \VD \ and \eqref{e:phi-scale}, 
\begin{align}
	\lambda_1(D) \ge \frac{c_1}{\phi(x_0,R)} \bigg( \frac{V(x_0,R)}{\mu(D)}\bigg)^\nu  \ge \frac{c_1}{\phi(x_0,R)} \bigg( \frac{\theta^2V(x_0,R)}{2\sI_0}\bigg)^\nu.
	\label{e:cond-32}\end{align}
 Applying Lemma \ref{l:L2-onestep} and using \eqref{e:cond-32} and \eqref{e:Caccio},  we arrive at 
	\begin{align*}
\sI_1&\le  \frac{c_1(1+\phi(x_0,R)\lambda_1(D))^{1/(1+\gamma_1)}}{\lambda_1(D)}     \sE(u_\theta\vp,u_\theta\vp)\\
&\le \frac{c_2 \phi(x_0,R) \sI_0^{1+\gamma_1\nu/(1+\gamma_1)}}{ (\theta^2 V(x_0,R))^{\gamma_1\nu/(1+\gamma_1)}}   \bigg(  \frac{1}{\phi(x_0,R)} \bigg(\frac{R}{r}\bigg)^{\beta_2}  +  \frac{
	\NT(u_+;x_0,R+r/2,R+r)}{\theta}  \bigg).
\end{align*}
\qed

We will use the following elementary iteration lemma.
\begin{lem}[\!\!{\cite[Lemma 4.9]{CKW-memo}}]\label{l:iteration}
 Let $(a_n)_{n \ge 0}$ be a sequence of positive real numbers such that
 $$
 a_{n+1} \le c_0 b^n a_n^{1+\delta} \quad \text{for all $n \ge 0$},
 $$
 for some constants $\delta>0$, $b>1$ and $c_0>0$. If $a_0 \le c_0^{-1/\delta}\,b^{-1/\delta^2}$, 
 then $\lim_{n \to \infty}a_n=0$.
\end{lem}

We now establish the  $L^2$-weighted mean value inequality for subharmonic functions.

\begin{prop}\label{p:L2mean} 
 	Suppose that \VD, \Tail, 
 	\CS \  and \FK \  hold.  Let   $x_0 \in M$ and   $R \in (0,R_0/q_1)$ where $q_1 \ge1$ is the constant in \FK. For any $p>0$, there exists a constant $C=C(p)>0$ independent of $x_0$ and $R$ such that  for any bounded subharmonic function $u$  on $B(x_0,R)$,
	\begin{align*}
		\esssup_{ B(x_0, R/2)} u& \le C \big( (1+\delta^{-1})^{1/(2\nu_0)} \, \sI + \delta \sT(p) \big)\quad \text{for all} \;\, \delta>0,
	\end{align*}
where
\begin{align}\label{e:L2-Caccio}
	\sI&:=\Big(  \frac{1}{V(x_0,R)}\int_{B(x_0,R)} u(x)^2\p(x,R)\mu(dx)  \Big)^{1/2},\nn\\
	\sT(p)&:= \phi(x_0,R) \sup \left\{({s}/{R})^p\,\NT(u_+;x_0,r,r+s): R/2< r\le R,\,0<s\le R/4\right\},
\end{align}
and the constant $\nu_0>0$ is  defined by \eqref{e:def-nu0}.
\end{prop}
\pf Write $B_s:=B(x_0,s)$ for $s>0$. Let $u$ be a bounded subharmonic function  on $B_R$ and $\delta>0$. For a constant $\theta>0$ whose exact value will be determined later,  define 
\begin{align*}
 s_n=\frac{1}{2}(1+2^{-n})R, \quad \theta_n=(1-2^{-n})\theta \quad \text{and} \quad a_n= \int_{B_{s_n}} (u(x)-\theta_n)_+^2 \,\p(x,s_n) \mu(dx), \quad n\ge0.
\end{align*}
Note that $s_n \in (R/2, R]$ for all $n \ge 0$.  For all $n \ge 0$, by applying Lemma \ref{l:Gr0-Lemma3.2} (with $R=s_{n+1}$, $r=s_{n}-s_{n+1}$ and $\theta=\theta_{n+1}-\theta_n$),  we get that \begin{align}\label{e:iteration-1}
&	a_{n+1} = \int_{B_{s_{n+1}}} [(u(x)-\theta_n) - (\theta_{n+1}-\theta_n)]_+^2 \,\p(x,s_{n+1}) \mu(dx)\nn\\
	&\le \frac{c_1 a_n^{1+\nu_0}}{(\theta_{n+1}-\theta_n)^{2\nu_0}V(x_0,s_{n+1})^{\nu_0}}  \bigg[ \bigg(\frac{s_{n+1}}{s_{n}-s_{n+1}}\bigg)^{\beta_2}   + 
\frac{\phi(x_0,s_{n+1})}{\theta_{n+1}-\theta_n}\NT(u_+;x_0,(s_n+s_{n+1})/2,s_n)\bigg] \nn\\
	&\le \frac{c_1 a_n^{1+\nu_0}}{(\theta_{n+1}-\theta_n)^{2\nu_0}V(x_0,s_{n+1})^{\nu_0}}  \bigg[ \bigg(\frac{s_{n+1}}{s_{n}-s_{n+1}}\bigg)^{\beta_2}   + \frac{ (2R/(s_n-s_{n+1}))^p }{\theta_{n+1}-\theta_n}\sT(p)\bigg] \nn\\
	&\le \frac{c_2 a_n^{1+\nu_0}}{(2^{-n-1}\theta)^{2\nu_0}\,V(x_0,R)^{\nu_0}}  \bigg[(2^{n+1})^{\beta_2} +\frac{(2^{n+2})^p }{2^{-n-1}\theta}\sT(p)\bigg] \nn\\
	&\le \frac{2^{(2\nu_0 +\beta_2+p+1)  n } c_3 a_n^{1+\nu_0}}{\theta^{2\nu_0}V(x_0,R)^{\nu_0}}  \left[1 + \frac{\sT(p)}{\theta}\right].
\end{align}
Set $b:=2^{2\nu_0 +\beta_2+p+1}$ and  $c_4:=c_3^{1/\nu_0} b^{1/\nu_0^2}$, and  take
\begin{align*}
	\theta:=  \big({c_4 (1+\delta^{-1})^{1/\nu_0} a_0}/{V(x_0,R)}\big)^{1/2} +\delta \sT(p).
\end{align*}
Let $	c_0:=c_3\theta^{-2\nu_0}(1+ \theta^{-1}\sT(p)) V(x_0,R)^{-\nu_0} .$ By \eqref{e:iteration-1}, $a_{n+1} \le c_0 b^n a_n^{1+\nu_0}$ for all $n \ge 0$. Further,
\begin{align*}
	&c_0^{-1/\nu_0} b^{-1/\nu_0^2} = c_3^{-1/\nu_0} b^{-1/\nu_0^2} \frac{ \theta^{2} V(x_0,R)}{(1+\theta^{-1}\sT(p))^{1/\nu_0}}
\ge  c_3^{-1/\nu_0} b^{-1/\nu_0^2}\frac{c_4(1+\delta^{-1})^{1/\nu_0}a_0}{ (1 +\delta^{-1})^{1/\nu_0}}=a_0.
\end{align*}
Therefore, by Lemma \ref{l:iteration}, we  obtain
$\int_{B_{R/2}} (u(x)-\theta)_+^2\,\p(x,R/2) \mu(dx) \le \lim_{n \to \infty} a_n=0.$
Since $\p\ge 1$, we conclude that $\esssup_{B_{R/2}} u \le \theta$.  \qed

\begin{cor}\label{c:L2mean-diff-form}
		Suppose that \VD, \Tail, \CS \
	  and \FK \  hold.  Let   $x_0 \in M$,   $R \in (0,R_0/q_1)$ where $q_1 \ge 1$ is the constant in \FK.  For any $p>0$, there exists a constant $C=C(p)>0$ independent of $x_0$ and $R$ such that for any bounded subharmonic function $u$ on $B(x_0,R)$, 
\begin{align}\label{e:L2mean-diff-form-1}
	&\esssup_{B(x_0, R/2)} u \, \le C\sI^{\nu_1} \,( \sI \vee \sT(p))^{1-\nu_1}, 
\end{align}
where $\sI$ and $\sT(p)$  are defined by \eqref{e:L2-Caccio}, 
\begin{align}\label{e:def-nu1}
	\nu_1:=2\gamma_1\nu/(1+\gamma_1+2\gamma_1\nu),
\end{align}
and $\gamma_1$ and $\nu$ are the constants in \eqref{e:blow-up-integral-general} and  \FK, respectively. In particular, if $u$ is bounded in $M$, there exists $C>0$ independent of $x_0,R$ and $u$ such that
\begin{align}\label{e:L2mean-diff-form-2}
	&\esssup_{B(x_0, R/2)} u \, \le C\sI^{\nu_1} \,\Big( \sI \vee \esssup_{B(x_0,R/2)^c} u_+ \Big)^{1-\nu_1}.
\end{align}
\end{cor}
\pf Let $\nu_0:=\gamma_1\nu/(1+\gamma_1)$. Then we have $\nu_1=2\nu_0/(1+2\nu_0)$. Applying  Proposition \ref{p:L2mean} with $\delta=(\sI/\sT(p))^{\nu_1}$, we get that
\begin{align*}
	\esssup_{ B(x_0, R/2)} u \, &\le c_1   \big(  (1+\delta^{-1})^{1/(2\nu_0)} \sI  +  \delta \sT(p)\big) \\[-3mm]
	&\le c_1 \begin{cases}
		2^{1/(2\nu_0)} \sI  + \sI	&\mbox{ if }  \sI \ge \sT(p),\\
		2^{1/(2\nu_0)} (\sI/\sT(p))^{-\nu_1/(2\nu_0)}\sI + \sI^{\nu_1}\sT(p)^{1-\nu_1}		&\mbox{ if }  \sI < \sT(p)
	\end{cases}\\
	&= (2^{1/(2\nu_0)}+1)c_1 \sI^{\nu_1} (\sI \vee \sT(p))^{1-\nu_1}.
\end{align*}
Hence, \eqref{e:L2mean-diff-form-1} holds. For \eqref{e:L2mean-diff-form-2}, assume additionally   that $u$ is bounded in $M$.  For any  $R/2<r\le R$, $0<s\le R/4$ and $x \in B(x_0,r)$, by \Tail, the monotonicity of $\p$ and \eqref{e:phi2}, we have
\begin{align*}
&	\frac{\phi(x_0,R)}{\p(x,r)}\bigg(\frac{s}{R}\bigg)^{\beta_2}\int_{B(x_0,r+s)^c} u_+(y)J(x,dy) \le \frac{ \phi(x_0,R)}{\p(x,r)}\bigg(\frac{s}{R}\bigg)^{\beta_2}  J(x,B(x,s)^c) \esssup_{B(x_0,R/2)^c} u_+ \\
	&\le \frac{c_2\phi(x_0,R)}{\phi(x,s) }\bigg(\frac{s}{R}\bigg)^{\beta_2}  \esssup_{B(x_0,R/2)^c} u_+ \le c_3\esssup_{B(x_0,R/2)^c} u_+.
\end{align*}
Hence, $\sT(\beta_2)\le c_2k_0$ so that \eqref{e:L2mean-diff-form-2} follows from \eqref{e:L2mean-diff-form-1} (with $p=\beta_2$).	\qed

In the remainder of this section, we establish \mE \ under \VD,  \Tail, \CS \ and  \FK \  (Proposition \ref{p:CS->E}) through a ``lemma of growth". This terminology was introduced in \cite{GHH18}, with its original concept  tracing back to \cite{Lan98}, where second-order elliptic operators were studied.

\begin{lem}\label{l:LG}
{\rm \!(Lemma of growth)}	Suppose that \VD, \Tail, 
\CS \ and \FK \ hold. Then  there exist    $\delta_0,\eps_0 \in (0,1)$ such that the following  holds:  Let $x_0 \in M$ and $R \in (0,R_0/q_1)$ where $q_1 \ge 1$ is the constant in \FK. For any non-negative $u \in \sB_b(M)$ that is superharmonic in $B(x_0,R)$, if
	\begin{align}\label{e:LG-ass}
		\mu (  B(x_0,R)  \cap \{ u < a\} ) \le \delta_0 V(x_0,R)
	\end{align}
 for some constant $a>0$,	then
	\begin{align*}
		\essinf_{B(x_0,R/2)} u \ge \eps_0 a.
	\end{align*} 
\end{lem}
\pf Write $B_s:=B(x_0,s)$ for $s>0$. For a constant $h>0$ whose exact value is to be determined  later, we choose  a function $F_h\in C^2(\R)$ such that
\begin{align}\label{e:condition-F-1}
	F_h (s)=(s+h)^{-1}\;\; \text{for all $s>-h/2$}, \qquad 	F_h (s)\le -10 h^{-2}s \;\; \text{for all $s\le -h/2$}
\end{align}
and
\begin{align}\label{e:condition-F-2}
	\sup_{s \in \R} F_h '(s)\le 0, \quad \inf_{s \in \R} F_h ''(s)\ge 0, \quad \sup_{s \in \R} \,(|F_h '(s)| + |F_h ''(s)|)<\infty.
\end{align}
By \cite[(11.3)]{GHH}, $F_h(u)$ is subharmonic in $B$.
Applying Corollary \ref{c:L2mean-diff-form} to $F_h(u)$, we get
\begin{align}\label{e:LG-eq-0}
	\Big(\essinf_{B_{R/2}} u +h \Big)^{-1} = \esssup_{B_{R/2}}F_h(u) \le c_1\sI^{\nu_1} \Big( \sI \vee \esssup_{B_{R/2}^c} F_h(u) \Big)^{1-\nu_1},
\end{align}
where $\nu_1>0$ is the constant in \eqref{e:def-nu1} and
\begin{align*}
 \sI:=\bigg(  \frac{1}{V(x_0,R)}\int_{B_R} F_h(u(x))^2\p(x,R)\mu(dx)  \bigg)^{1/2}.
\end{align*} 
Since $u\ge 0$  in $M$, we see that $F_h(u(x)) \le h^{-1}$ for all $x \in M$. Thus, by \eqref{e:blow-up-integral}, we have 
\begin{align}\label{e:LG-eq-1}
  \sI \vee \esssup_{B_{R/2}^c} F_h(u)  \le h^{-1}\bigg( 1 \vee \frac{1}{V(x_0,R)}\int_{B_R} \p(x,R)\mu(dx)  \bigg)^{1/2}\le c_2h^{-1}.
\end{align}
On the other hand,  using \eqref{e:blow-up-integral0} and \eqref{e:LG-ass}, we get
\begin{align*}
	&\frac{1}{V(x_0,R)}\int_{B_{R} \cap \{ u < a\}} F_h(u(x))^2\p(x,R)\mu(dx)\le 	\frac{h^{-2}}{V(x_0,R)}\int_{B_{R} \cap \{ u < a\}} \p(x,R)\mu(dx)\\
	&\le C^{1/(1+\gamma_1)} \frac{h^{-2} 	\mu (  B_R  \cap \{ u < a\} )^{\gamma_1/(1+\gamma_1)}}{V(x_0,R)^{\gamma_1/(1+\gamma_1)}}\le c_3\delta_0^{\gamma_1/(1+\gamma_1)}h^{-2}.
\end{align*}
Using this and \eqref{e:blow-up-integral}, we obtain
\begin{align}\label{e:LG-eq-2}
	\sI^2 &\le \frac{1}{V(x_0,R)} \bigg( \int_{B_{R} \cap \{ u \le a\}} F_h(u(x))^2\p(x,R)\mu(dx)+ (a+h)^{-2}  \int_{B_R  \cap \{ u > a\}} \p(x,R)\mu(dx)\bigg)\nn \\
	&\le c_4\left( \delta_0^{\gamma_1/(1+\gamma_1)}h^{-2}+ (a+h)^{-2}  \right).
\end{align}

Now we take
$
h=\delta_0^{\gamma_1/(2+2\gamma_1)}a / (1-\delta_0^{\gamma_1/(2+2\gamma_1)})
$
so that $\delta_0^{\gamma_1/(1+\gamma_1)}h^{-2}= (a+h)^{-2} $. Combining \eqref{e:LG-eq-0} with  \eqref{e:LG-eq-1} and \eqref{e:LG-eq-2}, we arrive at
\begin{align*}
	\essinf_{B_{R/2}} u  & \ge  c_1^{-1}\sI^{-\nu_1} (c_2h^{-1})^{\nu_1-1} -h \ge  c_5 (\delta_0^{\gamma_1/(1+\gamma_1)}h^{-2})^{-\nu_1/2}\,h^{1-\nu_1} -h\\
	&= \frac{\delta_0^{\gamma_1/(2+2\gamma_1)}(c_5 \delta_0^{-\gamma_1\nu_1/(2+2\gamma_1)}-1)}{1-\delta_0^{\gamma_1/(2+2\gamma_1)}} a.
\end{align*}
By letting $\delta_0$ be a small constant satisfying $c_5 \delta_0^{-\gamma_1\nu_1/(2+2\gamma_1)}\ge 2$, we arrive at the  result. \qed

\begin{prop}\label{p:CS->E}
	The following implication holds: 	\begin{align*}		\text{\rm \VD \, $+$\, \Tail \, $+$\, \CS \, $+$\,   \FK \; $\Rightarrow$ \ \mE}.	\end{align*} 
\end{prop}
\pf \mEu \  follows from Proposition \ref{p:FK}. Moreover, by  Lemmas \ref{l:Cap} and  \ref{l:LG}, conditions (LG) and (Cap$_\le$) in \cite{GHH} hold. Hence, by \cite[Lemma 12.4]{GHH}, \mEl  \ holds. \qed

\section{Equivalences among \CS, \CSU, \Gcap \ and \GU}\label{s:proof-thm-1}
In this section, we show that all four weighed versions of functional inequalities \CS, \CSU, \Gcap \ and \GU\ are equivalent.

The next lemma is standard. See \cite[Proposition 13.4]{GHH} and \cite[Lemma 4.16]{CKW-memo} for the proof.

\begin{lem}\label{l:EP}
	Suppose that \mE \ holds.	Then there exist constants $\eps, \delta \in (0,1)$ such that  
	\begin{align*}
		\essinf_{z \in B(x_0,r/4)}	\P^z\left(\tau_{B(x_0,r)} > \eps \phi(x_0,r)\right) \ge \delta \quad \text{ for all $x_0 \in M$ and $r\in (0,R_0)$.} 
	\end{align*}
\end{lem}

Since $\p \ge 1$, by  \cite[Lemma 13.5]{GHH} and  Lemma \ref{l:EP}, we obtain the next lemma.

\begin{lem}\label{l:E->GU}
\mE \ implies \GU.  More precisely, if \mE \ holds, then there exists a constant $\kappa \ge 1$ such that for any $x_0 \in M$ and $0<r\le R$ with $R+r<R_0$, the function
\begin{align*}
\vp:=\kappa \lambda G_\lambda^{B(x_0,R+r)}\1_{B(x_0,R+r)} \quad \text{where} \;\; \lambda:= \Big(\inf_{z \in B(x_0,R)}\phi(z,r) \Big)^{-1}
\end{align*}
belongs to $\sF_b$, is  a $\kappa$-cutoff function for $B(x_0,R) \Subset B(x_0,R+r)$ and satisfies \eqref{e:Gcap} for all $f \in \sF_b'$.
\end{lem}

The following result is an analogue of \cite[Lemma 2.4]{GHH18} and  \cite[Proposition 2.5]{CKW-adv}.

\begin{lem}\label{l:GU->SCS}
Suppose that \VD \ and \Tail \ hold. Then \Gcap \ implies \CS,  and \GU \ implies \CSU. 
	More precisely, for any $x_0 \in M$, $0<r\le R$ with $R+2r<R_0$ and any $f \in \sF'_b$, if $\vp \in \sF_b$ is a $\kappa$-cutoff function for $B(x_0,R) \Subset B(x_0,R+r)$ satisfying \eqref{e:Gcap}, then $\vp \wedge 1 \in \sF_b$ is a cutoff function for $B(x_0,R) \Subset B(x_0,R+r)$ satisfying \eqref{e:CS}.
\end{lem}
\pf Since the proofs are similar, we only give the proof of the second implication.

Assume that  \GU \ holds.  Fix $x_0 \in M$ and $0<r \le R$ with $R+2r <R_0$.  Set $V_0:=B(x_0,R)$, $V_1:=B(x_0,R+r)$, $V_2:=B(x_0,R+3r/2)$ and $V_3:=B(x_0, R+2r)$. By \GU,  there exist a constant $\kappa \ge 1$ and  a $\kappa$-cutoff function $\vp \in \sF_b$ for $V_0\Subset V_1$ such that for any $f \in \sF'_b$,
\begin{align}\label{e:SCS-1}
	\sE(f^2 \vp, \vp) \le \sup_{z \in V_1} \frac{c_1}{\phi(z,r)} \int_{V_1} f(x)^2 \p(x,r) \mu(dx).
\end{align}
Set $\vp_1:=\vp \wedge 1$ and let $f \in \sF'_b$.  Since $\vp_1=0$ in $V_2^c$, by the strong local property,  we have
\begin{align}\label{e:SCS-2}
	\int_{V_3 \setminus V_2} f^2 d\Gamma^{(L)}(\vp_1,\vp_1)=0.
\end{align}
Using \Tail, \eqref{e:phi-scale} and the monotonicity of $\p$, since $\vp_1^2\le 1$, we obtain
\begin{align}\label{e:SCS-3}
	&\int_{(V_3 \times M) \setminus (V_2 \times V_2)} f(x)^2 (\vp_1(x)-\vp_1(y))^2 J(dx,dy)\nn\\
	&= \int_{(V_3 \setminus V_2) \times V_1}  f(x)^2 \vp_1(y)^2 J(dx,dy) + \int_{V_1 \times (M \setminus V_2)} f(x)^2 \vp_1(x)^2 J(dx,dy)\nn\\
	&\le \int_{V_3 \setminus V_2} f(x)^2 J (x, B(x,r/2)^c) \mu(dx)  + \int_{V_1} f(x)^2 J(x, B(x,r/2)^c) \mu(dx)\nn\\
	&\le  \sup_{z \in V_3} \frac{c_2}{\phi(z,r)} \int_{V_3} f(x)^2 \p(x,r) \mu(dx).
\end{align}
Further, using \cite[(2.5)]{CKW-adv} in the third line below and   $\vp \le \kappa \vp_1$ in the fourth, we get
\begin{align}\label{e:SCS-4}
	&\int_{V_2} f^2 \, d\Gamma^{(L)}(\vp_1,\vp_1) + \int_{V_2 \times V_2} f(x)^2 (\vp_1(x)-\vp_1(y))^2 J(dx,dy)\nn\\
		&\le \int_{V_2} f^2 \, d\Gamma^{(L)}(\vp,\vp) + \int_{V_2 \times V_2} f(x)^2 (\vp(x)-\vp(y))^2 J(dx,dy)\nn\\
	&\le  4 \left( \int_{V_2} \vp^2 d \Gamma^{(L)}(f,f)  + \int_{V_2 \times V_2} \vp(x)^2 (f(x)-f(y))^2 J(dx,dy)\right)  + 2\sE(f^2\vp,\vp)\nn\\
	&\le 4\kappa^2 \bigg(\int_{V_2} \vp_1^2 d\Gamma^{(L)}(f,f) + \int_{V_2 \times V_2} \vp_1(x)^2 (f(x)-f(y))^2 J(dx,dy) \bigg) + 2\sE(f^2\vp,\vp).
\end{align}
Combining \eqref{e:SCS-2}, \eqref{e:SCS-3} and  \eqref{e:SCS-4} with  \eqref{e:SCS-1},   we arrive at the desired result. \qed

The proof of Theorem \ref{t:main-1} is now immediate.

\medskip

\noindent \textbf{Proof of Theorem \ref{t:main-1}.}
Clearly, \CSU \ implies \CS, and \GU \ implies \Gcap.   Moreover, under \VD, \Tail \ and \FK, 
by  Proposition \ref{p:CS->E},  and Lemmas \ref{l:E->GU} and  \ref{l:GU->SCS}, we have
\begin{align*}
	\text{\rm 
		\Gcap
		\,}  \Rightarrow  \text{\rm  \
		\CS
		\,}  \Rightarrow \text{\rm \  \mE \,} \Rightarrow \text{\rm  \ \GU\,}\Rightarrow \text{ \CSU \  $+$ \Gcap}.
\end{align*}
This completes the proof. \qed

The original formulation of the cutoff Sobolev inequality  introduced in \cite{BB04, BBK06} includes the requirement of   H\"older continuity for  cutoff functions (see condition CS$(\Psi)$ in \cite{BB04, BBK06}). This regularity was proven using a parabolic Harnack inequality.  Under an additional assumption \PHR, we obtain the following interior H\"older regularity for cutoff functions in \CSU.

\begin{lem}\label{l:CS-continuity}
	 Suppose that \VD, \Tail, \mE \ and \PHR \ hold. Then there exist constants  $C>0$ and $\chi \in (0,1]$  such that for any $x_0 \in M$ and $0<r\le R$ with $R+2r<R_0$, there exists a  cutoff function $\vp \in \sF_b$ for $B(x_0,R) \Subset B(x_0,R+r)$ satisfying \eqref{e:CS} for all $f \in \sF'_b$ and
	 \begin{align}\label{e:CS-continuity}
	 	|\vp (x)-\vp(y)| \le \sup_{ z \in B(x_0,R)}\frac{C\phi(x,r)}{\phi(z,r)} \bigg( \frac{d(x,y)}{r}\bigg)^\chi \quad \text{for all} \;\; x,y \in B(x_0, R+r/2).
	 \end{align}
\end{lem}	 	 
	\pf		Let $x_0 \in M$ and $0<r\le R$ with $R+r<R_0$. Set $V_0:=B(x_0,R)$ and $V:=B(x_0,R+r)$. By Lemmas \ref{l:E->GU} and \ref{l:GU->SCS},  there exists a constant $\kappa \ge 1$ independent of $x_0,r$ and $R$ such that for $\lambda:= (\inf_{z \in V_0}\phi(z,r) )^{-1}$, the function $	\vp:=(\kappa \lambda G_\lambda^{V}\1_{V}) \wedge 1$
	belongs to $\sF_b$,	is a cutoff function for $V_0 \Subset V$ and satisfies  \eqref{e:CS} for all $f \in \sF'_b$.
	
	 Now we prove that \eqref{e:CS-continuity} holds. Let $x \in B(x_0,R+r/2)$. By \PHR, since the map $(t,z) \mapsto P_t^V \1_V(z)$ is  bounded above by 1 and  caloric in $(0,\infty) \times B(x, r/2)$, we get that
		\begin{equation}\label{e:CS-PHR}
			\left| P_t^V \1_V (x) - P_t^V \1_V(y)\right| \le c_1 \bigg(\frac{d(x,y)}{\phi^{-1}(x, t) \wedge (r/2)} \bigg)^\theta \quad \text{for all} \;\, t>0, \; y\in B(x, \delta r/2), 
		\end{equation}
		where $\delta \in (0,1)$ and $\theta>0$ are the constants in \PHR. Without loss of generality, we assume that $\theta<\beta_1$ where $\beta_1>0$ is the constant in \eqref{e:phi-scale}. Then by \eqref{e:phi-scale}, we have
		\begin{align}\label{e:CS-PHR-2}
			\int_{0}^{\phi(x,r/2)} \frac{dt}{\phi^{-1}(x,t)^{\theta}}  \le  \frac{c_2}{(r/2)^\theta} \int_0^{\phi(x,r/2)}  \bigg(\frac{\phi(x,r/2)}{t}\bigg)^{\theta/\beta_1} dt = \frac{c_3 \phi(x,r/2)}{r^\theta}.
		\end{align} 		
		 For all $y \in B(x,\delta r/2)$, using \eqref{e:CS-PHR} in the second line below, \eqref{e:CS-PHR-2} in the third, and the monotonicity of $\phi$ and \eqref{e:phi-scale} in the last, we obtain
\begin{align}\label{e:CS-PHR-3}
		&\left|  \vp(x) - \vp(y)\right| \le \kappa \lambda \left| \int_0^\infty e^{-\lambda t} (	P_t^V \1_V (x) - P_t^V \1_V(y)) dt\right|\nn\\
		&\le 2\kappa \lambda  \int_0^{\phi(x,d(x,y))}  dt +   c_1 \kappa \lambda d(x,y)^\theta \int_{\phi(x,d(x,y))}^{\phi(x,r/2)} \frac{dt}{\phi^{-1}(x,t)^{\theta}}  + \frac{c_1 \kappa \lambda d(x,y)^\theta}{(r/2)^\theta} \int_{\phi(x,r/2)}^\infty e^{-\lambda t}dt\nn\\
		&\le c_4    \lambda\phi(x,d(x,y))  +  c_4\lambda  \bigg(\frac{d(x,y)}{r}\bigg)^\theta \phi(x,r/2)  +c_4 \bigg(\frac{d(x,y)}{r}\bigg)^\theta  \nn\\
		&\le c_5   \bigg(\frac{d(x,y)}{r}\bigg)^{\beta_1 \wedge \theta}   \sup_{z \in V_0} \frac{\phi(x,r)}{\phi(z,r)}.
\end{align}		
On the other hand, for all $y \in V \setminus B(x,\delta r/2)$, we have
\begin{align*}
	|\vp (x)-\vp(y)| \le 2 \le 2 (2/\delta)^{\beta_1 \wedge \theta}  \bigg(\frac{d(x,y)}{r}\bigg)^{\beta_1 \wedge \theta}   \sup_{z \in V_0} \frac{\phi(x,r)}{\phi(z,r)}.
\end{align*}
Combining this with \eqref{e:CS-PHR-3}, we arrive at the result. \qed

\section{Weighted  weak elliptic Harnack inequality}\label{s:weHI}

Recall that 	$\NT(u;x_0,r,R)$,  the weighted non-local tail of $u$ with respect to $\p$,  is defined by \eqref{e:def-tail}.
Using $\NT(u;x_0,r,R)$, we introduce  weighted versions of Harnack inequalities.

\begin{defn}\label{d:WEHI}
	\rm (i) We say that \textit{weak elliptic Harnack inequality with weight $\p$} (\WEHIp) holds for  $(\sE,\sF)$, if there exist constants $\delta \in (0,1)$, $\lambda \ge 1$ and $C>0$ such that for all $x_0 \in M$,  $0<(\lambda+2)r<R<R_0$ and 		any $u \in \sB(M)$ which is non-negative and superharmonic in $B(x_0,R)$,
		\begin{align*}
				\bigg(\frac{1}{V(x_0,r)}
				\int_{B(x_0,r)} u^\delta d\mu \bigg)^{1/\delta}\le C\left( \essinf_{B(x_0,r)} u + \phi(x_0,r)\NT(u_-;x_0,\lambda r,R-2r) \right).
	\end{align*}

\noindent (ii) We say that \textit{weak elliptic Harnack inequality with weight $\p$} (\WEHI) holds  if (i) holds for any $u \in \sB(M)$ which is non-negative and harmonic in $B(x_0,R)$.
\end{defn}

\smallskip

Clearly, \WEHIp \ implies \WEHI. The formulation of Harnack inequalities for non-local forms without requiring  additional positivity on the whole space but  instead incorporating a non-local tail term first appeared in  \cite{K1}.

The goal of this  section is to establish the following theorem. 
\begin{thm}\label{t:WEHI}
	Suppose that \VD \ holds. Then we have 
\begin{align*}
	\text{\Tail \, $+$\, 
		\CS \, $+$\, \FK \, $+$\, \PI \; $\Rightarrow$ \ \WEHIp}.
\end{align*}
\end{thm}

 We first establish several lemmas that will be used in the proof of Theorem \ref{t:WEHI}. The proof of Theorem \ref{t:WEHI} will be given at the end of this section.

We start with the following logarithmic lemma, adapting the proof from \cite[Proposition 5.1]{CKW-adv}.

\begin{lem}\label{l:WEHI-1}
Suppose that \VD, \Tail \ and 
\CS \ hold.	Let $x_0 \in M$, $R \in (0,R_0)$ and  $u \in \sB(M)$ be such that non-negative and superharmonic in $B(x_0,R)$. There exists $C>0$ independent  of $x_0,R$ and $u$ such that for all $r \in (0,R/2)$ and $h>0$,
	\begin{align*}
	&\int_{B(x_0,r)} d \Gamma^{(L)}(\log(u+h), \log(u+h)) + 2 \int_{B(x_0,r) \times B(x_0,r)} \left[ \log \frac{u(x)+h}{u(y)+h}\right]^2 J(dx,dy)\\
	& \le \frac{CV(x_0,r)}{\phi(x_0,r)} \left(1 + \frac{1}{h}\phi(x_0,r)\NT(u_-,x_0,2r,R)\right).
	\end{align*}
\end{lem}
\pf Write $B_s:=B(x_0,s)$ for $s>0$. By Lemma \ref{l:Cap} and \eqref{e:phi-scale},   there exists a cutoff function $\vp \in \sF_b$ for $B_r \Subset B_{4r/3}$ such that
\begin{align}\label{e:p-4.12-1}
	\sE(\vp,\vp) \le 2 \Cap (B_r, B_{4r/3}) \le{c_1V(x_0,r)}/{\phi(x_0,r)}.
\end{align}
Fix $h>0$. Note that  $\vp^2/(u+h) \in  \sF^{B_R}$ and $\vp^2/(u+h) \ge 0$ in $B_R$. Thus, by Theorem \ref{t:harmonicity}, 
\begin{align}\label{e:p-4.12-2}
	\sE (u, {\vp^2}/({u+h})) \ge0.
\end{align}
By the inequality in the first display in  \cite[p. 59]{CKW-adv}, 
\begin{align}\label{e:p-4.12-3}
	\int_{B_r}  d \Gamma^{(L)}(\log(u+h), \log(u+h))
	&\le  4\sE^{(L)}(\vp,\vp) -  2 \sE^{(L)}\left(u, \vp^2/(u+h)\right).
\end{align}
Besides, following the arguments for \cite[Proposition 4.12]{CKW-jems}, we obtain
\begin{align}\label{e:p-4.12-4}
		\sE^{(J)}\left(u, \vp^2/(u+h)\right) &= \int_{B_{2r} \times B_{2r}} (u(x)-u(y)) \left( \frac{\vp(x)^2}{u(x)+h}-\frac{\vp(y)^2}{u(y)+h} \right) J(dx,dy)\nn\\
	&\quad + 2 \int_{B_{2r} \times B_{2r}^c} (u(x)-u(y)) \frac{\vp(x)^2}{u(x)+h}\,  J(dx,dy)\nn\\
	&\le  - \int_{B_{r}\times B_{r}} \left[\log \frac{u(y)+h}{u(x)+h}\right]^2 J(dx,dy)+  \int_{B_{2r}\times B_{2r}} (\vp(x)-\vp(y))^2J(dx,dy)\nn\\
	&\quad +  2 \int_{B_{2r}\times B_{2r}^c} \vp(x)^2 J(dx,dy) + \frac{2}{h} \int_{B_{2r}\times B_{R}^c} \vp(x)^2 u_-(y) J(dx,dy).
\end{align}
Since  $\vp=0$ in $B_{2r}^c$, we have $ \int_{B_{2r}\times B_{2r}^c} \vp(x)^2 J(dx,dy) =  \int_{B_{2r}\times B_{2r}^c} (\vp(x)-\vp(y))^2 J(dx,dy).$
On the other hand, since $\vp^2\le 1$, using  \eqref{e:blow-up-integral} and \VD,  we get
\begin{align*}
 \int_{B_{2r}\times B_{R}^c} \vp(x)^2 u_-(y) J(dx,dy)	& \le\NT(u_-,x_0,2r,R) \int_{B_{2r}} \p(x,2r) \mu(dx)\nn\\
&\le c_2V(x_0,r)\NT(u_-,x_0,2r,R) .
\end{align*}
Hence, we deduce from \eqref{e:p-4.12-4} that
\begin{align}\label{e:p-4.12-5}
	&\int_{B_r \times B_r} \left[ \log \frac{u(x)+h}{u(y)+h}\right]^2 J(dx,dy)\le -	\sE^{(J)}\left(u, \vp^2/(u+h)\right) \nn\\
	& \qquad+ \int_{(M\times M) \setminus (B_{2r}\times B_{2r})} (\vp(x)-\vp(y))^2J(dx,dy) +  c_2V(x_0,r)\NT(u_-,x_0,2r,R) /h\nn\\
	&= -	\sE^{(J)}\left(u, \vp^2/(u+h)\right)	 +\sE^{(J)}(\vp,\vp)+  {c_2}V(x_0,r)\NT(u_-,x_0,2r,R)/h.
\end{align}
Combining \eqref{e:p-4.12-1}, \eqref{e:p-4.12-2}, \eqref{e:p-4.12-3} and  \eqref{e:p-4.12-5}, we arrive at
\begin{align*}
	& \int_{B_r} d \Gamma^{(L)}(\log(u+h), \log(u+h)) + 2 \int_{B_r \times B_r} \left[ \log \frac{u(x)+h}{u(y)+h}\right]^2 J(dx,dy) \\
	&\le 4\sE(\vp,\vp)-  2 \sE\left(u, \vp^2/(u+h)\right)  	+   2c_2 V(x_0,r)\NT(u_-,x_0,2r,R) /h\\
	&\le {4c_1V(x_0,r)}/{\phi(x_0,r)} +	2c_2 V(x_0,r)\NT(u_-,x_0,2r,R)/h.
\end{align*}\qed

Recall that $\overline f_D=\mu(D)^{-1}\int_D f d\mu$ for a Borel set $D \subset M$ with $\mu(D) \in (0,\infty)$ and $f \in L^1(D;\mu)$.
\begin{lem}\label{l:WEHI-1+}
Suppose that \VD, \Tail, 
\CS \ 
and \PI \ hold.	Then there exists $C>0$ such that the following  holds:	Let $x_0 \in M$,  $R \in (0,R_0)$  and $u \in \sB(M)$ be such that non-negative and superharmonic in $B(x_0,R)$.  For given constants $a,h>0$ and $b>1$, define
		\begin{align*}
			v=v(a,h,b):=\left[ \log \frac{a+h}{u+h}\right]_+ \wedge \log b.
		\end{align*}
	Then for all $r \in (0,R/(2q_2))$, 
	where $q_2 \ge 1$ is the constant in \PI, it holds that
	\begin{align*}
	\frac{1}{V(x_0,r)} \int_{B(x_0,r)}  (v-\overline{v}_{B(x_0,r)})^2  \,d\mu \le C \left(1 +\frac1{h}\phi(x_0,r) \,\NT(u_-;x_0,2q_2r,R) \right).
	\end{align*}
\end{lem}
\pf Observe that $|v(x)-v(y)| \le |\log  ((u(x)+h)/(u(y)+h))|$ for all $x,y\in M$.  Moreover, by the Markovian property, we  have $	\int_{B(x_0,q_2 r)} d\Gamma^{(L)}(v,v)  \le\int_{B(x_0,q_2 r)} d\Gamma^{(L)}(\log(u+h),\log(u+h))$. Thus, by \PI, we obtain
\begin{align*}
	\int_{B(x_0,r)} (v- \overline{v}_{B(x_0,r)})^2 d\mu &\le c_1\phi(x_0,r) \bigg(\int_{B(x_0,q_2 r)} d\Gamma^{(L)}(\log(u+h),\log(u+h)) \\
	&\qquad \qquad \qquad + \int_{B(x_0,q_2 r) \times B(x_0,q_2 r)} \left[\log \frac{u(x)+h}{u(y)+h}\right]^2 J(dx,dy) \bigg). 
\end{align*}
Combining this with Lemma \ref{l:WEHI-1}, and using \VD \ and  \eqref{e:phi-scale}, we get the desired result. \qed

Following the proof of \cite[Lemma 3.3]{CKW-elp}, we deduce the next lemma from Lemma \ref{l:WEHI-1+}.

\begin{lem}\label{l:WEHI-2}
	 Suppose that \VD, \Tail, 
	 \CS \ and \PI \ hold.	Let  $x_0 \in M$,  $R \in (0,R_0)$, and 
	$u \in \sB(M)$ be such that    non-negative and   superharmonic in $B(x_0,R)$.
	 There exists  $C>0$ independent of $x_0,r$ and $R$ such that the following  holds:	  If
		\begin{align}\label{e:WEHI-2-1}
			\mu( B(x_0,r) \cap \{u < a\}) \le \delta_0 V(x_0,r),
		\end{align}
	for some $a>0$, $\delta_0 \in (0,1)$ and 	$r \in (0,R/(4q_2))$, where  $q_2 \ge 1$ is the constant in \PI,
	 then 
	\begin{align*}
		\frac{\mu\big( B(x_0,2r) \cap \{ u \le  \eps  a-\phi(x_0,r) \NT(u_-;x_0,4q_2r,R)\} \big)}{V(x_0,2r)} \le \frac{C}{(1-\delta_0) |\log \eps|} \quad \text{for all} \;\, \eps \in (0,1).
	\end{align*}
\end{lem}
\pf Write $B_s:=B(x_0,s)$ for $s>0$. Assume that \eqref{e:WEHI-2-1} holds. 
Fix $\eps \in (0,1)$ and define
$$
b:=1/\eps, \quad \; l_0:=\phi(x_0,r) \NT(u_-;x_0,4q_2r,R)  \quad \text{ and } \quad 	v:=\left[ \log \frac{a+l_0}{u+l_0}\right]_+ \wedge  \log b.
$$
By \eqref{e:WEHI-2-1}, we have
\begin{align}\label{e:WEHI-2-2}
	\mu(  B_{2r}  \cap \{ v =0\}) =	\mu(  B_{2r}  \cap \{ u\ge a\}) \ge 	\mu(  B_{r}  \cap \{ u\ge a\})\ge (1-\delta_0) V(x_0,r).
\end{align}
Using \eqref{e:WEHI-2-2} and \VD, we get 
\begin{align}\label{e:WEHI-2-3}
&(\log b)\, \mu( B_{2r}  \cap \{ v = \log b\}) = \int_{   \{y\in B_{2r}:\, v(y) =\log b\} }  v(y) \mu(dy) \nn\\
&= \int_{ \{y \in B_{2r}:\, v(y) =\log b\} }  \int_{  \{  z\in B_{2r}:\,v(z) =0\}} \frac{v(y) - v(z) }{	\mu(  B_{2r}  \cap \{ v =0\}) }  \mu(dz) \mu(dy) \nn\\
&\le \frac{1}{(1-\delta_0) V(x_0,r)}  \int_{ \{y \in B_{2r}:\, v(y) =\log b\} }   \int_{B_{2r}}  (v(y)-v(z)) \,\mu(dz)\mu(dy) \nn\\
&=\frac{c_1}{1-\delta_0}  \int_{ \{y \in B_{2r}:\, v(y) =\log b\} }   (v(y)-\overline v_{B_{2r}}) \mu(dy)  \le \frac{c_1}{1- \delta_0} \int_{ B_{2r}}    \left| v - \overline v_{B_{2r}} \right| d\mu.
\end{align}
Using the Cauchy-Schwarz inequality, Lemma \ref{l:WEHI-1+} and \eqref{e:phi-scale}, we get
\begin{align}\label{e:WEHI-2-4}
	 &\frac{1}{V(x_0,2r)} \int_{ B_{2r}}    \left| v - \overline v_{B_{2r}} \right| d\mu \le \bigg( \frac{1}{V(x_0,2r)} \int_{ B_{2r}}    ( v - \overline v_{B_{2r}})^2 d\mu \bigg)^{1/2}\nn\\
	 &\le  c_2\left( 1 + l_0^{-1} \phi(x_0,4r)\NT(u_-;x_0,4q_2r,R)\right)^{1/2} \le c_3.
\end{align}
Combining \eqref{e:WEHI-2-3} with \eqref{e:WEHI-2-4}, we arrive at
\begin{align*}
	\frac{\mu ( B_{2r}  \cap \{ v =\log b\} )}{V(x_0,2r)} \le \frac{c_1c_3}{(1-\delta_0) \log b}.
\end{align*}
Since $ \{ v = \log b\} = \{ u+l_0 \le b^{-1}(a+l_0)\}  \supset \{ u \le \eps a -l_0 \}$, this yields the desired result. \qed

\begin{lem}\label{l:WEHI-3}
	Suppose that \VD, \Tail, 
	\CS, \FK \ and \PI \ hold.  Let $x_0 \in M$, $R \in (0,R_0)$ and 
	$u \in \sB(M)$ be such that non-negative and superharmonic in $B(x_0,R)$.  Suppose that there are  constants $a>0$, $\delta_0 \in (0,1)$ and $r \in (0,R/(4q_2))$ such that
	\begin{align*}
		\mu( B(x_0,r) \cap \{u < a\}) \le \delta_0 V(x_0,r),
	\end{align*}
 where $q_2 \ge 1$ is the constant in \PI. Then there exists a constant $\eps_0\in(0,1)$ depending on $\delta_0$ but  independent of $x_0,r,R,u$ and $a$ such that 
	\begin{align}\label{e:PLG-claim}
		\essinf_{B(x_0,r)} u \ge \eps_0 a - \phi(x_0,r) \NT(u_-;x_0,4q_2r,R).
	\end{align} 
\end{lem}
\pf  Write $B_s:=B(x_0,s)$ for $s>0$.  Set $l_0:= \phi(x_0,r) \NT(u_-;x_0,4q_2r,R)$.  For a constant $h>0$ to be determined later, let $F_h\in C^2(\R)$ be a function satisfying \eqref{e:condition-F-1} and \eqref{e:condition-F-2}. Then $F_h(u)$ is bounded and subharmonic in $B_{2r}$. Applying \eqref{e:L2mean-diff-form-1} to $F_h(u)$, we get
  \begin{align}\label{e:WEHI-eq-0}
  	\Big(\essinf_{B_{r}} u +h \Big)^{-1} = \esssup_{B_{r}}F_h(u) \le c_1\sI^{\nu_1} ( \sI \vee \sT(\beta_2) )^{1-\nu_1},
  \end{align}
  where 
  \begin{align*}
  	\sI&:=\bigg(  \frac{1}{V(x_0,2r)}\int_{B_{2r}} F_h(u(x))^2\p(x,2r)\mu(dx)  \bigg)^{1/2},\\
  	\sT(\beta_2)&:= \phi(x_0,2r) \sup_{r<t\le 2r,\,0<s\le r/2} \bigg(\frac{s}{2r}\bigg)^{\beta_2}\,\NT(F_h(u);x_0,t,t+s)
  \end{align*}
and $\nu_1:=2\gamma_1\nu/(1+\gamma_1+2\gamma_1\nu)$. Here $\beta_2, \gamma_1$ and $\nu$ are the  constants in \eqref{e:phi-scale}, \eqref{e:blow-up-integral-general} and \FK, respectively. Since $u \ge 0$ in $B_{2r}$,  by \eqref{e:condition-F-1}, we have $F_h(u)\le h^{-1}$ in $B_{2r}$. Thus, by \eqref{e:blow-up-integral},
 \begin{align}\label{e:WEHI-eq-1}
 	\sI \le \frac1h\bigg(  \frac{1}{V(x_0,2r)}\int_{B_{2r}} \p(x,2r)\mu(dx)  \bigg)^{1/2} \le \frac{c_2}{h}.
 \end{align}
 Let $r<t\le 2r$, $0<s \le r/2$ and $x \in B_t$. By using \eqref{e:condition-F-1},   \Tail, the monotonicity of $\p$ and \eqref{e:phi-scale},  since $u \ge 0$ in $B_R$, we get
\begin{align*}
	&\frac{s^{\beta_2} \phi(x_0,2r)}{r^{\beta_2}\p(x_0,t)} \int_{B_{t+s}^c} F_h (u(y)) J(x,dy) \\
	&\le \frac{s^{\beta_2}\phi(x_0,2r)}{r^{\beta_2}\p(x_0,t)}\bigg( \frac{2}{h}\int_{ B_{t+s}^c \cap \{u> -h/2\}}J(x,dy) + \frac{10}{h^2}\int_{B_{t+s}^c \cap \{u\le -h/2\}} u_-(y) J(x,dy)\bigg)\\
	&\le \frac{s^{\beta_2}\phi(x_0,2r)}{r^{\beta_2}\p(x_0,t)}\bigg( \frac{2}{h} \int_{ B(x,s)^c}J(x,dy) +\frac{10}{h^2}\int_{B_R^c} u_-(y) J(x,dy)\bigg)\\
	&\le \frac{c_3s^{\beta_2}\phi(x_0,2r)}{h\phi(x_0,s)}    + \frac{10s^{\beta_2}\phi(x_0,2r)  \NT(u_-;x_0,2r,R)}{h^2r^{\beta_2}}\le c_4\left(\frac1h + \frac{l_0}{h^2}\right).
\end{align*}
Combining this with \eqref{e:WEHI-eq-1}, we get
\begin{align}\label{e:WEHI-eq-2}
	\sI \vee \sT(\beta_2) \le   \frac{c_5}{h}\left(1+ \frac{l_0}{h}\right).
\end{align}
Let $\eps \in (0,1)$ be a constant whose exact value is also to be determined later. By \eqref{e:condition-F-1}, we have
\begin{align*}
	\sI^2 &= \frac{1}{V(x_0,2r)}\int_{B_{2r}  \cap \{ u > \eps a-l_0\}} \frac{\p(x,2r)}{(u(x)+h)^2}\,\mu(dx) + \frac{1}{V(x_0,2r)}\int_{B_{2r} \cap \{ u \le \eps a - l_0\}} \frac{\p(x,2r)}{(u(x)+h)^2}\,\mu(dx) \nn\\
	&=:I_1+I_2.
\end{align*}
 By \eqref{e:blow-up-integral}, we have
\begin{align*}
	I_1 \le \frac{1}{(\eps a -l_0+h)^{2}V(x_0,2r)} \int_{B_{2r}}\p(x,2r)\mu(dx)\le \frac{c_6}{(\eps a -l_0+h)^{2}}.
\end{align*}
On the other hand, since $u \ge 0$ in $B_{2r}$,  using \eqref{e:blow-up-integral0} we see that
\begin{align*}
	I_2&\le  \frac{h^{-2}}{V(x_0,2r)} \int_{B_{2r} \cap \{ u \le \eps a - l_0\}} \p(x,2r)\mu(dx) \\
	&\le  c_7h^{-2}\bigg(\frac{\mu(B_{2r} \cap \{ u \le \eps a - l_0\} )}{V(x_0,2r)}\bigg)^{\gamma_1/(1+\gamma_1)}.
\end{align*}
Thus, by using Lemma \ref{l:WEHI-2}, we get
$
I_2 \le c_8 h^{-2} ((1-\delta_0)|\log \eps|)^{-\gamma_1/(1+\gamma_1)}.
$ Therefore,
\begin{align}\label{e:WEHI-eq-3}
\sI	&\le \frac{c_6}{(\eps a -l_0+h)^{2}} +  \frac{c_8}{h^{2} ((1-\delta_0)|\log \eps|)^{\gamma_1/(1+\gamma_1)}}.
\end{align}
Combining \eqref{e:WEHI-eq-0} with \eqref{e:WEHI-eq-2} and \eqref{e:WEHI-eq-3}, we arrive at
\begin{equation}\label{e:WEHI-eq-4}
	\essinf_{B_{r}} u   \ge c_9 h\left(\frac{1}{((1-\delta_0) |\log \eps| )^{\gamma_1/(1+\gamma_1)}} +\frac{h^2}{(\eps a-l_0+h)^2} \right)^{-\nu_1/2}  \left(1+\frac{l_0}{h}\right)^{\nu_1-1} -h.
\end{equation}

Now we choose a constant $\eps \in (0,1)$ satisfying
\begin{align}\label{e:WEHI-eq-5}
2^{\nu_1/2-1}c_9 ((1-\delta_0) |\log \eps| )^{\nu_1\gamma_1/(2+2\gamma_1)}    \ge 2,
\end{align}
and take
$$
	\eps_0=\frac{\eps}{((1-\delta_0) |\log \eps| )^{\gamma_1/(2+2\gamma_1)}} \quad \text{and} \quad h=\eps_0 a.
$$
 Clearly, \eqref{e:PLG-claim} holds  if $l_0 \ge h$ since $u \ge 0$ in $ B_{r}$. Assume that $l_0 <h$. Then  we get from \eqref{e:WEHI-eq-4} and \eqref{e:WEHI-eq-5} that 
\begin{align*}
		\essinf_{B_{r}} u  &\ge  2^{\nu_1-1}c_9 h\left(\frac{1}{((1-\delta_0) |\log \eps| )^{\gamma_1/(1+\gamma_1)}} +\frac{h^2}{(\eps a)^2} \right)^{-\nu_1/2}   -h\\
		&=2^{\nu_1/2-1}c_9 h\big(((1-\delta_0) |\log \eps| )^{\nu_1\gamma_1/(2+2\gamma_1)} \big)  -h \ge h.
\end{align*}
The proof is complete. \qed

We will use the next Krylov-Safonov type covering lemma in the proof of Theorem \ref{t:WEHI}.   The proof of the next lemma   closely follows that of \cite[Lemma 3.8]{CKW-elp}. We give the proof for completeness.

\begin{lem}\label{l:covering}
	Suppose that \VD \ holds. Let $x_0 \in M$, $r\in (0,R_0/3)$ and $E \subset B(x_0,r)$ be a measurable set. For  $\eps \in (0,1)$, define
	\begin{equation}\label{e:[E]}
	[E]_\eps= \bigcup_{s \in (0,r/2)} \Big\{ B(x,5s) \cap B(x_0,r): x \in B(x_0,r) \text{ and } \mu(E \cap B(x,5s))>\eps V(x,s)\Big\}.
	\end{equation}
Then for each $ \eps \in (0,1)$, either of the following statements holds:
\begin{align*}
	(1)\; [E]_\eps = B(x_0,r) \quad \text{or} \quad (2)\; \mu([E]_\eps) \ge \eps^{-1}\mu(E).
\end{align*}
\end{lem}
\pf Let $\eps \in (0,1)$ and suppose that $B(x_0,r) \setminus [E]_\eps \neq \emptyset$. For   $x \in [E]_{\eps}$, we write $r_x:=4^{-1}\textrm{dist}(x,B(x_0,r) \setminus [E]_\eps).$
  Since $[E]_\eps$ is open, $r_x>0$ for all $x \in [E]_\eps$. Moreover,  for all $x \in [E]_\eps$, we have $r_x<4^{-1}\text{diam}(B(x_0,r)) \le r/2$ and
 \begin{align}\label{e:covering-1}
 	B(x,5r_x) \cap (B(x_0,r) \setminus [E]_\eps) \neq \emptyset.
 \end{align}
By the Besicovitch covering lemma, which can be applicable by \VD \ and the fact that $r_x<R_0/6$, there are countably many pairwise disjoint balls $B(x_i, r_i)$, where $r_i=r_{x_i}$ for $i \ge 1$, such that $[E]_\eps \subset \cup_{i=1}^\infty B(x_i, 5r_i)$. See, e.g., \cite[Theorem 1.16]{He}. Observe that
\begin{align}\label{e:covering-2}
	\mu(E \cap B(x_i, 5r_i)) \le \eps V(x_i,r_i) \quad \text{for all } i \ge 1.
\end{align}
Indeed, if \eqref{e:covering-2} does not hold, then  $B(x_i, 5r_i) \cap B(x_0,r) \subset [E]_\eps$ so that $	B(x_i,5r_i) \cap (B(x_0,r) \setminus [E]_\eps)= \emptyset$. This contradicts \eqref{e:covering-1}. On the other hand, by the Lebesgue differentiation theorem  (see \cite[Theorem 1.8 and Remark 1.13]{He}), we have  $\mu(E\setminus [E]_\eps)=0$ for all $\eps \in (0,1)$. Using this  and \eqref{e:covering-2}, since  $B(x_i,r_i)$, $i \ge 1$, are pairwise disjoint, we arrive at
 \begin{align*}
 	\mu(E) = \mu(E \cap [E]_\eps) \le \sum_{i=1}^\infty \mu(E \cap B(x_i,5r_i)) \le \eps \sum_{i=1}^\infty V(x_i,r_i)= \eps \mu( \cup_{i=1}^\infty B(x_i,r_i)) \le \eps \mu([E]_\eps).
 \end{align*}
The proof is complete. \qed

We now present the proof of Theorem \ref{t:WEHI},   motivated by that of \cite[Theorem 3.1]{CKW-elp}.

\medskip

\noindent \textbf{Proof of Theorem \ref{t:WEHI}.} 
  Let $x_0 \in M$, $R \in (0,R_0)$, $r \in (0,R/(11q_2+2))$ and   $u \in \sB(M)$ be such that non-negative and superharmonic in $B(x_0,R)$. Write $B_s:=B(x_0,s)$ for $s>0$. Set
$$l_1:=\sup_{x \in B_r,  s \in (0,r/2)}\phi(x,5s)\, \NT(u_-;x,20q_2s,R-r).$$ 
By \VD, there exists $c_1 \in (0,1/2)$ such that
\begin{equation}\label{e:WEHI-VD}
	V(x,s) \ge 2c_1 V(x,5s) \quad \text{for all} \;\, x \in B_r \text{ and } 0<s\le r. 
\end{equation}
Let $\eps_0\in (0,1)$ be the constant given in Lemma \ref{l:WEHI-3}  corresponding to the case $\delta_0 = 1-c_1$.

 For each $b>0$,  define
\begin{align*}
	E_n(b):=B_r \cap \left\{ u\ge b\eps_0^n - (1-\eps_0)^{-1}l_1\right\}, \quad n \ge 0.
\end{align*}
We claim that 
\begin{align}\label{e:WEHI-claim}
\mu(E_{n+1}(b)) \ge \mu(	[E_n(b)]_{2^{-1}}) \quad \text{for all $n\ge 0$},
\end{align}
where the set $[E_n(b)]_{2^{-1}}$ is defined by \eqref{e:[E]}. Indeed,  let $n \ge 0$ and define
\begin{align*}
	A_n(b)=\left\{ (x,s) \in B_r \times (0,r/2) : \mu(E_n(b) \cap B(x,5s)) > 2^{-1} V(x,s)\right\}.
\end{align*}
By the definition \eqref{e:[E]}, we have $[E_n(b)]_{2^{-1}}=\cup_{(x,s) \in A_n(b)} (B(x,5s) \cap B_r) \subset \cup_{(x,s)\in A_n(b)} B(x,5s).$ 
Let $\eps>0$. Since $\mu$ is inner regular, there exist finitely many open balls $B(x_i,5s_i)$, $1\le i \le N$, with $(x_i,s_i) \in A_n(b)$ such that
 \begin{align}\label{e:cover-En(b)}
 \mu([E_n(b)]_{2^{-1}}) \le \mu( \cup_{i=1}^N B(x_i,5s_i)) + \eps.
 \end{align} For each $1\le i \le N$,  by the definition of $A_n(b)$ and \eqref{e:WEHI-VD},   we have
\begin{align*}
	\mu\left( E_n(b) \cap B(x_i,5s_i)\right)>2^{-1} V(x_i,s_i) \ge c_1 V(x_i,5s_i).
\end{align*}
Since $u$ is non-negative and superharmonic in $B(x_i,R-r)$ and $5s_i<5r/2<(R-r)/(4q)$, we can apply Lemma \ref{l:WEHI-3} with $a=b\eps_0^n-(1-\eps_0)^{-1}l_1$ to obtain
\begin{align*}
	\essinf_{B(x_i,5s_i)} u &\ge \eps_0 ( b\eps_0^n - (1-\eps_0)^{-1}l_1)  - l_1 = b\eps_0^{n+1} - (1-\eps_0)^{-1}l_1.
\end{align*}
Thus,  $\mu(\cup_{1\le i\le N}B(x_i,5s_i) \setminus E_{n+1}(b))=0$. Using this and \eqref{e:cover-En(b)}, we get
\begin{align*}
\mu(	[E_n(b)]_{2^{-1}}) - \mu(E_{n+1}(b)) 	\le \mu(\cup_{i=1}^N B(x_i, 5s_i) )- \mu(E_{n+1}(b))  +\eps \le \eps.
\end{align*}
Since $\eps>0$ was arbitrarily chosen, we deduce that   \eqref{e:WEHI-claim} holds true.

Let $n_0=n_0(b)\ge 0$ be such that 
\begin{align}\label{e:WEHI-n0}
	2^{-n_0-1}V(x_0,r) < \mu(E_0(b)) \le 2^{-n_0}V(x_0,r).
\end{align}
Then we have
\begin{align}\label{e:WEHI-claim2}
	\mu(E_{n_0}(b)) >  V(x_0,r)/2.
\end{align}
Indeed, if there exist $k \in \{0,\cdots, n_0\}$ such that $\mu(E_k(b)) \ge V(x_0,r)$, then  since $E_k(b) \subset E_{n_0}(k)$, \eqref{e:WEHI-claim2}  holds. Assume that $\mu(E_k(b)) <V(x_0,r)$ for all $0 \le k \le n_0$. Then $\mu([E_k(b)]_{2^{-1}})<V(x_0,r)$ for all $0 \le k <n_0$ by \eqref{e:WEHI-claim}. By Lemma \ref{l:covering}, it follows that $\mu([E_k(b)]_{2^{-1}}) \ge 2 \mu(E_k(b))$ for all $0\le k<n_0$. Thus, by \eqref{e:WEHI-n0}, we get $\mu(E_{n_0}(b)) \ge 2^{n_0} \mu(E_0(b)) > V(x_0,r)/2.$ By \eqref{e:WEHI-claim2}, we can apply Lemma \ref{l:WEHI-3} with $a=b\eps_0^{n_0}-(1-\eps_0)^{-1}l_1$ and $\delta_0=1/2$ to obtain
\begin{align}\label{e:WEHI-last}
	\essinf_{B_{r}} u &\ge c_2 (b\eps_0^{n_0}-(1-\eps_0)^{-1}l_1 )- l_1 \ge c_2b \left( \frac{\mu(E_0(b))}{V(x_0,r)}\right)^{|\log \eps_0|/(\log 2)}-c_3 l_1,
\end{align}
where $c_2$ and $c_3$ are the constants independent of $b$. We used  \eqref{e:WEHI-n0} in the last inequality above.

Set $\delta_1:=(\log 2)/|\log \eps_0|$ and $K:=c_2^{-1}(\essinf_{B_{5r}} u  + c_3l_1)$.
 By \eqref{e:WEHI-last}, we get that for any $\delta \in (0,\delta_1)$,
\begin{align*}
	\frac{1}{V(x_0,r)}\int_{B_r} u^\delta d\mu & = \delta \int_0^\infty b^{\delta-1} \frac{\mu(B_r \cap \{u>b\})}{V(x_0,r)} db\le \delta \int_0^K b^{\delta-1} db + \delta \int_K^\infty b^{\delta-1} \frac{\mu(E_0(b))}{V(x_0,r)} db \\
	&\le \delta \int_0^K b^{\delta-1} db + \delta\int_K^\infty b^{\delta-1} (K/b)^{\delta_1} db =   \frac{\delta_1}{\delta_1-\delta} K^\delta.
\end{align*}
Since $l_1\le c_4\phi(x_0,r)\, \NT(u_-;x_0,11q_2r,R-2r)$ by   \eqref{e:phi2}, this implies that \WEHIp \ holds with $\lambda=11q_2$. The proof is complete. \qed 

\section{Elliptic H\"older regularity}\label{s:eHR}

For a Borel set $D \subset M$ and $u \in L^\infty(D;\mu)$, define
$$\essosc_D u:=\esssup_D u - \essinf_D u.$$

\begin{prop}\label{p:EHR-osc}
Suppose that \VD \ and \WEHI \ hold. Then there exist constants $\theta>0$ and $C>0$ such that for all $x_0 \in M$, $R \in (0,R_0)$, and any $u \in \sB_b(M)$ which is harmonic in $B(x_0,R)$, 
\begin{align}\label{e:EHR-osc}
\essosc_{B(x_0,r)} u \le C  (r/R)^\theta \lVert u \rVert_{\infty} \quad \text{for all} \;\, r \in (0,R].
\end{align} 
\end{prop} 
\pf Without loss of generality, we assume that $\lVert u \rVert_\infty=1$.  Set $\eps:=(2\lambda+2)^{-1}$, where  $\lambda\ge1$  is the constant in the statement of \WEHI, and define $U_n:=B(x_0,\eps^n R)$ for $n\ge 0$.

Let $\theta \in (0,\gamma_0/2)$ be a constant  whose exact value is to be determined later, where $\gamma_0>0$ is the constant in  \eqref{e:blow-up-scale}.  We will construct a non-increasing sequence $(a_n)_{n \ge 0}$ and  a non-decreasing sequence $(b_n)_{n \ge 0}$  such that
\begin{align}\label{e:EHR-construct}
\esssup_{U_n} u \le a_n, \quad \essinf_{U_n} u \ge b_n \quad \text{and} \quad 	a_n-b_n= 2\eps^{n\theta} \quad \text{for all} \;\, n\ge 0.
\end{align}
  \eqref{e:EHR-construct} implies \eqref{e:EHR-osc} by a standard argument. Indeed,  for any $r \in (0,R]$, let $n_0 \ge 0$ be such that $\eps^{n_0+1}<r/R\le \eps^{n_0}$. Then since $B(x_0,r) \subset U_{n_0}$, by \eqref{e:EHR-construct}, we get
\begin{align*}
	\essosc_{B(x_0,r)} u \le a_{n_0} - b_{n_0} = 2\eps^{n_0 \theta}  \le  2\eps^{-\theta} (r/R)^\theta.
\end{align*}

In the following, we construct $(a_n)_{n \ge 0}$ and $(b_n)_{n \ge 0}$ inductively.  
Define $a_0=b_0=1$. Clearly,  \eqref{e:EHR-construct} holds for $n=0$. Suppose that  $(a_n)_{0 \le n \le k}$ and $(b_n)_{0 \le n \le k}$ are constructed to satisfy \eqref{e:EHR-construct} for all $0\le n \le k$. Consider the function
\begin{align*}
	w_k(x):=\eps^{-k\theta }\left(u(x) - (a_k+b_k)/2 \right).
\end{align*}
 By the induction hypothesis,   for all $0 \le j \le k$, since  $b_{k-j}\le b_{k}$ and $a_{k-j} \ge a_{k}$, we have
\begin{equation*}
		\esssup_{U_{k-j}} w_k \le \eps^{-k\theta }\left(a_{k-j} - (a_k+b_k)/2 \right) \le \eps^{-k\theta }\left(a_{k-j} - b_{k-j} -  (a_k-b_k)/2 \right) = 2\eps^{-j\theta} - 1
\end{equation*}
and
\begin{equation}\label{e:EHR-2}
	\essinf_{U_{k-j}} w_k \ge \eps^{-k\theta }\left(b_{k-j} - (a_k+b_k)/2 \right)\ge \eps^{-k\theta }\left(b_{k-j} - a_{k-j} + (a_k-b_k)/2 \right) = -2\eps^{-j\theta} +1.
\end{equation}
We consider the following two cases separately.

\smallskip

Case 1: $\mu(U_{k+1} \cap \{w_k \ge 0\}) \ge 2^{-1} \mu(U_{k+1})$. Define $F=1+w_k$. Then $F$ is harmonic in $U_k$ and is non-negative in $U_k$ by  \eqref{e:EHR-2}. Since $\eps^k R < \eps^{k-1}R/(\lambda+2)$, by applying \WEHI, we get
\begin{align}\label{e:EHR-onestep}
	\bigg(\frac{1}{\mu(U_{k+1})}\int_{U_{k+1}} F^\delta d\mu \bigg)^{1/\delta}\le c_1\left( \essinf_{U_{k+1}} F + \phi(x_0,\eps^{k+1}R)\NT(F_-;x_0,\lambda \eps^{k+1}R,(\eps^k-2\eps^{k+1}) R) \right),
\end{align}
 where the constants $c_1>0$ and $\delta>0$ are independent of $x_0,r,R$ and $u$.   
 Since $F\ge 0$ in $U_{k+1}$ and  $\mu(U_{k+1} \cap \{w_k \ge 0\}) \ge 2^{-1} \mu(U_{k+1})$, we have
\begin{align}\label{e:EHR-main1}
	\bigg(\frac{1}{\mu(U_{k+1})}\int_{U_{k+1}} F^\delta d\mu \bigg)^{1/\delta} \ge \bigg(\frac{\mu( U_{k+1} \cap \{F \ge 1\})}{\mu(U_{k+1})}\bigg)^{1/\delta} \ge 2^{-1/\delta}.
\end{align}
Moreover, since $F\ge 0$ in $U_k \setminus B(x_0,(\eps^k-2\eps^{k+1}) R)$, we  also have that
\begin{align}\label{e:EHR-main2}
&\NT(F_-;x_0,\lambda \eps^{k+1}R,(\eps^k-2\eps^{k+1}) R) = \NT(F_-;x_0,\lambda \eps^{k+1}R,\eps^k R) \nn\\
&\le  \sup_{x \in B(x_0, \lambda \eps^{k+1}R)} \frac{1}{\p(x, \lambda \eps^{k+1}R)} \sum_{j=0}^{k} \int_{U_{k-j-1} \setminus U_{k-j}} F_-(y) J(x,dy),
\end{align}
where $U_{-1}:=M$. Let $x  \in B(x_0, \lambda \eps^{k+1}R)$. For  $0\le j \le k$,  using \eqref{e:EHR-2} and the fact  that  $\lVert u \rVert_\infty =1$ in the first inequality below, \Tail \ in the third,   \eqref{e:blow-up-scale} in the fourth, and \eqref{e:phi2}  in the last,   we get
\begin{align}\label{e:EHR-main3}
	&\int_{U_{k-j-1} \setminus U_{k-j}} F_-(y) J(x,dy) \le 
	2(\eps^{-(j+1)\theta}-1)    J(x,U_{k-j-1} \setminus U_{k-j})  \nn\\
	&	\le 	2(\eps^{-(j+1)\theta}-1)    J(x,B(x, \eps^{k-j}R/2)^c) \le 	c_2\frac{(\eps^{-(j+1)\theta}-1)\,\p(x, \eps^{k-j} R/2)}{\phi(x, \eps^{k-j} R/2)} \nn\\
		& \le c_3	\frac{ (2\lambda \eps^{j+1})^{\gamma_0} \,(\eps^{-(j+1)\theta}-1)\,\p(x, \lambda \eps^{k+1} R)}{\phi(x, \lambda \eps^{k+1} R)} \le c_4\frac{ \eps^{(j+1)\gamma_0} (\eps^{-(j+1)\theta}-1)\,\p(x, \lambda\eps^{k+1} R)}{\phi(x_0, \eps^{k+1} R)}.
\end{align}
Combining \eqref{e:EHR-main2} with  \eqref{e:EHR-main3}, since $\theta<\gamma_0/2$ and   $a^b-1\le b a^b \log a$ for all $a>1$ and $b>0$,  we get
\begin{align}\label{e:EHR-main5}
&\phi(x_0,\eps^{k+1}R)\NT(F_-;x_0,\lambda \eps^{k+1}R,(\eps^k-2\eps^{k+1}) R)\le  c_4\sum_{j=0}^{k} \eps^{(j+1)\gamma_0}(\eps^{-(j+1)\theta}-1)\nn\\
&\le  c_4\theta |\log \eps| \sum_{j=0}^{k} (j+1) \eps^{(j+1)(\gamma_0-\theta)} \le c_4\theta |\log \eps| \sum_{j=0}^{\infty} (j+1) \eps^{(j+1)\gamma_0/2} = c_5(\eps)  \theta.
\end{align}
By \eqref{e:EHR-onestep}, \eqref{e:EHR-main1} and \eqref{e:EHR-main5}, since $\eps^{k\theta}=(a_k-b_k)/2$, we obtain
\begin{align}\label{e:EHR-main6}
&c_1^{-1}2^{-1/\delta} - c_5(\eps)\theta \le  \essinf_{U_{k+1}} F = 2+ \eps^{-k\theta} \Big( \essinf_{U_{k+1}}u - a_k\Big).
\end{align}
By choosing a sufficiently small value for $\theta$ to  satisfy $c_5(\eps)\theta <c_1^{-1} 2^{-1-1/\delta} $ and $ \eps^\theta \ge 1- c_1^{-1}2^{-2-1/\delta}$, we deduce from \eqref{e:EHR-main6} that
\begin{align*}
	 \essinf_{U_{k+1}}u    \ge a_k - \eps^{k\theta}(2-c_1^{-1}2^{-1/\delta} + c_5(\eps)\theta) \ge  a_k - \eps^{k\theta} (2-c_1^{-1}2^{-1-1/\delta} )   \ge  a_k- 2\eps^{(k+1)\theta}.
\end{align*}
Finally, by setting $a_{k+1}=a_{k}$ and $b_{k+1}=a_{k}-2\eps^{k\theta}$, we can conclude that  \eqref{e:EHR-construct} holds for $k+1$.

\smallskip

Case 2: Suppose that  $\mu(U_{k+1} \cap \{w \ge 0\}) < 2^{-1} \mu(U_{k+1})$. Then  $\mu(U_{k+1} \cap \{-w \ge 0\}) \ge 2^{-1} \mu(U_{k+1})$. Using $F=1-w$ instead of $1+w$ and following the arguments above,  one can see that  \eqref{e:EHR-construct} holds for $k+1$ with $a_{k+1}=b_{k} + 2\eps^{k\theta}$ and $b_{k+1}=b_{k}$.

\smallskip

The proof is complete. \qed

\begin{cor}\label{c:EHR}
	Suppose that \VD \ holds. Then \WEHI \ implies \EHR.
\end{cor} 
\pf  Let $x_0 \in M$,  $R \in (0,R_0)$ and $u$ be a bounded harmonic function on $B(x_0,R)$. For almost all $x,y \in B(x_0,R/5)$, since $u$ is harmonic in $B(x,4R/5)$, we get from Proposition \ref{p:EHR-osc} that 
\begin{align*}
|u(x)-u(y)| \le 	\essosc_{B(x,2d(x,y))} u \le c_1 \bigg( \frac{5d(x,y)}{4R}\bigg)^\theta \lVert u \rVert_{\infty},
\end{align*} 
where $\theta>0$ is the constant in Proposition \ref{p:EHR-osc}. Hence, \EHR \ holds.  \qed 

\section{Near-diagonal lower estimates}\label{s:NDL}

In this section, we establish H\"older regularity of the Dirichlet heat kernels and  \NDL. 
Our  approach  follows \cite[Section 4.2]{CKW-jems} which presents similar results for pure jump-type Dirichlet forms when the scale function is independent of the spatial variable,  $R_0=\infty$ and $\p=1$.

We begin with a condition stronger than \NDL, which ensures H\"older regularity for the  heat kernels.
\begin{defn}
	\rm   We say that  \sNDL \ holds if \NDL \ holds, and for all   $x_0 \in M$ and $r\in (0,R_0/q_0)$,  the heat kernel $p^{B(x_0,r)}$ is locally H\"older continuous in $(0,\infty) \times B(x_0,r) \times B(x_0,r)$.
\end{defn}

For a non-empty open set $D \subset M$ and $f\in L^2(D;\mu)$, a function $u \in \sF$ is called a weak solution to the  equation
\begin{align}\label{e:Poisson-equation}
	\sL u = f \quad \text{in } D,
\end{align}
if $\sE(u,\vp)=\la f, \vp \ra$ for all $\vp \in \sF^D$.

Recall that $G^D$ is the Green operator on $D$ defined by \eqref{e:def-Green}. When $\esssup_D 	G^D \1_D<\infty$, the  operator $G^D$  uniquely extends to $L^p(D;\mu)$ for any $p \in [1,\infty]$. See \cite[Lemma 3.2]{GT12}.

\begin{prop}\label{p:weak-solution}{\rm (\!\!\cite[Proposition 6.6]{GH})}.
	Let $D \subset M$ be a non-empty subset and $f \in L^2(D;\mu)$.	
	Suppose that	$\esssup_{x \in D} \E^x \tau_D<\infty$. Then  for any weak solution $u$ to \eqref{e:Poisson-equation}, the function	  $u-G^Df$ is $\sE$-harmonic in $D$.
\end{prop}

\begin{lem}\label{l:l.3.9}
	Suppose \EHR \ and \mEu  \ hold. 
	There exist constants $a_1,a_2>0$ such that the following  holds: For all $x_0 \in M$, $R_1 \in (0,R_0)$ 
	and any bounded measurable function $f$ on 	$D:=B(x_0,R_1)$,
	 if $u\in \sF^{D}$ is a weak solution to \eqref{e:Poisson-equation}, then
	\begin{align}\label{e:EHR-weak-solution}
		\essosc_{B(x_0,r)} u \le  a_1(r/R)^\theta \esssup_{D} |u| +  a_2\phi(x_0,R) \esssup_{B(x_0,R)} |f|    \quad \text{for all} \;\, 0<r \le R \le R_1,
	\end{align}
where $\theta>0$ is the  constant in \EHR.
\end{lem}
\pf Without loss of generality, we assume that  $u$ is bounded. 
Let $0<r \le R\le R_1$ and  $q_3 \ge 1$ be the constant in \mEu. Set $B':=B(x_0,r)$ and  $B:=B(x_0,R)$. We consider the following two cases separately.

\medskip 

Case 1:  $R<R_0/q_3$. Define $v=u-G^Bf$. By \eqref{e:mean-exit-time} and \mEu, we have
\begin{align}\label{e:green-1}	
	\esssup_D |G^B f| \le  \esssup_B |f| \, \esssup_B G^B \1_B \le c_1\phi(x_0,R) \esssup_B |f|.\end{align}
In particular, $v$ is bounded. Since \mEu \ holds, by Proposition \ref{p:weak-solution} and Theorem  \ref{t:harmonicity}, we see that $v$ is harmonic in $B$.  If $r\ge R/2$, then using $\essosc_{B'} u \le 2 \esssup_{D} |u|$ and choosing $a_1$ larger than $2^{1+\theta}$,  we obtain \eqref{e:EHR-weak-solution}. If $r<R/2$, then using  \EHR \ and Remark \ref{r:EHR} in the second inequality below, and \eqref{e:green-1} in the fourth, we get
\begin{align*}
 \essosc_{B'} 	u &\le \essosc_{B'} v  + 	\essosc_{B'} G^B f \le c_2 (r/R)^\theta  \esssup_{D} |v|  + 2\esssup_{B'} |G^B f| \\
  &\le c_2 (r/R)^\theta  \esssup_{D} |u|  +(2^{-\theta}c_2+2)\esssup_{D} |G^B f| \\
   &\le c_2 (r/R)^\theta  \esssup_{D} |u| + c_3\phi(x_0,R) \esssup_B |f|.
\end{align*}

Case 2: $R_0<\infty$ and $R\ge R_0/q_3$.  If $r \ge R_0/(2q_3)$, then 
by taking $a_1$ larger than $2(2q_3m)^\theta$, \eqref{e:EHR-weak-solution}  holds. If $r<R_0/(2q_3)$, then by applying Case 1 with $R=(R_0/(2q_3)) \wedge R_1$, we obtain 
	\begin{align*}
	\essosc_{B'} u &\le  c_4r^\theta \esssup_{D} |u| +  c_5\phi(x_0,R_0/(2q_3)) \esssup_{B} |f|   \\
	 &\le  c_4R_0^\theta (r/R)^\theta \esssup_{D} |u| +  c_5\phi(x_0,R) \esssup_{B} |f|   .
\end{align*}
The proof is complete. \qed

\begin{lem}\label{l:EHR+E=>FK}
	Suppose that \VD, \EHR \ and \mEu \ hold. Then \FK \ holds.
\end{lem}
\pf The proof is essentially due to \cite[Lemma 4.7]{CKW-jems}. Set
$$
\nu:=\frac{\beta_1 \wedge \theta }{2d_2 + \beta_1 \wedge \theta},
$$
where $\beta_1,d_2$ and $\theta$ are positive constants in \eqref{e:phi-scale}, \VD \ and \EHR, respectively. Let $q_3\ge 1$ be the constant in \mEu. By \eqref{e:mean-exit-time} and  \eqref{e:lambda1}, it suffices to show that there exists $c_1>0$  such that  for all $x_0 \in M$, $r \in (0,R_0/(3q_3))$ and any non-empty open set $D \subset B(x_0,r)$,
\begin{equation}\label{e:E-FK-claim}
	\esssup_{D}	G^D\1_D \le c_1 \phi(x_0,r) \bigg(\frac{\mu(D)}{V(x_0,r)} \bigg)^\nu.
\end{equation}

Let $x_0 \in M$, $r \in (0,R_0/(3q_3))$  and $D$ be a non-empty open set with $D\subset B(x_0,r)$. Define
$$
\eps:=\frac13\bigg( \frac{\mu(D)}{V(x_0,r)}\bigg)^{\nu/(\beta_1 \wedge \theta)}.
$$
Since \VD \ holds,  
there exist  $N\ge1$, independent of $x_0$ and $r$, and  $y_1,\cdots, y_N \in B(x_0,r)$ such that $B(x_0,r) \subset \cup_{i=1}^NB(y_i, 3\eps^2 r)$ (see the proof of \cite[Lemma 3.1(ii)]{KS05}).  Fix $i \in \{1,\cdots, N\}$ and define
\begin{align*}
	U_0=B(y_i,3r), \quad U_1=B(y_i,3\eps r), \quad U_2=B(y_i,3\eps^2 r) \quad \text{and} \quad  u=G^{U_0}\1_D - G^{U_1}\1_D.
\end{align*}
Note that $D\subset B(x_0,r) \subset U_0$ and $u$ is non-negative in $M$. By the strong Markov property, $u$ is harmonic in $U_1$. Further, by \mEu \ and \eqref{e:phi2}, we see that
\begin{align}\label{e:u-sup}
	\esssup_{ M}\, u \le  \esssup_{M} E^{U_0} =  \esssup_{U_0} E^{U_0} \le c_1\phi(y_i,3r) \le c_2\phi(x_0,r).
\end{align}
Applying \EHR \ with $u$ and using \VD, the $\mu$-symmetry of $G^{U_0}$ and \eqref{e:u-sup}, we obtain 
\begin{align}\label{e:EHR-FK-1}
	\esssup_{U_{2}} u &\le \essinf_{U_{2}} u + c_3 \eps^{\theta} \esssup_{M} u\nn\\
	&\le \frac{1}{V(y_i,3\eps^2r)} \int_{U_2}\1_{U_0}(y) \, G^{U_0}\1_D(y) \, \mu(dy)  + c_3 \eps^{\theta} \esssup_{M} u\nn\\
	&\le \frac{c_4\eps^{-2d_2}}{V(y_i,3r)} \int_{U_{2}}\1_D(y)  \, E^{U_0}(y)  \, \mu(dy)  + c_3 \eps^{\theta} \esssup_{M} u\nn\\
	&\le \frac{c_4\eps^{-2d_2} \mu(U_2 \cap D)}{V(y_i,3r)} \esssup_{U_2} E^{U_0} + c_3 \eps^{\theta} \esssup_{M} u \nn\\ 
	&\le  c_5\phi(x_0, r) \left(\frac{\eps^{-2d_2}\mu(D)}{V(x_0,r)}+ \eps^\theta \right).
\end{align} On the other hand, since $G^{U_1}\1_M=G^{U_1}\1_{U_1}=E^{U_1}$, by \mEu, \eqref{e:phi-scale} and \eqref{e:phi-comp}, we have
\begin{equation}\label{e:EHR-FK-2}
	\esssup_{x \in U_2} G^{U_1} \1_D \le\esssup_{x \in U_1} E^{U_1}\le c_6 \phi(y_i,3\eps r) \le c_7\eps^{\beta_1} \phi(y_i,r) \le c_8\eps^{\beta_1} \phi(x_0,r)  . 
\end{equation}
Combining \eqref{e:EHR-FK-1} with \eqref{e:EHR-FK-2}, we arrive at
\begin{align*}
	&\esssup_{x \in U_2}	G^D \1_D  \le \esssup_{x \in U_2}	G^{U_0} \1_D \le \esssup_{x \in U_2} u + \esssup_{x \in U_2} G^{U_1} \1_D \nn\\
	&\le  c_9\phi(x_0, r) \left(\frac{\eps^{-2d_2}\mu(D)}{V(x_0,r)}  + \eps^{\beta_1 \wedge \theta} \right)  =c_{10} \phi(x_0,r)\left(   \frac{\mu(D)}{V(x_0,r)}\right)^\nu,
\end{align*}
proving that \eqref{e:E-FK-claim} holds. \qed

\begin{lem}\label{l:HK-osc-weak}
	Suppose that \VD, \EHR \ and \mE \ hold.  Then there exists $q_0 \ge 1$ such that for  all $x_0 \in M$ and $r\in (0,R_0/q_0)$, the heat kernel $p^{B(x_0,r)}$ exists on $(0,\infty)\times (B(x_0,r)\setminus \sN') \times (B(x_0,r)\setminus \sN')$ exists for  a properly exceptional set $\sN'\supset \sN$.   Moreover, there exist  $\nu,C>0$ independent of $x_0$ and $r$  such that the following  hold:
	
	\smallskip

	\noindent	(i) For all $t>0$, 
	\begin{align*}
		\sup_{x \in B(x_0,r)\setminus \sN'}	\esssup_{y \in B(x_0,r)}	|\partial_t p^{B(x_0,r)}(t,x,y)|\le  \frac{C \phi(x_0,r)^{1/\nu}}{V(x_0,r)}\,t^{-1-1/\nu}.
	\end{align*}

	\noindent (ii) 
	For all $0<t\le \phi(x_0,r)$ and  $\eta \in (0,1]$,
	\begin{align*}
		\sup_{x \in B(x_0,\eta \phi^{-1}(x_0,t))\setminus\sN'} \,	\essosc_{y \in B(x_0,\eta \phi^{-1}(x_0,t))} p^{B(x_0, r)} (t,x,y) \le \frac{C \eta^{\beta_1 \theta /(\beta_1+\theta)}\phi(x_0,r)^{1/\nu}}{V(x_0,r)}\,t^{-1/\nu},
	\end{align*}
	where $\beta_1$ and $\theta$ are the constants in \eqref{e:phi-scale} and \EHR \ respectively.
\end{lem}
\pf By Lemma \ref{l:EHR+E=>FK}, the implication (1) $\Rightarrow$ (3) in Proposition \ref{p:FK} and \cite[Corollary 2.7]{GT12}, there exists  $q_0 \ge 1$ such that for any $B:=B(x_0,r)$ with $x_0 \in M$ and $r \in (0,R_0/q_0)$, the  heat kernel $p^{B}:(0,\infty) \times (B\setminus \sN') \times (B\setminus \sN')\to [0,\infty)$ exists  for  a properly exceptional set $\sN' \supset \sN$. Moreover, there exists $c_1>0$ such that $f(t):=	\sup_{x,y \in B\setminus \sN'} p^{B}(t,x,y)$ satisfies that
\begin{align}\label{e:NDL-osc-0}
	f(t) \le \frac{c_1\phi(x_0,r)^{1/\nu}}{V(x_0,r)}t^{-1/\nu} \quad \text{for all } \, t>0.
\end{align}
Observe that $f$ is non-increasing. Indeed, for all $0<t \le t'$, by the semigroup property, we have\begin{align*}	f(t')= 	\sup_{x,y \in B\setminus \sN'} \int_B p^B(t'-t,x,z) p^B(t,z,y) \mu(dz) \le f(t) 	\sup_{x \in B\setminus \sN'} \int_B p^B(t'-t,x,z) \mu(dz)  \le f(t).\end{align*}
For each $x \in B \setminus \sN'$,  define $q_x(t,y)=p^{B}(t,x,y)$.  Following the argument of \cite[Lemma 4.8]{CKW-jems}, we see that  $\partial_t q_x(t,\cdot)$ is well-defined as an element of $L^2(B;\mu)$,
and  for all $t>0$ and a.e. $y \in B \setminus \sN'$,
\begin{align}\label{e:NDL-osc-3}
	|	\partial_t q_x(t,y)| \le  2t^{-1}f(t/2).\end{align}
The result follows from \eqref{e:NDL-osc-0} and \eqref{e:NDL-osc-3}.

(ii)  Let  $0<t\le \phi(x_0,r)$ and  $\eta \in (0,1]$. Set 
 $$ r_1:= \eta \phi^{-1}(x_0,t)^{\beta_1} \quad \text{and} \quad 
r_2:= \eta^{\theta/(\beta_1+\theta)} \phi^{-1}(x_0,t)^{\beta_1}.$$
 Choose any $x \in B(x_0,r_1) \setminus \sN'$. Note that  $q_x(t,\cdot)$ is a weak solution to \eqref{e:Poisson-equation} with $f=\partial_t q_x(t,\cdot)$. By applying Lemma \ref{l:l.3.9} with $u=q_x(t,\cdot)$, $f=\partial_t q_x(t,\cdot)$, $D=B$ and $R=r_2$, and using \eqref{e:NDL-osc-3} and \eqref{e:NDL-osc-0},  we obtain
\begin{align}\label{e:sNDL-1}
&	\essosc_{y \in B(x_0,r_1)} p^B(t,x,y) 	 \le  c_2\bigg(\frac{r_1}{r_2}\bigg)^\theta \esssup_{y\in B} q_x(t,y)+ c_3\phi(x_0, r_2) \esssup_{y \in B(x_0,r_2)} |\partial_t q_x(t,y)| \nn\\
	&\le  c_2 \bigg(\frac{r_1}{r_2}\bigg)^\theta f(t) + \frac{ 2c_3  \phi(x_0,r_2)}{t}f(t/2) \le c_4 \bigg[ \bigg(\frac{r_1}{r_2}\bigg)^\theta + \frac{ \phi(x_0,r_2)}{t}\bigg]\frac{\phi(x_0,r)^{1/\nu}}{V(x_0,r)}t^{-1/\nu}.
\end{align}
By  \eqref{e:phi-scale},  we have
\begin{align*}
		\bigg(\frac{r_1}{r_2}\bigg)^\theta +  \frac{ \phi(x_0,r_2)}{t}&\le   \bigg( \frac{r_1}{r_2}\bigg)^{\theta}  + c_5\bigg(\frac{r_2}{\phi^{-1}(x_0,t)}\bigg)^{\beta_1} = (1+c_5) \eta^{\beta_1\theta/(\beta_1+\theta)}.  
\end{align*}
Combining this with \eqref{e:sNDL-1}, we arrive at the result. 
  \qed

\begin{prop}\label{p:EHR->NDL}
	Suppose that \VD, \EHR \ and \mE \ hold. Then \sNDL \ holds.
\end{prop}
\pf Let $x_0 \in M$ and $R \in (0,R_0/q_0)$ where $q_0\ge 1$ is the constant in Lemma \ref{l:HK-osc-weak}.  By following the arguments in \cite[Lemma 5.13]{GHH18}, using Lemma \ref{l:HK-osc-weak}, we can assume that for all $r \in (0,R]$, the heat kernel
 $p^{B(x_0,r)}$ exists and is locally H\"older continuous.
By \eqref{e:phi2} and Lemma \ref{l:EP},	 there exist $\eps_0 \in (0,1/2)$ and $\delta \in (0,1)$ independent of $x_0$ such that  for all  $r\in (0,R]$,
\begin{align}\label{e:NDL-EP}
&\inf_{y \in B(x_0,r/4)} \int_{B(x_0,r)}  p^{B(x_0,r)}(\eps_0 \phi(x_0,r), y, w)\mu(dw) \nn\\
&\ge \inf_{y \in B(x_0,r/4)} \int_{B(x_0,r)}  p^{B(x_0,r)}(c_1\eps_0 \phi(y,r), y, w)\mu(dw) \ge \delta.
\end{align}
Define $t_r=2\eps_0 \phi(x_0,r)$ for $r \in (0,R]$. For all $r \in (0,R]$, using the  semigroup property,  symmetry,  the Cauchy-Schwarz inequality and \eqref{e:NDL-EP}, we see that
\begin{align}\label{e:NDL-1}
	&\inf_{y \in B(x_0,r/4)}  p^{B(x_0,r)}(t_r, y, y)  =	\inf_{y \in B(x_0,r/4)} \int_{B(x_0,r)} p^{B(x_0,r)}(t_r/2, y,w)^2 \mu(dw)\nn\\
	&\ge  	\inf_{y \in B(x_0,r/4)} \frac{1}{V(x_0,r)}\bigg( \int_{B(x_0,r)} p^{B(x_0,r)}(t_r/2,y,w) \mu(dw)\bigg)^2 \ge \frac{\delta^2}{V(x_0,r)}.
\end{align}
On the other hand, by    Lemma \ref{l:HK-osc-weak}(ii),  there exists $\nu>0$ such that for all $r \in (0,R]$, $\eta \in (0,1]$ and $y,z\in B(x_0,\eta \phi^{-1}(x_0,t_r))$,
\begin{align}\label{e:NDL-2}
 | p^{B(x_0,r)} (t_r,y,y) -  p^{B(x_0,r)} (t_r,y,z) |
&\le \frac{c_2 \eta^{\beta_1\theta/(\beta_1+\theta)} }{(2\eps_0)^{1/\nu}V(x_0,r) }, 
\end{align}
where $\beta_1$ and $\theta$ are the constants in \eqref{e:phi-scale} and \EHR \ respectively.

Take $\eta_2\in (0,1]$ small enough so that  $2c_2\eta_2^{\beta_1\theta/(\beta_1+\theta)}\le (2\eps_0)^{1/\nu}\delta^2$. Let $t \in (0,2 \eps_0 \phi(x_0,R)]$ and $r \in (0,R]$ be such that $t=t_r$.  Using the domain monotonicity of heat kernels (\!\!\cite[Theorem 2.12(b)]{GT12}),  \eqref{e:NDL-1}, \eqref{e:NDL-2}, \VD \ and \eqref{e:phi-scale}, we get that for all $y,z \in B(x_0,\eta_2 \phi^{-1}(x_0,t))$, 
\begin{align*}
	p^{B(x_0,R)}(t,y,z) &\ge p^{B(x_0,r)} (t,y,z) \ge  p^{B(x_0, r)}(t_r, y,y) -  |p^{B(x_0, r)}(t_r, y,y) - p^{B(x_0,r)}(t_r,y,z) |\\
	&\ge  \frac{\delta^2}{2V(x_0,r )} = \frac{\delta^2}{2V(x_0,\phi^{-1}(x_0, t/(2\eps_0) ))} \ge  \frac{c_3}{V(x_0, \phi^{-1}(x_0, t) )}.
\end{align*}
The proof is complete.\qed

\section{Main result and its proof}\label{s:PHR}

In this section, we establish our main result, incorporating additional equivalent hypotheses.
\begin{thm}\label{t:main-2}
	Suppose that \VD, \RVD \ and \Tail \  hold. Then the following equivalences hold:
		\begin{align*}
		\text{\rm \sNDL  \ } 	\Leftrightarrow \  \text{\rm \NDL  \ } 		& \Leftrightarrow \
		\text{\rm \PHR \ + \mE}\\
		&\Leftrightarrow \ 	\text{\rm  \EHR  \ + \mE \   } \\
		&\Leftrightarrow  \	\text{\rm  \WEHIp  \ + \mE \   }  \\
		&\Leftrightarrow \ 	\text{\rm  \WEHI  \ + \mE \   }  \\
		 &\Leftrightarrow  \	\text{\rm  \CS \ + \PI \  }\\
	&\Leftrightarrow  \	\text{\rm \Gcap \ + \PI}.
	\end{align*}
\end{thm}

The proof of Theorem \ref{t:main-2} will be provided at the end of this section.

\smallskip

Under \Tail, one sees that  $X$ satisfies the following L\'evy system formula (see e.g., \cite[Theorem A.3.21]{FOT} and the argument on \cite[p. 40]{CK03}): For any $x \in  M_0$,  non-negative measurable function $f$ on $M \times M_\partial$ vanishing on the diagonal, and  predictable stopping time $\tau$,
\begin{align}\label{e:Levysystem}
	\E^x\bigg[  \sum_{s \le \tau} f( X_{s-}, X_s)  \bigg] = \E^x \left[ \int_0^\tau   \int_{M}  f( X_s,y) J(x,dy)  \, ds \right].
\end{align}

\begin{prop}\label{p:EP}
	Suppose that \VD, \Tail, 
	\CS  \ and \FK \ hold.  There exist constants $b_1>0$ and  $C>0$ such that for all $x_0\in M$, $R \in (0,R_0/q_1)$ and $r \in (0,R/2]$,
	\begin{align*}
		\esssup_{y \in B(x_0,r/2)} \P^y\big(X_{\tau_{B(x_0,r)}} \in B(x_0,R)^c\big) \le C(r/R)^{b_1},
	\end{align*}
 where $q_1 \ge 1$ is the constant in \FK.
\end{prop}
\pf Write $B_s:=B(x_0,s)$ for $s>0$. Fix  $r \in (0,R/2]$ and define $u(y)=	\P^{y}( X_{\tau_{B_r}} \in B_R^c)$. Then $0\le u \le 1$ in $M$ and,  by the strong Markov property,  $u$ is harmonic in $B_r$. By \eqref{e:L2mean-diff-form-2}, we get
\begin{align}\label{e:EP-bound}
	\esssup_{B(x_0, r/2)} u \, \le c_1\sI^{\nu_1},
\end{align}
where $\nu_1>0$ is the constant in \eqref{e:def-nu1} and 
\begin{align*}
	\sI&:=\bigg(  \frac{1}{V(x_0,r)}\int_{B_r} u(x)^2\p(x,r)\mu(dx)  \bigg)^{1/2}.
\end{align*}
Using H\"older inequality and \eqref{e:blow-up-integral-general}, since $0\le u\le 1$, we see that
\begin{align}\label{e:EP-1}
	\sI^2 &\le \frac{1}{V(x_0,r)} \bigg(\int_{B_{r}} u(x)^{2(1+\gamma_1)/\gamma_1}\mu(dx)\bigg)^{\gamma_1/(1+\gamma_1)}  \bigg(\int_{B_{r}} \p(x,r)^{1+\gamma_1}\mu(dx)\bigg)^{1/(1+\gamma_1)}\nn\\
	 &\le \frac{c_2}{V(x_0,r)^{\gamma_1/(1+\gamma_1)}} \bigg(\int_{B_{r}} u(x)\mu(dx)\bigg)^{\gamma_1/(1+\gamma_1)} .
\end{align}

Let $x \in B_r \setminus \sN$. By (1) $\Rightarrow$ (3) in Proposition \ref{p:FK}, the heat kernel $p^{B_r}$ exists. Hence, using \eqref{e:Levysystem} in the first line below, \Tail \ in the third,  \eqref{e:blow-up-scale} in the fourth and \eqref{e:phi-comp} in the last, we get
\begin{align}\label{e:EP-2}
u(x) &
= \E^{x} \bigg[ \int_0^{\tau_{B_r}} J(X_s, B_R^c)ds \bigg] = \int_0^{\infty} \int_{B_r} p^{B_r}(s, x, z) J(z, B_R^c) \,\mu(dz)\, ds \nn\\
& \le \int_0^{\infty} \int_{B_r} p^{B_r}(s, x, z) J(z, B(z,R/2)^c) \,\mu(dz)\, ds\nn\\
&\le  c_3\int_0^{\infty} \int_{B_r} p^{B_r}(s, x, z)\frac{\p(z, R/2)}{\phi(z,R/2)}  \mu(dz)\, ds\nn\\
&\le  c_4 \bigg(\frac{r}{R}\bigg)^{\gamma_0}\int_0^{\infty} \int_{B_r} p^{B_r}(s, x, z)\frac{\p(z,r)}{\phi(z,r)}  \mu(dz)\, ds\nn\\
&\le  \frac{c_5}{\phi(x_0,r)} \bigg(\frac{r}{R}\bigg)^{\gamma_0 }\int_0^{\infty} \int_{B_r} p^{B_r}(s, x, z)\p(z,r) \mu(dz)\, ds.
\end{align}
By Proposition \ref{p:FK}, \SP  \ holds. Using \eqref{e:EP-2} and the symmetry of $p^B$ in the first inequality below, \SP \ in the second and \eqref{e:blow-up-integral} in the third, we obtain
\begin{align*}
&\frac{1}{V(x_0,r)}	\int_{B_r}u(x) \mu(dx) \le  \frac{c_5}{V(x_0,r)\phi(x_0,r)} \bigg(\frac{r}{R}\bigg)^{\gamma_0 }\int_0^{\infty} \int_{B_r} \p(z,r)\int_{B_r}p^{B_r}(s, z, x) \mu(dx)\, \mu(dz)\, ds\nn\\
	&\le  \frac{c_6}{V(x_0,r)\phi(x_0,r)} \bigg(\frac{r}{R}\bigg)^{\gamma_0 }\int_0^{\infty} e^{-c_7s/\phi(x_0,r)}ds \int_{B_r} \p(z,r) \mu(dz) \le c_8\bigg(\frac{r}{R}\bigg)^{\gamma_0 } .
\end{align*}
Combining this with \eqref{e:EP-1}, we obtain the desired result  from  \eqref{e:EP-bound}. \qed

\begin{prop}\label{p:PHR}
	The following implication holds:
	\begin{align*}
	\text{\rm \VD \, $+$\,   \RVD \, $+$\, \Tail \, $+$\, \NDL \; $\Rightarrow$  \ \PHR}.
	\end{align*}
\end{prop}
\pf   Let $x_0 \in M$,  $t_0\ge 0$, $R \in (0,R_0)$ and  $q$ be a bounded measurable function which is caloric in $[t_0,t_0+\phi(x_0,R)] \times B(x_0,R)$. Without loss of generality, we  assume that $t_0=0$ and   $\esssup_{[0, \phi(x_0,R)]\times M }|q|=1$. Let  $\delta \in (0,1/q_1]$ be a constant  to be determined later, where $q_1 \ge 1$ is the constant in \FK.  Fix $(t,x) \in  [\phi(x_0, R) - \phi(x_0,  \delta R),    \phi(x_0,  R)] \times B(x_0, \delta R)$,
and set 
$$S_0:=[0,t] \times M \quad \text{ and } \quad S_n:= [t- \phi(x,  \delta^n R),t]  \times B(x, \delta^nR) \quad \text{for} \;\, n \ge 1.$$
By \eqref{e:phi-scale} and \eqref{e:phi-comp}, we have $\phi(x,\delta R) \le c_1 \phi(x_0,\delta R) \le  c_2\delta^{\beta_1}\phi(x_0,R)$ for some $c_2>1$ independent of $x_0,x$ and $R$. Hence, by assuming $\delta \le c_2^{-1/\beta_1}$, we have  $S_n \subset S_0$ for all $n \ge 1$.  Let $\eps \in (0,1/4)$ be a constant whose  exact value is  to be also determined later.  We will construct a non-increasing sequence $(a_n)_{n \ge 0}$ and  a non-decreasing sequence $(b_n)_{n \ge 0}$  such that 
\begin{align}\label{e:PHR-construct}
	\sup_{S_n} q \le a_n, \quad 	\inf_{S_n} q \ge b_n \quad \text{and} \quad 	a_n-b_n= 2\eps^{n-1} \quad \text{for all} \;\, n\ge 0.
\end{align}
Assuming \eqref{e:PHR-construct} for the moment, we get   \PHR  \   by a standard argument. Indeed, for all $s \in [\phi(x_0,R)-\phi(x_0,\delta R), \phi(x_0,R)]$  and  $y \in B(x_0,\delta R)$ with $s\le t$ and $(s,y) \neq (t,x)$, there exists $n_0\ge 1$ such that $(s,y) \in S_{n_0}$ and $(s,y) \notin S_{n_0+1}$. Since either $s<t- \phi(x,\delta^{n_0+1}R)$ or $y \notin B(x,\delta^{n_0+1}R)$, we deduce from \eqref{e:PHR-construct} that
\begin{align*}
	&\frac{\phi^{-1}(x_0,|s-t|)+d(x,y)}{R} \\
	&\ge \delta^{n_0+1} \ge  \delta  (2/\eps)^{-\log \delta / \log \eps}(a_{n_0}-b_{n_0})^{\log \delta / \log \eps} \ge \delta  (2/\eps)^{-\log \delta / \log \eps} |	q(t,x)  - 	q(s,y) |^{\log \delta/ \log \eps},
\end{align*}
proving that \PHR \ holds.

We now construct monotone sequences $(a_n)_{n \ge 0}$ and $(b_n)_{n \ge 0}$ satisfying \eqref{e:PHR-construct}.  Set $a_0=\eps^{-1}$, $b_0=-\eps^{-1}$, $a_1=1$ and $b_1=-1$.  Then \eqref{e:PHR-construct} holds for $n=0$ and $n=1$. Suppose that \eqref{e:PHR-construct} holds for all $n \le k$ for some $k \ge1$.  Define
\begin{align*}
	S_k^-&:= [t-3\phi(x,\delta_2 \delta^k R),\, t-2\phi(x,\delta_2\delta^k R)] \times B(x,\delta_1 \delta^k R),\\
		E_k^0&:=\{z \in S_k^- :q(z) \le (a_{k-1}+b_{k-1})/2\},
\end{align*}
where $0<\delta_1<\delta_2\le 1/4$ are the constants in Lemma \ref{l:PHR-1}. We consider   two separate cases.

\medskip

Case 1:   $\wt \mu (E^0_k) \ge \wt\mu (S_k^-)/2$. 
Let $E_k$ be a compact subset  of $E^0_k$ such that  $\wt\mu (E_k) \ge \wt\mu (S_k^-)/3$. Fix  $z=(t,y) \in S_{k+1}$. Since $q$ is caloric in $[t_0,t_0+\phi(x_0,R)] \times B(x_0,R)$ and $b_j \le b_k$ for all $0\le j \le k$, by the induction hypothesis, we get that
\begin{align}\label{e:PHR-induction-1}
	q(z) -b_{k} 	&=   \E^{z} \big[q(Z_{\sigma_{E_k} }) -b_{k}  :  \sigma_{E_k} <\wh \tau_{S_{k}}\big] +   \E^{z} \big[q(Z_{\wh \tau_{S_{k}}})-b_{k} : \, \sigma_{E_k}  > \wh \tau_{S_{k}}, \, Z_{\wh \tau_{S_{k}}} \in S_{k-1}\big]\nn\\
	&\quad + \sum_{j=1}^{k-1}  \E^{z} \big[q(Z_{\wh \tau_{S_{k}}})-b_{k} :  \sigma_{E_k}  > \wh \tau_{S_{k}}, \, Z_{\wh \tau_{S_{k}}} \in S_{k-j-1} \setminus S_{k-j}\big]\nn\\
		&\le   \E^{z} \big[q(Z_{\sigma_{E_k} }) -b_{k-1}  :  \sigma_{E_k} <\wh \tau_{S_{k}}\big] +   \E^{z} \big[q(Z_{\wh \tau_{S_{k}}})-b_{k-1} : \, \sigma_{E_k}  > \wh \tau_{S_{k}}, \, Z_{\wh \tau_{S_{k}}} \in S_{k-1}\big]\nn\\
	&\quad + \sum_{j=1}^{k-1}  \E^{z} \big[q(Z_{\wh \tau_{S_{k}}})-b_{k-j-1} :  \sigma_{E_k}  > \wh \tau_{S_{k}}, \, Z_{\wh \tau_{S_{k}}} \in S_{k-j-1} \setminus S_{k-j}\big]\nn\\
	&\le  \frac{(a_{k-1}-b_{k-1})}{2}\P^{z} ( \sigma_{E_k} <\wh \tau_{S_{k}}) +  (a_{k-1}-b_{k-1})\left(1- \P^{z} ( \sigma_{E_k}  < \wh \tau_{S_{k}})\right) \nn\\
	&\quad + \sum_{j=1}^{k-1}  (a_{k-j-1} -b_{k-j-1})\,\P^{z} \big ( Z_{\wh \tau_{S_{k}}} \in S_{k-j-1} \setminus S_{k-j} \big )\nn\\
	&\le \eps^{k-2}\left (2-\P^{z} ( \sigma_{E_k}<\wh \tau_{S_{k}})\right )+  2\sum_{j=1}^{k-1} \eps^{k-j-2}\P^{z} \big ( Z_{\wh \tau_{S_{k}}} \in S_{k-j-1} \setminus S_{k-j} \big ).
\end{align}
By Lemma \ref{l:PHR-1}, \VD \ and \eqref{e:phi-scale}, we have
\begin{align}\label{e:PHR-induction-2}
\P^{z} \big ( \sigma_{E_k}<\wh \tau_{S_{k}}\big ) \ge \frac{c_3\wt \mu(E_k)}{V(x, \delta^k R)\phi(x, \delta^kR)} \ge \frac{c_3\wt \mu(S_k^-)}{3V(x, \delta^k R)\phi(x, \delta^kR)} \ge c_4,
\end{align}
where the constant $c_4 \in (0,1)$ depends only on the constants in  \VD, \NDL \ and  \eqref{e:phi-scale}.  On the other hand, note that   \CS \ and \FK \  hold by 
Lemmas \ref{l:NDL-1}, \ref{l:E->GU} and \ref{l:GU->SCS}. Since $\delta R<R_0/q_1$, by using Proposition \ref{p:EP}, we get that for all $1\le j \le k-1$,
\begin{align}\label{e:PHR-induction-3}
	\P^{z} \big ( Z_{\wh \tau_{S_{k}}} \in S_{k-j-1} \setminus S_{k-j} \big ) \le  \P^{y} \big ( X_{\tau_{B(x, \delta^k R)}} \in B(x, \delta^{k-j}R)^c \big ) \le c_5 \delta^{b_1 j}.
\end{align}
 Now, we take
$$
\eps:=(1- 4^{-1}c_4)^{1/2} \quad \text{and} \quad \delta:=q_1^{-1} \wedge c_2^{-1/\beta_1} \wedge  ((c_4\eps)/(c_4+4c_5))^{1/b_1}.
$$
Combining \eqref{e:PHR-induction-1}, \eqref{e:PHR-induction-2} and  \eqref{e:PHR-induction-3}, we obtain
\begin{align*}
	q(z)& \le b_k + 2\eps^{k-2}(1-c_4/2) +2c_5 \eps^{k-2} \sum_{j=1}^{k-1} (\delta^{b_1}/\eps)^j\\
	& \le b_k + 2\eps^{k-2} (\eps^2 - c_4/4) +2c_5 \eps^{k-2} \sum_{j=1}^\infty (c_4/ (c_4+4c_5))^j = b_k+2\eps^{k}.
\end{align*}
Therefore, \eqref{e:PHR-construct} holds for $n=k+1$ with $b_{k+1}=b_k$ and $a_{k+1}= b_k+2\eps^{k}$.

\smallskip

Case 2: $\wt \mu (E^0_k) < \wt\mu (S_k^-)/2$. Set $E^1_k:=\{z\in S_k^-: -q(z) < -(a_{k-1}+b_{k-1})/2\}$. Then $\wt \mu(E_k^1)   \ge \wt \mu(S_k^-)/2$. Take a compact subset $E'_k$ of $E_k^1$ such that $\wt \mu(E'_k) \ge \wt \mu(S_k^-)/3$. Then \eqref{e:PHR-induction-2}  holds with $E'_k$ instead of $E_k$.  Following the arguments in \eqref{e:PHR-induction-1} for $-q$ instead of $q$, we deduce  that $a_k-q(z)\le 2\eps^k$ for all $z \in S_{k+1}$.
Taking $a_{k+1}=a_k$ and $b_{k+1}=a_k-2\eps^k$,  we obtain \eqref{e:PHR-construct}  for $n=k+1$.

\smallskip

 The proof is complete. \qed

\noindent \textbf{Proof of Theorem \ref{t:main-2}.} Clearly, \sNDL \ implies \NDL. We have
\begin{align*}
	\text{\NDL \,} \Rightarrow \, \text{\mE \ +  \PI}  \qquad &\text{(by Lemma \ref{l:NDL-1})},\\
	\text{\mE \,} \Rightarrow \, \text{\Gcap \,} \Rightarrow \,  \text{\CS} \qquad &\text{(by Lemmas \ref{l:E->GU} and \ref{l:GU->SCS})},
\end{align*}
and
\begin{align*}
\text{\rm \CS \ + \PI \ } &\Rightarrow  \	\text{\rm \CS \ + \FK \  + \PI \ } \qquad \text{(by Corollary \ref{c:PI-FK})} \\
& \Rightarrow \	\text{\rm  \WEHIp  \ + \mE} \qquad  \qquad \quad \;\,\text{(by Theorems \ref{t:main-1} and \ref{t:WEHI})}\\
& \Rightarrow \ 	\text{\rm  \WEHI  \ + \mE \   } \\
& \Rightarrow \ 	\text{\rm  \EHR  \ + \mE}  \qquad  \qquad \qquad \;\;\;\text{(by Corollary \ref{c:EHR})} \\
& \Rightarrow \ \text{\rm  \sNDL}  \qquad  \qquad \qquad \qquad \qquad \text{(by Proposition \ref{p:EHR->NDL})}.
\end{align*}
It remains to prove that \PHR \ + \mE \ is equivalent to other sets of conditions.  By   Proposition \ref{p:PHR}, \NDL \ implies \PHR. 
Since \PHR \ implies \EHR, this completes the proof. \qed

\section{Examples}\label{s:example}

Recall that $\phi_\beta(x,r)=r^\beta$, and  the definition of a $\kappa$-fat open set is given in Definition \ref{d:fat}. 
\begin{prop}\label{ex:reflecting}  Let $\beta \in (0,2)$ and  $D$ be a non-empty $\kappa$-fat open subset of $\R^d$ satisfying \eqref{ex:e:AS} for some $\eta \in (0,d]$.  Consider a symmetric measurable function $\sB$ on $D \times D$ satisfying
\begin{align}\label{e:condition-reflecting}	
	C^{-1} \le \sB(x,y) \le C  \bigg(1+ \frac{|x-y|}{\delta_D(x)}\bigg)^{\theta}  \bigg(1+ \frac{|x-y|}{\delta_D(y)}\bigg)^{\theta}  \quad \text{for all} \;\, x,y \in D,
\end{align}
for some $C>1$ and $\theta \in [0, \beta \wedge \eta)$. Define a bilinear  form $(\sE,\sF)$ on $L^2(D;dx)$  by
		\begin{align*}
		\sE(f,g)&=\int_{D \times D} (f(x)-f(y))(g(x)-g(y)) \frac{\sB(x,y)}{|x-y|^{d+\beta}} dxdy,\qquad 
		 \sF=\overline{ \text{\rm Lip}_c(\overline D)}^{\sE_1}.
		\end{align*}
 Then $(\sE,\sF)$ is a regular Dirichlet form satisfying {\rm NDL}$(\phi_\beta)$ and {\rm PHR}$(\phi_\beta)$ with $R_0=\diam(D)$.
\end{prop}
  \pf
  By \eqref{e:d-set}, $D$ satisfies \VD \ and \RVD. Define
\begin{align}\label{e:blow-up-example-boundary}
\p(x,r):=\bigg(1+ \frac{r}{\delta_D(x)}\bigg)^\theta, \quad x \in D, \, r>0.
\end{align}
By Proposition \ref{ex:Hardy-boundary}, $\p$   is an admissible weight function for $(\sE, \phi_\beta)$  with $R_0= \diam(D)$. By \eqref{e:condition-reflecting} and Lemma \ref{l:Ju-Tail}, $(\sE,\sF)$ satisfies TJ$^\p_\le(\phi_\beta)$. Since $\beta<2$,  by Proposition \ref{ex:beta<2},  $(\sE,\sF)$ is a regular Dirichlet form satisfying  CSU$^\p(\phi_\beta)$.

 Consider a bilinear form $(\sE^{(\beta)}, \sF^{(\beta)})$ on $L^2(\overline D;dx)$  defined by
 \begin{align*}
 	\sE^{(\beta)}(f,g)&:= \int_{D \times D} \frac{(f(x)-f(y))(g(x)-g(y))}{ |x-y|^{d+\beta}} dxdy, \quad		\sF^{(\beta)}:=\big\{f \in L^2(D;dx): \sE^{(\beta)}(f,f)<\infty\big\}.
 \end{align*} According to \cite[Theorem 1.1]{CK03}, $(\sE^{(\beta)}, \sF^{(\beta)})$ is a regular Dirichlet form satisfying condition HK$(\phi_\beta)$ in \cite{CKW-jems}. Thus, by \cite[Corollary 1.3 and Theorem 1.18]{CKW-jems},  PI$(\phi_\beta)$ holds for $(\sE^{(\beta)}, \sF^{(\beta)})$. By \eqref{e:condition-reflecting}, it follows that
  PI$(\phi_\beta)$ holds  for $(\sE,\sF)$. The desired result now follows from Theorem \ref{t:main-2}. \qed

 \begin{prop}\label{ex:reflecting-mixed}
 	Let $\beta \in (0,2)$,  $D$ be a non-empty $\kappa$-fat open subset of $\R^d$ satisfying \eqref{ex:e:AS} for some $\eta \in (0,d]$, and $\sB$ be a symmetric measurable function  on $D \times D$ satisfying \eqref{e:condition-reflecting}
 	for some $C>1$ and $\theta \in [0, 2 \wedge \eta)$. Let $(a_{ij}(x))_{1\le i,j\le d}$
 	be a symmetric $d\times d$-matrix valued measurable function on $D$ satisfying the uniform ellipticity condition: There exists $\Lambda>1$ such that 
 	\begin{align}\label{ex:reflecting-mixed-uniform-ellipticity}
 		\Lambda^{-1}|\xi|^2 \le 	\sum_{i,j=1}^d a_{ij}(x) \xi^i \xi^j \le \Lambda |\xi|^2 \quad \text{for all $x \in D$ and $\xi=(\xi^1,\cdots,\xi^d)\in \R^d$}.
 	\end{align}	 Define a bilinear  form $(\wt\sE,\wt\sF)$ on $L^2(\overline D;dx)$  by
 	\begin{align*}
 		\wt\sE(f,g)&=\sum_{i,j=1}^d\int_D a_{ij}(x) \partial_i f(x)  \partial_j g(x)  dx +\int_{D \times D} (f(x)-f(y))(g(x)-g(y)) \frac{\sB(x,y)}{|x-y|^{d+\beta}} dxdy,\\
 		\wt	\sF&=\overline{ \text{\rm Lip}_c(\overline D)}^{\sE_1}.
 	\end{align*}
 	Then $(\wt\sE,\wt\sF)$ is a regular Dirichlet form satisfying {\rm NDL}$(\phi_2)$ and {\rm PHR}$(\phi_2)$  with $R_0=1 \wedge \diam(D)$.
 \end{prop}
 \pf Since  $\wt \sE(f,f)<\infty$ for all $f \in  \text{\rm Lip}_c(D)$, $(\wt\sE,\wt\sF)$ is a regular Dirichlet form. Define $\p$ by \eqref{e:blow-up-example-boundary} and set $R_0:=1\wedge \diam(D)$. By Proposition \ref{ex:Hardy-boundary}, $\p$   is an admissible weight function for $(\sE, \phi_2)$ with $R_0$.   Further, by \eqref{ex:reflecting-mixed-uniform-ellipticity}, PI$(\phi_2)$ holds for $(\wt\sE,\wt\sF)$. 
 Let $x_0 \in M$ and $0<r \le R$ with $R+2r<R_0$.  
 Define a cutoff function $\vp$ for $B(x_0,R) \Subset B(x_0,R+r)$ by \eqref{e:Lipschitz-cutoff}. Then $\vp \in  \text{\rm Lip}_c(\overline D)\in \wt \sF$ and  $|\nabla \vp| \le 1/r$. 
 By Lemma \ref{l:Ju-Tail}, TJ$^\p_\le(\phi_\beta)$ holds for $(\wt \sE,\wt\sF)$. Hence, since $\beta<2$, using Lemma \ref{l:Tail-integral}, we get
 \begin{equation}\label{e:Lipschitz-cutoff-1}
 	\int_D (\vp(x)-\vp(y))^2 \frac{\sB(x,y)}{|x-y|^{d+\beta}}  dy\le \int_D \bigg( 1 \wedge \frac{d(x,y)}{r}\bigg)^2 \frac{\sB(x,y)}{|x-y|^{d+\beta}}  dy\le \frac{c_1\p(x,r)}{r^\beta} \quad \text{for all} \;\, x\in D.
 \end{equation}
 Using  \eqref{ex:reflecting-mixed-uniform-ellipticity} and \eqref{e:Lipschitz-cutoff-1}, since   $|\nabla \vp |\le 1/r$  and $\p \ge 1$, we get that for all $f \in \wt\sF'_b$,
 \begin{align*}
 	&\int_{B(x_0,R+2r)} f^2\, d\Gamma(\varphi,\varphi) \\
 	&\le \Lambda \int_{B(x_0,R+r)} f(x)^2 |\nabla\vp(x)|^2 dx +\int_{B(x_0,R+2r) \times D} f(x)^2 (\vp(x)-\vp(y))^2  \frac{\sB(x,y)}{|x-y|^{d+\beta}} dxdy
 	\\
 	&\le \frac{c_2}{r^2} \int_{B(x_0,R+r)} f(x)^2  dx +\frac{c_1}{r^\beta}\int_{B(x_0,R+2r)} f(x)^2 \p(x,r) dx\le \frac{c_1+c_2}{r^2} \int_{B(x_0,R+r)} f(x)^2\p(x,r)  dx.
 \end{align*}
 Thus, CSU$^\p(\phi_2)$ holds for $(\wt\sE,\wt\sF)$. The  result now follows from Theorem \ref{t:main-2}. \qed

  \begin{example}
\rm (Traces of isotropic $\beta$-stable process on the exterior of bounded $C^{1,1}$ open sets). Let $d\ge 1$, $\beta \in (0,2)$ and $Y$ be an isotropic $\beta$-stable process on $\R$ with generator $\Delta^{\beta/2}$.
Let  $D\subset \R^d$ be  the complement of a bounded $C^{1,1}$ open set and $X$ be the trace process of $Y$ on $D$. That is, $X_t=Y_{\tau_t}$ for $t \ge 0$, where $\tau_t:=\inf\{s>0:A_s>t\}$ and $A_t:=\int_0^t \1_{\{Y_s \in \overline D \} }$. By \cite[Theorem 6.2.1]{FOT}, $X$ is a Hunt process associated with a  regular Dirichlet form $(\sE,\sF)$ given by
\begin{equation}\label{e:def-trace-process}
	\sE(f,g)=\frac12\int_{D\times D} (f(x)-f(y))(g(x)-g(y))J_D(x,y)dxdy,  \qquad  \sF=\overline{\text{Lip}_c(\overline D)}^{\sE_1},
\end{equation}
where
\begin{align}\label{e:def-trace-process-jump-kernel}
J_D(x,y)= \frac{\sA(d, -\beta)}{|x-y|^{d+\beta}} + \int_{D^c \times D^c} \frac{\sA(d, -\beta)}{|x-z|^{d+\beta}} G_{D^c}(z,w) \frac{\sA(d, -\beta)}{|w-y|^{d+\beta}}dzdw,
\end{align}
$G_{D^c}(z,w)$ is the Green function of $Y$ on $D^c$ and 
$
\sA(d, -\beta)=\frac{\beta 2^{\beta-1}\Gamma((d+\beta)/2)}{\pi^{d/2}\Gamma(1-\beta)}.
$
Note that $D$ satisfies \eqref{ex:e:AS} with $\eta=1$. 
Further, by \cite[Theorem 2.6]{BGPR} and \cite[Proposition 3.6]{GH24}, we have for all $x,y \in D$,
\begin{align*}
&\int_{D^c \times D^c} \frac{\sA(d, -\beta)}{|x-z|^{d+\beta}} G_{D^c}(z,w) \frac{\sA(d, -\beta)}{|w-y|^{d+\beta}}dzdw\\
 &\le \frac{c_1}{\delta_D(x)^{\beta/2} (1+\delta_D(x))^{\beta/2} \delta_D(y)^{\beta/2} (1+\delta_D(y))^{\beta/2} (|x-y|+\delta_D(x) + \delta_D(y) + \delta_D(x)\delta_D(y))^d}\\
	&\le c_1\bigg(\frac{|x-y|^2}{\delta_D(x)\delta_D(y)} \bigg)^{\beta/2} \frac{1}{|x-y|^{d+\beta}}.
\end{align*}
Hence, $(\sE,\sF)$ satisfies \eqref{e:condition-reflecting} with $\theta=\beta/2 <\beta \wedge 1$. Applying Proposition \ref{ex:reflecting}, we deduce that {\rm NDL}$(\phi_\beta)$ and {\rm PHR}$(\phi_\beta)$ hold for $(\sE,\sF)$ with $R_0=\infty$.
\end{example}

\begin{example}
	\rm (Traces of isotropic $\beta$-stable process on the half-space). Let $d\ge 3$, $\beta \in (0,2)$ and $Y$ be an isotropic $\beta$-stable process on $\R$ with generator $\Delta^{\beta/2}$.
	Let   $X$ be the trace process of $Y$ on the half-space $\mathbb H:=\{(\wt x, x_d)\in \R^d: x_d>0\}$.  Then  $X$ is a Hunt process associated with  the regular Dirichlet form $(\sE,\sF)$ given by \eqref{e:def-trace-process} and \eqref{e:def-trace-process-jump-kernel} with $D=\mathbb H$. Clearly, $\mathbb H$ satisfies \eqref{ex:e:AS} with $\eta=1$. 
	Moreover, by \cite[Theorem 6.1]{BGPR}, it holds that for all $x,y \in \mathbb H$,
	\begin{align*}
		&\int_{\mathbb H^c \times \mathbb H^c} \frac{\sA(d, -\beta)}{|x-z|^{d+\beta}} G_{\mathbb H^c}(z,w) \frac{\sA(d, -\beta)}{|w-y|^{d+\beta}}dzdw\\
		&\le \frac{c_1}{\delta_{\mathbb H}(x)^{\beta/2}\delta_{\mathbb H}(y)^{\beta/2} (|x-y|+\delta_{\mathbb H}(x) + \delta_{\mathbb H}(x))^d}\le c_1\bigg(\frac{|x-y|^2}{\delta_{\mathbb H}(x)\delta_{\mathbb H}(y)} \bigg)^{\beta/2} \frac{1}{|x-y|^{d+\beta}}.
	\end{align*}
	Hence, $(\sE,\sF)$ satisfies \eqref{e:condition-reflecting} with $\theta=\beta/2 <\beta \wedge 1$. By Proposition \ref{ex:reflecting}, we conclude that {\rm NDL}$(\phi_\beta)$ and {\rm PHR}$(\phi_\beta)$ hold for $(\sE,\sF)$ with $R_0=\infty$.
\end{example}

\begin{example}
	\rm (Dirichlet forms with non-local Neumann conditions). Let $d\ge 1$, $\beta \in (0,2)$ and $D\subset \R^d$ be a Lipschitz domain. Then $D$ satisfies \eqref{ex:e:AS} with $\eta=1$. Consider a Dirichlet form $(\sE,\sF)$ on $L^2(D;dx)$ defined by
	\begin{equation*}
		\sE(f,g)=\frac12\int_{D\times D} (f(x)-f(y))(g(x)-g(y))J^D(x,y)dxdy,  \qquad  \sF=\overline{\text{Lip}_c(\overline D)}^{\sE_1},
	\end{equation*}
	where
	\begin{align*}
		J^D(x,y)= \frac{\sA(d, -\beta)}{|x-y|^{d+\beta}} +\bigg( \int_D \frac{dw}{|z-w|^{d+\beta}} \bigg)^{-1} \int_{D^c} \frac{\sA(d, -\beta)}{|x-z|^{d+\beta}|z-y|^{d+\beta}}  dz.
	\end{align*}
	The Dirichlet form $(\sE,\sF)$ arises in the study of non-local Neumann problem introduced in \cite{DRV17}. For further details, we also refer to \cite{Ab20, Von21}. 	 By  \cite[Proposition 2.1]{AFR23}, we have for all $x,y \in D$,
	\begin{align*}
		J^D(x,y)\le  \frac{c_1}{|x-y|^{d+\beta}} \log\bigg( e + \frac{|x-y|}{\delta_D(x) \wedge \delta_D(y)}\bigg).
	\end{align*}	Hence, since $\sup_{r>0}r^{\beta/2}\log(e+r)<\infty$,  it holds that for all $x,y \in D$,
		\begin{align*}
		J^D(x,y)\le   \frac{c_2}{|x-y|^{d+\beta}} \bigg( e + \frac{|x-y|}{\delta_D(x) \wedge \delta_D(y)}\bigg)^{\beta/2}\le   \frac{c_2}{|x-y|^{d+\beta}} \bigg( e + \frac{|x-y|}{\delta_D(x)}\bigg)^{\beta/2}\bigg( e + \frac{|x-y|}{\delta_D(y)}\bigg)^{\beta/2},
	\end{align*}
	proving that $(\sE,\sF)$ satisfies \eqref{e:condition-reflecting} with $\theta=\beta/2 $. By  Proposition \ref{ex:reflecting},
	 {\rm NDL}$(\phi_\beta)$ and {\rm PHR}$(\phi_\beta)$ hold for $(\sE,\sF)$ with $R_0=\diam(D)$. Note that elliptic regularity  for $(\sE,\sF)$  was established in \cite{AFR23}.
	 
	 We next  consider a Dirichlet form $(\wt \sE, \wt \sF)$ defined by
	 	\begin{equation*}
	 	\wt\sE(f,g)=\frac{a}{2}\int_{D} \nabla f \cdot \nabla g \, dx + \sE(f,g),   \qquad  \wt\sF=\overline{\text{Lip}_c(\overline D)}^{\sE_1},
	 \end{equation*}
	 where $a>0$ is a constant. This Dirichlet form arises in the study of mixed local and non-local dispersal models with a certain Neumann condition. See \cite{DPV23,DV21}. By  Proposition \ref{ex:reflecting-mixed}, 
	 {\rm NDL}$(\phi_2)$ and {\rm PHR}$(\phi_2)$ hold for $(\sE,\sF)$ with $R_0=1\wedge \diam(D)$.
\end{example}

 \begin{prop}\label{ex:Hardy}
 	Suppose that \VD \ and \RVD \ hold with $R_0=\infty$, and \Exi \ holds. 	 Let $o \in M$,  $\beta \in (0, \alpha \wedge d_1)$ and $\lambda \in (0, \beta)$, where $d_1$ is the constant in \RVD.  Consider a bilinear form
 \begin{align*}
 		\sE(f,g)&=\int_{M \times M} (f(x)-f(y))(g(x)-g(y)) J(x,y) \mu(dx)\mu(dy), \\
 		\sF&=\left\{ f \in L^2(M;\mu):\sE(f,f)<\infty\right\},
 \end{align*}
 		where  $J$ is a symmetric measurable function on $M \times M$  satisfying
 		\begin{align}\label{e:condition-Hardy}
 		\frac{	C^{-1}}{V(x,d(x,y))d(x,y)^{\beta}}\le J(x,y) \le 	\frac{	C}{V(x,d(x,y))d(x,y)^{\beta}}\bigg(1 + \frac{d(x,y)^{2\lambda}}{d(x,o)^{\lambda}d(y,o)^{\lambda}} \bigg)   
 		\end{align}
 	for all $x,y \in M$, with some  $C>1$.  Then $(\sE, \sF)$ is a regular Dirichlet form on $L^2(M;\mu)$ satisfying
 		{\rm NDL}$(\phi_\beta)$ and {\rm PHR}$(\phi_\beta)$
 		 with $R_0=\infty$.
 \end{prop}
\pf  Define the Dirichlet form  $(\sE^\beta,\sF^\beta)$  as in \eqref{e:def-sE-beta}. By \eqref{e:CS-Hardy-domain} and Proposition \ref{p:subordinate-Exi}, since $(\sE^\beta,\sF^\beta)$ is regular, we get
\begin{align*}
	\sF&= \left\{ f \in L^2(M;\mu):
	W_{\beta/2}(f)
<\infty\right\}= \overline{\left\{ f \in C_c(M):
	W_{\beta/2}(f)
	<\infty\right\}}^{\sE_1^\beta}  \\
& = \overline{\left\{ f \in C_c(M):
	W_{\beta/2}(f)
	<\infty\right\}}^{\sE_1} = \overline{\left\{ f \in C_c(M):
	\sE(f,f)
	<\infty\right\}}^{\sE_1},
\end{align*}
proving that $(\sE,\sF)$ is a regular Dirichlet form on $L^2(M;\mu)$. Here $W_{\beta/2}(f)$ denotes the Sobolev norm defined in \eqref{e:def-Sobolev}. \bk  By  Propositions \ref{ex:Hardy-origin-2} and  \ref{ex:CS-Hardy}, there exists a constant $\theta \in [\lambda, \beta  )$ such that 
$\p(x,r):=\left(1+r/d(x,o)\right)^\theta$ is an admissible weight function for $(\sE, \phi_\beta)$ and CSU$^\p(\phi_\beta)$ holds for $(\sE, \sF)$.  Further, since $\theta \ge \lambda$, by \eqref{e:condition-Hardy}, we see that J$^\p_\le(\phi_\beta)$ holds for $(\sE,\sF)$. Hence, by Lemma \ref{l:Ju-Tail}, TJ$^\p_\le(\phi_\beta)$ holds  for $(\sE,\sF)$.  Recall from  Proposition \ref{p:subordinate-Exi} that $(\sE^\beta,\sF^\beta)$  satisfies J$_\le(\phi_\beta)$ and PI$(\phi_\beta)$.  Thus,  since J$_\ge(\phi_\beta)$  holds for $(\sE,\sF)$ by \eqref{e:condition-Hardy}, we deduce that  PI$(\phi_\beta)$ holds for $(\sE,\sF)$. By Theorem \ref{t:main-2}, we conclude that 	{\rm NDL}$(\phi_\beta)$ and {\rm PHR}$(\phi_\beta)$
 hold for $(\sE,\sF)$. \qed

\begin{example}
	\rm 
Suppose that \VD \ and \RVD \ hold with $d_2,d_1>0$ and $R_0=\infty$, and \Exi \ holds. 	 Let $o \in M$, $\beta \in (0, \alpha \wedge d_1)$ and $\lambda \in (0, \beta)$. Consider a symmetric Dirichlet form 
	\begin{align*}
		\sE(f,g)&= \frac12 \int_{M\times M}  \frac{(f(x)-f(y))(g(x)-g(y))}{V(x,y)d(x,y)^{\beta}} \mu(dx)\mu(dy)\nn\\
		&\quad\;  +  \frac12 \int_{M\times M}  \frac{(f(x)-f(y))(g(x)-g(y))}{V(x,y)d(x,y)^{\beta-2\lambda} d(x,o)^{\lambda}d(y,o)^{\lambda}} \mu(dx)\mu(dy),\\
		\sF&=\left\{ f\in L^2(M;\mu): \sE(f,f)<\infty \right\},
	\end{align*}
	where $V(x,y):=V(x,d(x,y)) + V(y,d(x,y))$ for $x,y \in M$.
By Proposition \ref{ex:Hardy}, $(\sE,\sF)$ is a regular Dirichlet form on $L^2(M;\mu)$ satisfying {\rm NDL}$(\phi_\beta)$ and {\rm PHR}$(\phi_\beta)$  with $R_0=\infty$. Note that the (formal) generator associated with $(\sE,\sF)$ is of the form
	\begin{align*}
		Lf(x) &= \lim_{\eps\to 0} \int_{y\in M: d(x,y)>\eps} \frac{f(y)-f(x)}{V(x,y)d(x,y)^{\beta}}\mu(dy)\nn\\
		&\quad  +\lim_{\eps\to 0} \int_{y\in M: d(x,y)>\eps} \frac{f(y)-f(x)}{V(x,y)d(x,y)^{\beta-2\lambda} d(x,o)^\lambda d(y,o)^\lambda}\mu(dy),
	\end{align*}
	which serves as a non-local analogue of the Laplacian with singular drift,  $\Delta + |x|^{-\lambda} \nabla$. 
\end{example}

\section{Appendix}\label{s:Appendix}

\subsection{Proofs of Proposition \ref{p:FK} and Corollary \ref{c:PI-FK}}\label{s:A}
We introduce two additional  Nash-type inequalities for the regular Dirichlet form $(\sE,\sF)$.

\begin{defn}
	\rm 
	\noindent	(i) We say that  \textit{Kigami-Nash inequality} \KN \ holds if there exists $C>0$ such that for all $x_0 \in M$,  $r>0$, $s\in (0,R_0)$  and $f \in \sF^{B(x_0,r)}\cap L^\infty(M;\mu)$,
	\begin{align*}
		\lVert f \rVert_{2}^{2} \le C \bigg( \sup_{z \in B(x_0,r)} \phi(z,s)  \sE(f,f) + \frac{\lVert f \rVert_{1}^2 }{\inf_{z \in B(x_0,r)} V(z,s)}   \bigg).
	\end{align*}

	\noindent (ii) 	We say that   \textit{local Nash inequality} \LN \ holds if there exist  $C,\nu>0$   such that for all $x_0 \in M$, $r \in (0,R_0/2)$  and $f \in \sF^{B(x_0,r)}\cap L^\infty(M;\mu)$,
	\begin{align}\label{e:LN}
		\lVert f \rVert_{2}^{2+2\nu}\, \lVert f \rVert_{1}^{-2\nu}\le C \frac{\phi(x_0,r) }{V(x_0,r)^\nu} \left(  \sE(f,f)  +	\frac{\lVert f \rVert_{2}^2}{\phi(x_0,R_0)} \right).
	\end{align} 
\end{defn}

The proof of the next lemma is standard. See, e.g. \cite[Theorem 2.1]{Sa}.

\begin{lem}\label{l:PI-KN}
	Suppose that \VD \ holds. Then \PI \ implies \KN.
\end{lem}
\pf Let $x_0 \in M$, $r>0$, $s \in (0,R_0)$ and $f \in \sF^{B(x_0,r)} \cap L^\infty(M;\mu)$.  Define
\begin{align*}	f_s(x)=\frac{1}{V(x,s/(4q_2))}\int_{B(x,s/(4q_2))} f(z) \mu(dz), \quad x \in M,
\end{align*}
where $q_2\ge 1$ is the constant in \PI. By the Cauchy-Schwarz inequality,
\begin{align}\label{e:KN-eq-0}
	\lVert f \rVert_2^2 \le 	2\int_{B(x_0,r)} |f(x)-f_s(x)|^2 \mu(dx) +2\lVert f_s \rVert_2^2 .
\end{align}

 By the Fubini's theorem  and \eqref{e:VD2}, we have 
\begin{align}\label{e:KN-eq-1}
	\lVert f_s \rVert_{1} 	&\le c_1 \int_M  \frac{|f(z)|}{V(z,s/(4q_2))}\int_{B(z,s/(4q_2))}   \mu(dx) \mu(dz)= c_1 \lVert f \rVert_{1}.
\end{align} 
Note that $f_s(x)=0$ if $B(x,s/(4q_2)) \cap  B(x_0, r )=\emptyset$. If there exists $y \in B(x,s/(4q_2)) \cap  B(x_0, r )$, then by \eqref{e:VD2}, $V(x,s/(4q_2)) \ge  c_2V(y,s/2) \ge c_3 V(y,s) \ge c_3\inf_{z \in B(x_0,r)} V(z,s).$ 
Thus, we get
\begin{align}\label{e:KN-eq-2}
	\lVert f_s \rVert_{\infty} \le \frac{\lVert f \rVert_1}{c_3\inf_{z \in B(x_0,r)} V(z,s) }.
\end{align}
Combining \eqref{e:KN-eq-1} with \eqref{e:KN-eq-2}, by using H\"older inequality, we get
\begin{align}\label{e:KN-eq-3}	
	\lVert f_s \rVert_{2}^2 \le \lVert f_s \rVert_{1}\,\lVert f_s \rVert_{\infty} \le \frac{c_4\lVert f \rVert_1^2}{\inf_{z \in B(x_0,r)} V(z,s) }.
\end{align}

Next, we estimate the integral $\int_{B(x_0,r)} |f(x)-f_s(x)|^2 \mu(dx)$.
By the standard covering theorem (see \cite[Theorem 1.2]{He}) and compactness of $\overline{B(x_0,r)}$, 
 there exists a family of  finite pairwise disjoint open balls $\sV=\{B(y_i,s/(20q_2) ) : y_i \in \overline{B(x_0,r)}, \, 1\le i\le N\}$ such that  $
 \overline{B(x_0,r)} \subset \cup_{i=1}^N	B(y_i,s/(4q_2) )$. 
 Define  $U_i:=B(y_i,s/(4q_2))$,   $U_i':=B(y_i,s/(2q_2))$ and $W_i:=B(y_i,s/2)$ 
 for $1\le i\le N$.  We assume, without loss of generality, that $U_i \cap B(x_0,r) \neq \emptyset$ for all $1\le i \le N$. 
  Observe that there exists $N_0\ge1$ independent of $x_0,r$ and $s$ such that \begin{align}\label{e:KN-eq-covering}
	\sum_{i=1}^N \1_{W_i}(z)\le N_0 \quad \text{in}\;\, M. 
	\end{align}
Indeed, since  $\sV$ is pairwise disjoint, using \VD, we see that  for all $z \in M$,\begin{align*}	
		\sum_{i=1}^N \1_{W_i}(z) &
		\le \sum_{1\le i \le N: z \in W_i} \frac{V(y_i, s)}{V(z,s/2)} \le \frac{c_5}{V(z,s)}\sum_{1\le i\le N: z \in W_i} V(y_i, s/(20q_2))\nn\\	& \le \frac{c_5}{V(z,s)} \mu \left( \cup_{1\le i\le N: y_i \in B(z,s/2)} B(y_i, s/(20q_2)) \right) \le  \frac{c_5V(z,(10q_2+1)s/ (20q_2))}{V(z,s)} \le c_5.
	\end{align*}
Since $\{U_i:1\le i\le N\}$ is an open covering of $B(x_0,r)$, \bk we have
\begin{align}\label{e:KN-eq-4}
	&\int_{B(x_0,r)} |f(x)-f_s(x)|^2 \mu(dx) \le \sum_{i=1}^N \int_{U_i}  |f(x)-f_s(x)|^2 \mu(dx)\nn\\
	&\le 2\sum_{i =1}^N \int_{U_i}  |f(x)-\overline f_{U_i'}|^2  \mu(dx) + 2\sum_{i = 1}^N \int_{U_i} |f_s(x)-\overline f_{U_i'}|^2  \mu(dx).
\end{align} 
By the Cauchy-Schwarz inequality and \eqref{e:VD2}, we  get that for all $1\le i\le N$,
\begin{align}\label{e:KN-eq-5}
	&\int_{U_i}  |f_s(x)-\overline f_{U_i'}|^2  \mu(dx) \le 	\int_{U_i} \frac{1}{V(x,s/(4q_2))} \int_{B(x,s/(4q_2))}  |f(z)-\overline f_{U_i'}|^2  \mu(dz) \mu(dx) \nn\\
	&\le \frac{c_6}{V(y_i,s/(4q_2))}\int_{U_i}   \mu(dx) \int_{U_i'}  |f(z)-\overline f_{U_i'}|^2  \mu(dz)  = c_6 \int_{U_i'}  |f(z)-\overline f_{U_i'}|^2  \mu(dz).   
\end{align}
 Further, for all $1\le i\le N$,  since there is  $z_i \in U_i \cap B(x_0,r)$, \bk we get from \eqref{e:phi2} that
\begin{align}\label{e:KN-eq-4-2}
	\phi(y_i,s) \le c_6 \phi(z_i,s) \le c_7\sup_{z \in B(x_0,r)} \phi(z,s).
\end{align}
Using \PI, \eqref{e:KN-eq-covering} and \eqref{e:KN-eq-4-2}, we obtain
\begin{align}\label{e:KN-eq-6}
	&\sum_{i = 1}^N \int_{U_i}  |f(x)-\overline f_{U_i'}|^2  \mu(dx)  \le \sum_{i = 1}^N \int_{U_i'}  |f(x)-\overline f_{U_i'}|^2  \mu(dx) \nn\\
	&\le 	c_8\sum_{i = 1}^N \phi(y_i, s)\bigg(\int_{W_i} d\Gamma^{(L)}(f,f) + \int_{W_i \times W_i} (f(x) - f(y))^2 J(dx,dy) \bigg)\nn\\
	&\le 	c_8N_0\max_{1\le i\le N} \phi(y_i, s)\bigg(\int_{M} d\Gamma^{(L)}(f,f) + \int_{M \times M} (f(x) - f(y))^2 J(dx,dy) \bigg)\nn\\
	&\le c_7c_8  N_0 \sup_{z \in B(x_0,r)} \phi(z,s)  \sE(f,f).
\end{align}
By  \eqref{e:KN-eq-4}, \eqref{e:KN-eq-5} and \eqref{e:KN-eq-6}, we deduce that
\begin{align}\label{e:KN-eq-7}
	\int_{B(x_0,r)} |f(x)-f_s(x)|^2 \mu(dx) \le  c_9N_0\sup_{z \in B(x_0,r)} \phi(z,s)  \sE(f,f).
\end{align}
Combining \eqref{e:KN-eq-0}, \eqref{e:KN-eq-3} and \eqref{e:KN-eq-7}, we conclude that \KN \ holds.  \qed 

The proofs of the next two lemmas are originally due to \cite[Propositions 3.1.4 and 3.4.1]{BCS} where the case of $R_0=\diam(M)=\infty$ and   $\phi(x,r)=r^2$ was considered.

\begin{lem}\label{l:KN-LN}
	Suppose that \VD \ and \RVD \ hold. Then \KN \ implies \LN.
\end{lem}
\pf Let $x_0 \in M$, $r \in (0,R_0/2)$ and  $f \in \sF^{B(x_0,r)}\cap L^\infty(M;\mu)$. Without loss of generality, we assume that $\lVert f \rVert_1=1$. Define
$$H(f)=\sE(f,f)   + \frac{\lVert f \rVert_2^2}{\phi(x_0,R_0)}.$$
 By \KN, there exists $c_1\ge 1$ such that for all $s>0$,
\begin{align}\label{e:KN-global}
	\lVert f \rVert_{2}^{2}&\le  c_1 \1_{\{s<R_0\}}\bigg( \sup_{z \in B(x_0,r)} \phi(z,s)  \sE(f,f) + \frac{1 }{\inf_{z \in B(x_0,r)} V(z,s)}   \bigg) +  \1_{\{s\ge R_0\}} \sup_{z \in B(x_0,r)} \frac{\phi(z,s)	\lVert f \rVert_{2}^{2} }{\phi(x_0,R_0)} \nn\\
	&\le c_1 \bigg(  \sup_{z \in B(x_0,r)} \phi(z,s) H(f) + \frac{\1_{\{s<R_0\}}}{\inf_{z \in B(x_0,r)} V(z,s)}   \bigg).
\end{align} 
Set $\nu:=\beta_1/d_2$ and $\nu':=\beta_2/d_1$ where $d_1,d_2$ and $\beta_1,\beta_2$ are the constants in \VD \ and \eqref{e:phi-scale}, respectively.  We consider the following three cases separately.

\smallskip

Case 1: Assume $\phi(x_0,r)H(f)V(x_0,r)>1$. 
 By \eqref{e:KN-global}, \eqref{e:phi-scale}, \eqref{e:phi-comp} and \VD,  for all $s\le r$,
\begin{align*}
	\lVert f \rVert_{2}^{2} &\le  c_2 \bigg(  \bigg(\frac{s}{r}\bigg)^{\beta_1}  \sup_{z \in B(x_0,r)} \phi(z,r) H(f) + \bigg(\frac{2r}{s}\bigg)^{d_2}\frac{1}{\inf_{z \in B(x_0,r)} V(z,2r)}   \bigg)\nn\\
	&\le c_3 \bigg(  \bigg(\frac{s}{r}\bigg)^{\beta_1}  \phi(x_0,r) H(f) + \bigg(\frac{r}{s}\bigg)^{d_2}\frac{1}{V(x_0,r)}   \bigg).
\end{align*} 
Taking 
$s=(\phi(x_0,r)H(f)V(x_0,r))^{-1/(d_2+\beta_1)}r$ in the above, we obtain
$$
\lVert f \rVert_{2}^{2+2\nu} \le c_4\bigg(\frac{\phi(x_0,r)^{d_2/(d_2+\beta_1)}H(f)^{d_2/(d_2+\beta_1)}}{V(x_0,r)^{\beta_1/(d_2+\beta_1)}} \bigg)^{1+\nu} =\frac{c_4\phi(x_0,r)H(f)}{V(x_0,r)^\nu}.
$$

Case 2: Assume  $\phi(x_0,r)H(f)V(x_0,r)\le 1$ and  $(\phi(x_0,r)H(f)V(x_0,r))^{-1/(d_1+\beta_2)}r<R_0$. By H\"older inequality,  
\begin{align}\label{e:KN-case2-1}
	\lVert f \rVert_{2}^{2+2\nu} \le  V(x_0,r)^{\nu'-\nu}\lVert f \rVert_{2}^{2+2\nu'}.
\end{align}
Using \eqref{e:KN-global}, \eqref{e:phi2} and \eqref{e:VD2},  we see that for all $r\le s<R_0$, 
\begin{align*}	\lVert f \rVert_{2}^{2+2\nu'} &\le  c_5 \bigg(  \bigg(\frac{s}{r}\bigg)^{\beta_2}  \phi(x_0,r) H(f) + \bigg(\frac{r}{s}\bigg)^{d_1}\frac{1}{V(x_0,r)}   \bigg)^{1+\nu'}.\end{align*} 
Letting $s=(\phi(x_0,r)H(f)V(x_0,r))^{-1/(d_1+\beta_2)}r$ in the above, we obtain$$\lVert f \rVert_{2}^{2+2\nu'}  \le \frac{c_6\phi(x_0,r)H(f)}{V(x_0,r)^{\nu'}}.$$
Combining this with \eqref{e:KN-case2-1}, we get \eqref{e:LN}.

Case 3: Assume  $\phi(x_0,r)H(f)V(x_0,r)\le 1$ and  $(\phi(x_0,r)H(f)V(x_0,r))^{-1/(d_1+\beta_2)}r\ge R_0$.  Note that \eqref{e:KN-case2-1} is still valid. Using \eqref{e:KN-case2-1}, \eqref{e:KN-global} with $s=(\phi(x_0,r)H(f)V(x_0,r))^{-1/(d_1+\beta_2)}r$ and \eqref{e:phi2},  we obtain
\begin{align*}&	\lVert f \rVert_{2}^{2+2\nu} \le  V(x_0,r)^{\nu'-\nu}\lVert f \rVert_{2}^{2+2\nu'}\nn\\
	& \le  c_7V(x_0,r)^{\nu'-\nu} \bigg(  \bigg(\frac{(\phi(x_0,r)H(f)V(x_0,r))^{-1/(d_1+\beta_2)}r}{r}\bigg)^{\beta_2}  \phi(x_0,r) H(f) \bigg)^{1+\nu'} = \frac{c_7\phi(x_0,r)H(f)}{V(x_0,r)^{\nu}}.\end{align*}
  The proof is complete.  \qed

\begin{lem}\label{l:LN-LN} 
	Suppose that \VD \ and \RVD \ hold. Then \LN \ implies \HLN.
\end{lem}
\pf Let  $q \ge 1$ be a constant  to be determined later, $x_0 \in M$, $r \in (0,R_0/(2q))$ 
and  $f \in \sF^{B(x_0,r)} \cap L^\infty(M;\mu)$. Applying \eqref{e:LN} with $f$ 
and $r$ replaced by $qr$, by \RVD \ and \eqref{e:phi-scale}, we get
\begin{equation}\label{e:LN-HLN}
	\lVert f \rVert_{2}^{2+2\nu}\, \lVert f \rVert_{1}^{-2\nu} \le \frac{c_1\phi(x_0,qr) }{V(x_0,qr)^\nu} \left(  \sE(f,f)  +	\frac{\lVert f \rVert_{2}^2}{\phi(x_0,qr)} \right) \le \frac{c_2q^{\beta_2-d_1\nu}\phi(x_0,r) }{V(x_0,r)^\nu}  \sE(f,f) + \frac{c_3q^{-d_1\nu}\lVert f \rVert_{2}^2}{V(x_0,r)^\nu}.
\end{equation} 
Set $q:=(2c_3)^{1/(d_1\nu)}$. By H\"older inequality, we have
\begin{align*}
	\frac{c_3q^{-d_1\nu}\lVert f \rVert_{2}^2}{V(x_0,r)^\nu} \le 	\frac{c_3q^{-d_1\nu}\lVert f \rVert_{2}^{2}}{V(x_0,r)^\nu} \bigg( \frac{V(x_0,r)\lVert f\rVert_2^2 }{\lVert f \rVert_{1}^{2} }\bigg)^\nu = \frac{1}{2}\lVert f \rVert_{2}^{2+2\nu}\, \lVert f \rVert_{1}^{-2\nu} .
\end{align*}
Combining this with \eqref{e:LN-HLN}, we conclude that \HLN \ holds with $q_1=2q$. \qed

We now present the proofs of Proposition \ref{p:FK} and Corollary \ref{c:PI-FK}.

\medskip

\noindent \textbf{Proof of Proposition \ref{p:FK}.}
 Choose any $x_0\in M$ and $r \in (0,R_0/q_1)$. Write $B:=B(x_0,r)$.

\smallskip

(1) $\Rightarrow$ (2) $\Rightarrow$ (3): These implications can be proved  similarly to \cite[Lemmas 5.4 and 5.5]{GH} by choosing  $a=cV(x_0,r)^\nu/\phi(x_0,r)$ in their respective proofs.

\smallskip

(3) $\Rightarrow$ (1): Let $D \subset B$ be a non-empty open set and  $\lambda:=\phi(x_0,r)(\mu(D)/V(x_0,r))^\nu$.
Using \eqref{e:FK-Dirichlet}, we see that for all $a>1$, 
\begin{equation}\label{e:FK-Dirichlet-1}
\esssup_{y\in D} \P^y(\tau_D>a\lambda) =\esssup_{y \in D}\int_D p^{B}(a\lambda,y,z)\,\mu(dz) \le \frac{c_1\mu(D)}{V(x_0,r)} \bigg(\frac{\phi(x_0,r)}{a\lambda} \bigg)^{1/\nu} = \frac{c_1}{a^{1/\nu}}.
\end{equation} 
Set $a_0:=(2c_1)^\nu$.  By the Markov property and \eqref{e:FK-Dirichlet-1}, we get that for all $m \ge 1$,
\begin{align}\label{e:FK-Dirichlet-2}
		\esssup_{y \in D}\,	\P^y(\tau_D>m a_0\lambda) &\le  \Big(	\esssup_{y \in D}  \P^y(\tau_D>a_0\lambda) \Big)^m  \le 2^{-m}.\end{align}
Hence, for a.e. $y \in U$, 
\begin{align*}
	\E^y \tau_D=\sum_{m=0}^\infty \int_{ma_0\lambda}^{(m+1)a_0\lambda}  \P^y(\tau_D>t)dt \le a_0\lambda \sum_{m=0}^\infty   \P^y(\tau_D>ma_0\lambda)  \le a_0\lambda \sum_{m=0}^\infty   2^{-m} =2a_0\lambda.
\end{align*}
By \eqref{e:lambda1}, it follows that $\lambda_1(D) \ge 1/(2a_0\lambda)$. 

\smallskip

Clearly, \SP \ implies \mEu. To finish the proof, we assume (3) and show that \SP \ holds. Let $x_0 \in M$ and $r \in (0,R_0/q_1)$. Taking $D=B(x_0,r)$ and $\lambda=\phi(x_0,r)$ in \eqref{e:FK-Dirichlet-2}, we get that 
\begin{align}\label{e:FK-Dirichlet-3}
	\esssup_{y \in B(x_0,r)} \P^y (\tau_{B(x_0,r)} >m a_0 \phi(x_0,r)) \le 2^{-m} \quad \text{for all} \;\, m \ge 1. 
\end{align}
Fix  $t>0$ and let $m_0\ge 1$ be such that $m_0a_0 \le t/\phi(x_0,r) < (m_0+1)a_0$. By \eqref{e:FK-Dirichlet-3}, we get
\begin{align*}
	\esssup_{y \in B(x_0,r)} \P^y (\tau_{B(x_0,r)} >t) \le 	\esssup_{y \in B(x_0,r)} \P^y (\tau_{B(x_0,r)} >m_0 a_0 \phi(x_0,r)) \le 2^{-m_0} \le e   \exp \left(-\frac{(\log 2 )\, t}{a_0\phi(x_0,r)}\right).
\end{align*}
Thus, \SP \ holds. The proof is complete. \qed

\noindent \textbf{Proof of Corollary \ref{c:PI-FK}.} 
Combining Lemmas \ref{l:PI-KN}, \ref{l:KN-LN}, \ref{l:LN-LN} and Proposition \ref{p:FK},  we get	\begin{align*}
	\text{\PI \; $\Rightarrow$ \ \KN  \; $\Rightarrow$ \  \LN  \; $\Rightarrow$ \  \HLN  \; $\Leftrightarrow$ \ \FK}.
\end{align*}
\qed

\subsection{Proofs of Lemmas \ref{l:PHR-1}--\ref{l:NDL-1}}\label{s:A2}

In this subsection, we let $\eta_1,\eta_2 \in (0,1)$ and $q_0 \ge1$ be the constants in \NDL.

\medskip

\noindent \textbf{Proof of Lemma \ref{l:PHR-1}.}   By \eqref{e:phi-scale}, there is $\delta_2 \in (0,1/4]$  such that  
$
3\phi(x,\delta_2 r) \le  \eta_1\phi(x, r/q_0)
$  for all $x\in M$ and $r>0$.
Let $\delta_1:=\eta_2 \delta_2$, $x_0 \in M$, $r \in (0,R_0)$, $t \ge \phi(x_0,r)$ and $E \subset [t-3\phi(x_0,\delta_2 r), t-  2\phi(x_0,\delta_2 r)] \times B(x_0, \delta_1 r)$ be a compact subset. Note that for all $s \in [t-\phi(x_0,\delta_2 r),t]$,
$$
3\phi(x_0,\delta_2 r)- (t-s)  \le 3\phi(x_0,\delta_2 r) \le \eta_1 \phi(x_0,r/q_0)  \quad \text{and} \quad 2\phi(x_0,\delta_2 r)- (t-s)  \ge \phi(x_0, \delta_1r/\eta_2).
$$
Denote $\tau=\wh \tau_{[t-\phi(x_0,r), t] \times B(x_0,r)}$ and  $E_a=\{x\in B(x_0, \delta_1r):(a,x)\in E\}$ for $a>0$. Observe that for all $(s,y) \in [t-\phi(x_0,\delta_2 r),t] \times B(x_0,\delta_1 r)$, 
\begin{align}\label{e:NDL-PHR-1}
	3 \phi(x_0,\delta_2 r)\,	\P^{(s,y)}(\sigma_E<\tau) &= \P^{(s,y)} \left(\int_0^\tau \1_E(s-a, X_a) \,da>0 \right) \int_0^{3 \phi(x_0,\delta_2 r)} dw\nn\\
	&\ge \int_0^{3 \phi(x_0,\delta_2 r)} \P^{(s,y)} \left(\int_0^\tau \1_E(s-a, X^{B(x_0,r/q_0)}_a) \,da>w \right)dw\nn\\
	& =  \E^{(s,y)} \left[ \int_0^\tau \1_E(s-a,X^{B(x_0,r/q_0)}_a)da \right]\nn\\
	&=\int_{2\phi(x_0,\delta_2 r)- (t-s)}^{3\phi(x_0,\delta_2 r)- (t-s)} \int_{E_{s-a}} p^{B(x_0,r/q_0)}(a, y, z) \mu(dz)\, da.
\end{align}
By \NDL, \VD \ and \eqref{e:phi-scale}, we have for all $s \in [t-\phi(x_0,\delta_2 r),t]$ and $y \in B(x_0,\delta_1r)$,
\begin{align}\label{e:NDL-PHR-2}
	&\int_{2\phi(x_0,\delta_2 r)- (t-s)}^{3\phi(x_0,\delta_2 r)- (t-s)} \int_{E_{s-a}} p^{B(x_0,r/q_0)}(a, y, z) \,\mu(dz) da\ge \int_{2\phi(x_0,\delta_2 r)- (t-s)}^{3\phi(x_0,\delta_2 r)- (t-s)} \frac{c_1 \mu(E_{s-a})}{V(x_0, \phi^{-1}(x_0, a))}  da \nn\\
	&\ge \frac{c_2}{V(x_0,r)}  \int_{2\phi(x_0,\delta_2 r)- (t-s)}^{3\phi(x_0,\delta_2 r)- (t-s)} \mu(E_{s-a}) da = \frac{c_2}{V(x_0,r)}  \int_{t-3\phi(x_0,\delta_2 r)}^{t-2\phi(x_0,\delta_2 r)} \mu(E_a) da =  \frac{c_2 \wt\mu(E)}{V(x_0,r)}. \qquad
\end{align}
Combining \eqref{e:NDL-PHR-1} with \eqref{e:NDL-PHR-2} and using \eqref{e:phi-scale}, we arrive at the desired result. \qed

\medskip

\noindent \textbf{Proof of Lemma \ref{l:NDL-1}.}   By \eqref{e:phi-scale}, there is $q_2 \ge q_0$ such that
\begin{align}\label{e:NDL-PI}
	\phi(x, q_2r) \ge \eta_1^{-1}\phi(x,q_0r/\eta_2) \quad \text{for all} \;\, x \in M \text{ and } r>0.
\end{align}

(i) For \PI, we follow the proof of \cite[Proposition 3.5(i)]{CKW-jems}, which  is originally due to \cite{KS87}. Let $x_0 \in M$ and $r \in (0,R_0/q_2)$. Set $V:=B(x_0,r)$ and $U:=B(x_0,q_2r/q_0)$. Consider a bilinear form
\begin{align*}
	\overline \sE(f,g)&:=\int_{U} d\Gamma^{(L)}(f,g) + \int_{U \times U} (f(y) - f(x)) (g(y) - g(x)) J(dx,dy),\\
	\overline \sF&:=\{ f \in L^2(U;\mu): \overline \sE(f,f)<\infty\}.
\end{align*}
Let $\overline \sF^U$ be the closure of $\overline \sF \cap C_c(U)$ in $L^2(U;\mu)$. Following the argument in the  proof of  \cite[Proposition 3.5(i)]{CKW-jems}, one sees that the Hunt process $\overline X^U$ associated with $(\overline \sE, \overline \sF^U)$ has a transition density function $\overline p^U(t,x,y)$ with respect to $\mu$, 
\begin{align}\label{e:PI-1}
	\overline p^U(t,x,y) \ge p^U(t,x,y) \quad \text{for all} \;\, t>0 \text{ and a.e. }x,y \in U,
\end{align}
and 
\begin{align}\label{e:PI-2}
	\overline \sE(f,f) 	&\ge \frac{1}{2} \int_{U \times U} \frac{\overline p^U(\phi(x_0,r/\eta_2),x,y)}{\phi(x_0,r/\eta_2)}(f(x)-f(y))^2 \mu(dx)\mu(dy) \quad \text{for all}\;\, f \in \overline \sF.
\end{align}
By  \eqref{e:NDL-PI}, we get from \NDL \ that  $p^U(\phi(x_0,r/\eta_2),x,y) \ge c_1/V(x_0,r/\eta_2)$ for all $x,y \in V$. By \eqref{e:PI-1} and \VD, it follows that  for a.e. $x,y \in V$,
\begin{align*}
		\overline p^U(\phi(x_0,r/\eta_2),x,y) \ge \frac{c_1}{V(x_0,r/\eta_2)} \ge  \frac{c_2}{V(x_0,r)}.
\end{align*}
Hence, using  \eqref{e:phi-scale} and \VD,  we get from \eqref{e:PI-2} that for all $ f\in   \sF_b$, 
\begin{align*}
	\phi(x_0,r)	\,\overline \sE(f,f) &\ge c_3 \int_{V \times V} \overline p^U(\phi(x_0,r/\eta_2),x,y)(f(x)-f(y))^2 \mu(dx)\mu(dy) \\ 
	&\ge  \frac{c_4}{V(x_0, r)} \int_{V \times V}  (f(x)-f(y))^2 \mu(dx)\mu(dy)\ge  c_4\inf_{a \in \R} \int_{V}  (f(x)-a)^2 \mu(dx).
\end{align*}
Since $\inf_{a \in \R} \int_{V}  (f(x)-a)^2 \mu(dx)= \int_{V}  (f(x)-\overline f_{V})^2 \mu(dx)$, this proves \PI.

(ii) By \eqref{e:NDL-PI}, we can apply \NDL \ to get that  for all $x_0 \in M$ and $ r\in (0,R_0)$,
\begin{align*}
	&\P^{x_0} (\tau_{B(x_0,r)}  \ge  \phi(x_0,r/(q_2\eta_2))) \\
	& \ge \int_{B(x_0,\eta_1 r/(q_2\eta_2))} p^{B(x_0, r/q_0)} ( \phi(x_0,r/(q_2\eta_2)), x_0, y)\mu(dy) \ge \frac{c_1V(x_0,\eta_1 r/(q_2\eta_2))}{V(x_0, r/(q_2\eta_2))} \ge c_2.
\end{align*}
We used \VD \ in the last inequality above. By \eqref{e:phi-scale}, it follows that for all $x_0 \in M$ and $r \in (0,R_0)$,
$$
\E^{x_0}[\tau_{B(x_0,r)}] \ge \phi(x_0,r/(q_2\eta_2))\,	\P^{x_0} (\tau_{B(x_0,r)}  \ge  \phi(x_0,r/(q_2\eta_2))) \ge c_2 \phi(x_0,r/(q_2\eta_2)) \ge c_3 \phi(x_0,r).
$$

(iii) Suppose also that \RVD \ holds. By Corollary \ref{c:PI-FK} and Proposition \ref{p:FK}, since \PI \ holds by (i),  we see that \FK \ and \mEu \ hold.  Since \mEl \ holds by (ii), the proof is complete. \qed

\end{document}